\documentclass[11pt]{article}
\usepackage{graphicx}
\usepackage{xcolor}
\usepackage[scale=0.8]{geometry}
\usepackage{amsmath,amssymb,amsthm,graphicx,float,url}
\usepackage{mathtools}
\usepackage{upgreek}
\usepackage{makecell}
\usepackage{tikz,mathdots,yhmath,cancel,siunitx,array,multirow,amssymb,gensymb,tabularx,extarrows,booktabs}
\usetikzlibrary{fadings}
\usetikzlibrary{patterns}
\usetikzlibrary{shadows.blur}
\usetikzlibrary{shapes}
\usepackage{pgfplots,graphicx,float,svg}
\pgfplotsset{compat=1.18}
\usepackage{algorithm}
\usepackage{algpseudocode}
\usepackage{comment}
\usepackage{subcaption}

\usepackage[colorlinks=true]{hyperref}
\usepackage[nameinlink,capitalize,noabbrev]{cleveref}
\crefname{section}{Section}{Sections}
\crefname{subsection}{Section}{Sections}
\crefname{subsubsection}{Section}{Sections}

\usepackage{pdfsync}
\usepackage{float}
\usepackage{mathrsfs}
\title{\sf Volume-Preserving Geometric Shape Optimization of the Dirichlet Energy Using Variational Neural Networks}
\author{Amaury Bélières Frendo\footnote{ IRMA, Université de Strasbourg, CNRS UMR 7501, Inria, 7 rue René Descartes, 67084 Strasbourg, France. ({\tt amaury.belieres@math.unistra.fr})}
\and Emmanuel Franck\footnote{Inria, IRMA, Université de Strasbourg, CNRS UMR 7501, 7 rue René Descartes, 67084 Strasbourg, France. ({\tt emmanuel.franck@inria.fr})}
\and  Victor Michel-Dansac\footnote{Inria, IRMA, Université de Strasbourg, CNRS UMR 7501, 7 rue René Descartes, 67084 Strasbourg, France. ({\tt victor.michel-dansac@inria.fr})}
\and Yannick Privat\footnote{Universit\'e de Lorraine, CNRS, Institut Elie Cartan de Lorraine, Inria, BP 70239 54506
Vandœuvre-l\`es-Nancy Cedex, France. ({\tt yannick.privat@univ-lorraine.fr}).}~\footnote{Institut Universitaire de France (IUF)}
}

\date{\today}

\renewcommand{\leq}{\leqslant}

\newtheorem{theorem}{Theorem}

\newtheorem{definition}{Definition}
\newtheorem{lemma}{Lemma}

\theoremstyle{definition}\newtheorem{remark}{Remark}

\numberwithin{equation}{section}

\definecolor{review1}{RGB}{0,0,0}
\definecolor{review2}{RGB}{0,0,0}
\definecolor{review3}{RGB}{0,0,0}
\definecolor{review4}{RGB}{0,0,0}
\colorlet{allreviews}{review1!25!review2!50!review3!75!review4}

\newcommand{\ra}{\color{review1}}
\newcommand{\rb}{\color{review2}}
\newcommand{\rc}{\color{review3}}

\newcommand{\rall}{\color{allreviews}}

\begin{document}

\maketitle

\begin{abstract}
    In this work, we explore the numerical solution of geometric shape optimization problems using neural network-based approaches. This involves minimizing a numerical criterion that includes solving a partial differential equation with respect to a domain, often under geometric constraints like a constant volume. {\rc We successfully develop a proof of concept} using a flexible and parallelizable methodology to tackle these problems.
    We focus on a prototypal problem: minimizing the so-called Dirichlet energy with respect to the domain under a volume constraint, involving Poisson's equation in $\mathbb{R}^2$. We use variational neural networks
    to approximate the solution to Poisson's equation on a given domain, and represent the shape through a neural network that approximates a volume-preserving transformation from an initial shape to an optimal one. These processes are combined in a single optimization algorithm that minimizes the Dirichlet energy. {\rc A significant advantage of this approach is its inherent parallelizability}, which makes it easy to handle the addition of parameters. Additionally, it does not rely on shape derivative or adjoint calculations.
    Our approach is tested on Dirichlet and Robin boundary conditions, parametric right-hand sides, and extended to Bernoulli-type free boundary problems.
    The source code for solving the shape optimization problem is open-source and freely available.
\end{abstract}



\section{Introduction}

Throughout this article, let $V_0>0$ be a fixed real number and $n\in \mathbb{N}^*$. We will denote by $\mathbb{B}^n\subset \mathbb{R}^n$ the $n$-dimensional open ball with volume $V_0$ and by $\mathbb{S}^{n-1}$ its boundary. For simplicity, we will omit explicit notation of the dependency on $V_0$, as all the considered open sets share the same volume.

\subsection{Shape Optimization using Learning Techniques}
Shape optimization is a field of mathematics in which we seek to determine, if it exists, the shape of a domain that minimizes some numerical criterion. There are many such problems, ranging from the most fundamental to the most applied \cite{delfour2011shapes,henrot-pierre}. For example, on the one hand, optimizing the eigenvalues of the Laplacian or Schrödinger operators has a long history, dating back to at least as far as the work of Lord Rayleigh or Kac \cite{kac1966can}. The aim is to understand the very implicit links between the eigenmodes of these operators and the geometry. On the other hand, from a practical point of view, shape optimization is crucial for engineering applications. Let us mention, for instance, geometry optimization to increase the performance of a mechanical component, typically the optimal design of an aircraft wing \cite{allaire2007conception,mohammadi2009applied}.

In this article, we focus on the numerical implementation of shape optimization problems. There are many tools and algorithms available for determining (at least local) minimizers of such problems. Many approaches are based on the use of an ad-hoc notion of differential, either called shape derivative or derivative in the sense of Hadamard, enabling the implementation of gradient methods (see e.g. \cite{allaire2021shape}).

However, such numerical methods often suffer from several drawbacks.
Indeed, they are based on a derivative calculation that can become highly complex, for example when the criterion involves the solution of a multiphysics model. This approach often requires the determination of an adjoint problem and a descent step, which can make it very costly in terms of computation time and memory allocation.
Moreover, since shape optimization problems frequently have numerous local minimizers, these approaches tend to be highly local.
A major drawback of these methods is their general lack of parallelizability. This makes them particularly time-consuming and poorly suited for scenarios where model parameters need to be varied. In contrast, the method we introduce is highly parallelizable, offering a significant advantage in efficiency and adaptability.
Finally, such approaches are not suitable for all physical models. One example is models that are not ``well-posed'', such as the turbulent Navier-Stokes equations.

Our goal is to present a methodology that integrates a well-suited problem formulation with a parallel, mesh-free numerical method, such as a neural algorithm.
Let us mention \cite{odot2023real}, in which a first step is taken in this direction: the equation of state (hyper-elastic problem) is solved using the feed-forward neural network and the so-called {\it adjoint problem} is obtained through the backpropagation of the network.
{\ra %
Over the last few years, work has been undertaken on shape optimization
applied to mechanical engineering,
where the focus is placed on elliptic equations.
The main idea is to apply the SIMP (Solid Isotropic Material with Penalization) method
in its parameterized version (see e.g.~\cite{ChaSur2020,ZehLiCorTho2021,JEONG2023115484} or the review paper~\cite{WolAagBaeSig2022}),
replacing classical algorithms with physics-informed neural networks.
Moreover, in the same context, other authors parametrized level-set functions
with neural networks, see~\cite{ZhaYaoLiZhoChe2023}.
This parametric version of classical algorithms has undeniable computational and parallelizability advantages,
but the level of detail one can expect in the final shape
remains intrinsically hidden by the parameterization.
We can also mention \cite{CuoGiaIzzNitPicTro2022} that solves the Bernoulli free boundary problem with a PINNs, rephrasing it like an overdetermined PDE constrained by the measure of the positive part of its solution.
Note that this method is only applicable to such overdetermined problems,
and suffers from a lack of generalizability.
}

    {\rc We develop a proof of concept, highlighting the potential of this approach.
        It is validated} on a very simple model: Dirichlet energy minimization with homogeneous Dirichlet or Robin boundary conditions. The aim is to optimize, with respect to the domain, the ``natural'' energy associated with Poisson's equation, which takes the form
\begin{equation*}
    -\Updelta u=f\quad \text{in }\Upomega,
\end{equation*}
where $\Upomega \in \mathbb R^n$ is an open bounded connected set.
We add a classical volume constraint, which models a manufacturing cost: $|\Upomega|\leq V_0$, with $V_0>0$ a fixed parameter.
We deliberately focus on this prototypical problem, which has been the subject of much work and is now very well understood.

Note that there are various types of shape optimization problems, including parametric, geometric, and topological problems. The latter allow for topological changes (such as the number of holes). In this article, we focus exclusively on {\it geometric optimization problems}, where the desired shape must be homeomorphic to a given initial shape, such as a ball.

To describe our strategy,
note that our shape optimization problem consists in
minimizing the energy associated to a PDE
while also optimizing the shape of the space domain.
To that end, we introduce two neural networks.
The first one, a physics-informed neural network
(PINN, see~\cite{raissi2019physics}),
approximates the PDE solution.
The second one, a symplectic neural network
(SympNet, see~\cite{JIN2020166}),
takes care of the domain.
Indeed, symplectic maps
are examples of volume-preserving diffeomorphisms.
Hence, they are suitable candidates for representing the optimal shape,
which will be defined as the image, by the SympNet,
of a given initial domain.
As a consequence, in our strategy,
the constant volume constraint is automatically satisfied,
and it does not need to be added to the loss function.
Moreover, regarding the PDE approximation,
boundary conditions will be directly imposed in the network
rather than in the loss function,
in the spirit of~\cite{LagLikFot1998}.
Both of these remarks mean that we will only have to solve
a single optimization problem on the trainable parameters
of the PINN and of the SympNet, with a single loss function,
corresponding to the Dirichlet energy of the problem.

Conventional methods necessitate running a new simulation for each set of physical parameters or source terms. In contrast, a neural network can learn a parametric family of optimal shapes for a corresponding family of physical parameters or source terms, all within a time frame comparable to solving a single non-parametric problem.

Alongside this paper, we provide a turnkey open source code to solve such shape optimization problems.
The code comes as a ready-to-use \texttt{Python} library that the user can freely download on GitHub,
with a detailed documentation.
The code architecture is inspired by the AvaFrame python framework
\cite{tonnel2023avaframe}. The code is divided into three computational modules: \texttt{com1PINNs}, for PDE resolution with PINNs, \texttt{com2SympNets} for learning a given shape with SympNets, and \texttt{com3DeepShape} for learning-based
shape optimization.
To run each computational module, the user will find appropriate run scripts in the folder \texttt{examples}. In each run script, the user will find parameters that can be freely modified: number of layers and neurons in the network, shape of the initial domain, boundary conditions, source term, \dots

More information on the code can be found in \cref{sec:codeSource}.
The paper is then organized as follows.
First, \cref{sec:problem} introduces the shape optimization problem and reviews some theoretical results. Next, \cref{sec:nn} presents the neural networks and the joint optimization algorithm based on the use of appropriate symplectomorphisms. \cref{sec:numerics} focuses on the numerical validation of the code with Dirichlet boundary conditions, while \cref{sec:robin} addresses Robin boundary conditions.


\subsection{Main contributions and source code}\label{sec:codeSource}

We present an approach based on dual minimization: the first ensures the recovery of the manipulated PDE solution, while the second guarantees the minimization of the shape functional. It is important to note that, although we use neural networks for the minimization procedures, our approach does not involve traditional, data-driven learning. The PDE solution is obtained by minimizing the natural energy associated with the variational formulation. The optimal domain is sought among topological balls, i.e., the image of ball with volume~$V_0$ under volume-preserving transformations known as symplectomorphisms.

The resulting algorithm is also easily parallelizable. Another significant aspect of our approach, which makes it suitable for problems in higher spatial dimensions, is its mesh-free nature. It relies on collocation points and uses neural networks to evaluate and differentiate the quantities with respect to the parameters. Notably, the ``PDE resolution'' part, currently achieved at convergence by minimizing the energy associated with the variational formulation, could be replaced by any other mesh-free solving method. We have chosen this approach due to its potential for generalization to more complex scenarios and problems in higher spatial dimensions.

The source code \cite{amaury_belieres_frendo_2024_13133269} for solving the shape optimization problem is open and available at the following link:
\begin{center}
    \url{https://github.com/belieresfrendo/GeSONN}.
\end{center}
A comprehensive documentation explaining how to use it is also provided.

\subsection{Presentation of the shape optimization problem}
\label{sec:problem}

We first introduce a problem involving a PDE with Dirichlet conditions, but our approach generalizes to other conditions. For instance, the last section of this manuscript is dedicated to Robin boundary conditions.

Let $\Upomega$ be an open bounded connected set in $\mathbb R^n$ and $f\in H^{-1}({\Upomega})$.
In the whole article, we will denote by $u_\Upomega^f$ the unique solution in $H^1_0(\Upomega)$ of the Poisson problem
\begin{equation}
    \label{eq:poisson}
    \begin{cases}
        -\Updelta u^f_\Upomega = f & \text{in }\Upomega,         \\
        u^f_\Upomega=0             & \text{on }\partial\Upomega.
    \end{cases}
\end{equation}
For simplicity, we will only deal with homogeneous Dirichlet boundary conditions on $\partial\Upomega$.
It is well-known that $u_\Upomega^f$ can also be defined as the unique solution of the variational
problem, find $u^f_\Upomega \in H^1_0(\Upomega)$, such that
\begin{equation}
    \label{fv}
    \mathcal{J}(\Upomega, u^f_\Upomega) =
    \inf_{} \lbrace\mathcal{J}(\Upomega, u), \, u\in H^1_0(\Upomega)\rbrace,
\end{equation}
with
\begin{equation}\label{eq:J}
    \mathcal{J} (\Upomega, u) =
    \frac12\int_\Upomega \lvert\nabla u\rvert^2
    - \langle f,u\rangle_{H^{-1}(\Upomega),H^1_0(\Upomega)}, \quad  \forall u \in H^1_0(\Upomega),
\end{equation}
where $\langle\cdot,\cdot\rangle_{H^{-1}(\Upomega),H^1_0(\Upomega)}$ denotes the duality bracket between $H^{-1}(\Upomega)$ and $H^1_0(\Upomega)$.
Recall that the minimization problem above can be interpreted as an energetic formulation of the PDE \eqref{eq:poisson}.
Furthermore, if~$f$ enjoys additional regularity (say $f\in L^2(\Upomega)$), the PDE \eqref{eq:poisson} is satisfied at least in a pointwise manner.
Let us now introduce the so-called Dirichlet energy $\mathcal{E}$, a shape functional we will deal
with throughout this article. It is given by
\begin{equation}
    \label{eq:energy}
    \mathcal E(\Upomega) \coloneqq \inf_{u\in H^1_0(\Upomega)}\mathcal{J}(\Upomega, u)
\end{equation}
and it is notable that $\mathscr{E}(\Upomega)=\mathscr{J}(\Upomega,u_\Upomega^f)$.
Minimizing the Dirichlet energy within sets of given volume is a prototypical problem in shape optimization. It reads:
\begin{equation}
    \label{eq:optim_ener}
    \boxed{	\inf \lbrace\mathcal E(\Upomega),
        \, \Upomega \text{ bounded open set of }\mathbb R^n, \text{ such that }
        \lvert \Upomega \rvert=V_0\rbrace}.
\end{equation}
The analysis of such a problem is a long story. We refer for instance to \cite{henrot-pierre} for a detailed review of results, related to the existence of an optimal shape, but also to geometric properties of minimizers.

\begin{remark}[Some comments about existence and regularity of optimal shapes]
    Existence issues in shape optimization are generally difficult. It is standard to introduce a relaxed formulation of the initial problem posed among open sets. This often leads to considering a set of admissible forms living in the set of quasi-open sets. Recall that a set $A$ is said to be quasi-open if there exist sets $\omega$ of capacity as small as desired, such that $A \setminus \omega$ can be extended into an open by adding points of $\omega$.

    In the case we are interested in, it is easy, by ad-hoc extension of the definition of the functional space $H^1_0(\Upomega)$ to quasi-open $\Upomega$, to establish the existence of an optimal solution for the problem
    \begin{equation*}
        \inf_{} \lbrace\mathcal E(\Upomega),
        \, \Upomega \text{ quasi open set of } D, \text{ such that }
        \lvert \Upomega \rvert=V_0\rbrace,
    \end{equation*}
    where $D$ denotes a given compact set of $\mathbb R^n$.
    We refer to \cite[Chapter 4]{henrot-pierre} for detailed explanations about such notions.

    It is worth noting that in some cases (e.g. when considering shape optimization problems involving Dirichlet eigenvalues), the box constraint on $D$ can be removed, using concentration-compactness arguments originally established by Lions.
    Once such an existence result has been established, an a posteriori analysis can sometimes be carried out to establish the existence of regular openings. We refer, for example, to \cite[Chapters~2 and 3]{antunesshape}.

    In this article, the question of existence is not central. From the above-mentioned references, it is well-known that, even if we consider families of domains included in a fixed compact set $D$, there exists an optimal shape among the quasi-open sets. In the following, we will no longer mention this question, as our main problem here concerns the numerical determination of solutions to this problem, using techniques based on the use of neural networks.
\end{remark}

Another very useful result is the following characterization of optimal shapes.
\begin{theorem}(\cite[Section~6.1.3]{henrot-pierre})\label{thm:caracterisation_formes_optimales}
    Let $f\in L^2_\text{loc}(\mathbb R^n)$.
    If $\Upomega$ is a solution to \eqref{eq:optim_ener} with a $C^2$ boundary, then there exists $c>0$
    such that
    \begin{equation}\label{cion:surdet}
        \lvert \nabla u_\Upomega^f \rvert = c \;\; \text{on} \;\; \partial \Upomega.
    \end{equation}
\end{theorem}

It is noteworthy that this result provides another angle of attack for shape optimization problems. Indeed, rather than minimizing the shape functional $\mathcal E$ given by~\eqref{eq:energy} with respect to shapes $\Upomega$ of prescribed volume, we can seek to numerically solve the overdetermined PDE :
\begin{equation}\label{overdet:PDE}
    \begin{cases}
        -\Updelta u = f                   & \text{in }\Upomega,         \\
        \hfill u=0                        & \text{on }\partial\Upomega, \\
        \hfill \lvert \nabla u \rvert = c & \text{on }\partial\Upomega,
    \end{cases}
\end{equation}
where $c>0$ is a given parameter implicitly encoding the volume constraint.

In what follows, we will use the optimality condition provided by \cref{thm:caracterisation_formes_optimales} as a criterion for validating the shapes obtained. Indeed, if the algorithm we present leads to a domain $\Upomega_0$, we will consider it a good candidate for solving Problem~\eqref{eq:optim_ener} whenever Condition~\eqref{cion:surdet} is satisfied.

The overdetermined equation \eqref{overdet:PDE} will also serve as the starting point for solving Bernoulli problems, which are addressed in \cref{sec:bernoulli}.

\section{Solving the shape optimization problem with Variational Neural Networks} \label{sec:nn}

As outlined in \cref{sec:problem}, we will focus here on the Poisson problem \eqref{eq:poisson}. To identify areas for improvement, we briefly describe the main loop of classical shape optimization algorithms.

Initially, several PDEs are numerically solved to compute the state, the adjoint state if necessary, the shape gradient, and the gradient-descent step. Next, the computational domain is deformed according to the direction of the shape gradient. This is done through meshing, which requires a procedure known as {\it extension-regularization}, whose implementation is complex. This process is iterated until convergence, which can take hours, days, or even weeks for industrial applications \cite{allaire2002level,allaire2007conception}. Essentially, this method is not parallelizable. However, these challenges are not unavoidable and depend on the problem's formulation.


In what follows, we tackle these issues using neural networks, which present several advantages over traditional numerical methods for PDEs. For instance, by combining neural networks with the Monte-Carlo algorithm \cite{caflisch1998monte} for approximating integrals in the loss function, we can efficiently handle parameter-dependent problems on complex domains. Moreover, {\rb neural networks} allow for joint gradient descent on multiple interdependent networks. The shape derivative computation is effectively hidden behind an automatic differentiation process of the network. This is significant because calculating the shape derivative is complex and sometimes impossible for certain multiphysics problems. This advancement opens up new prospects for interdisciplinary research, especially in the life sciences. This means we can simultaneously train one network to represent the PDE solution and another to represent the computational domain, resulting in a parallelizable algorithm for solving shape optimization problems.

We now present a methodology for developing a neural network algorithm to solve the shape optimization problem \eqref{eq:optim_ener}. This approach leverages a variational Neural Network called DeepRitz \cite{e2017deep} to represent the solution to the Poisson problem \eqref{eq:poisson}, and SympNets \cite{JIN2020166} to represent the computational domain.

\subsection{DeepRitz, a subclass of PINNs methods}
\label{sec:pinns}

Let us outline the architecture of the neural network representing the solution to the PDE within the computational domain. We will first describe a fully connected neural network and then explain how to adapt it for solving PDEs, specifically the Poisson problem \eqref{eq:poisson}. For a comprehensive introduction to fully connected neural networks, PINNs, and DeepRitz, see \cite{Goodfellow-et-al-2016,e2017deep,raissi2019physics}.
A distinctive feature of these networks is the strictly physical nature of their loss function; no data is used for training the networks.



    {\rb
        In this work, we employ a fully connected neural network to represent the solution to the PDE. The neural networks are tailored to minimize the Dirichlet energy \eqref{eq:energy}, which is the natural energy associated with the variational formulation. This approach offers a significant advantage over classical Physics-Informed Neural Networks (PINNs \cite{raissi2019physics}) because it handles the PDE in its inherent weak form, capturing the problem's true nature.

        \begin{remark}
            The loss function $\mathcal{J}$ defined in \eqref{eq:J}, coupled to a fully connected neural network, is a particular case of a DeepRitz network introduced in~\cite{e2017deep}. If we had chosen the residual of the PDE as the loss function, it would have been a PINN, as defined in~\cite{raissi2019physics}.
            Nevertheless, for simplicity, we will call our {\rb neural network} representing the solution of the Poisson problem \eqref{eq:poisson} a PINN, because of its essential physics-informed character.
            We chose DeepRitz rather than classical PINNs for two reasons.
            First, the minimization of the Dirichlet energy problem only makes sense from a variational point of view.
            Indeed, the formulation \eqref{eq:optim_ener}, which includes two infima, prohibits the use of classical PINNs,
            since it would require taking the state equation as a constraint.
            This would serve no purpose other than to complicate the initial, simple problem.
            Second, classical PINNs require the computation of second derivatives,
            which is time-consuming and could preferably be avoided in this context.
        \end{remark}
    }


\subsubsection{\texorpdfstring{Fully connected {\rb neural networks}}{Fully connected neural networks}}
\label{sub:fcnn}

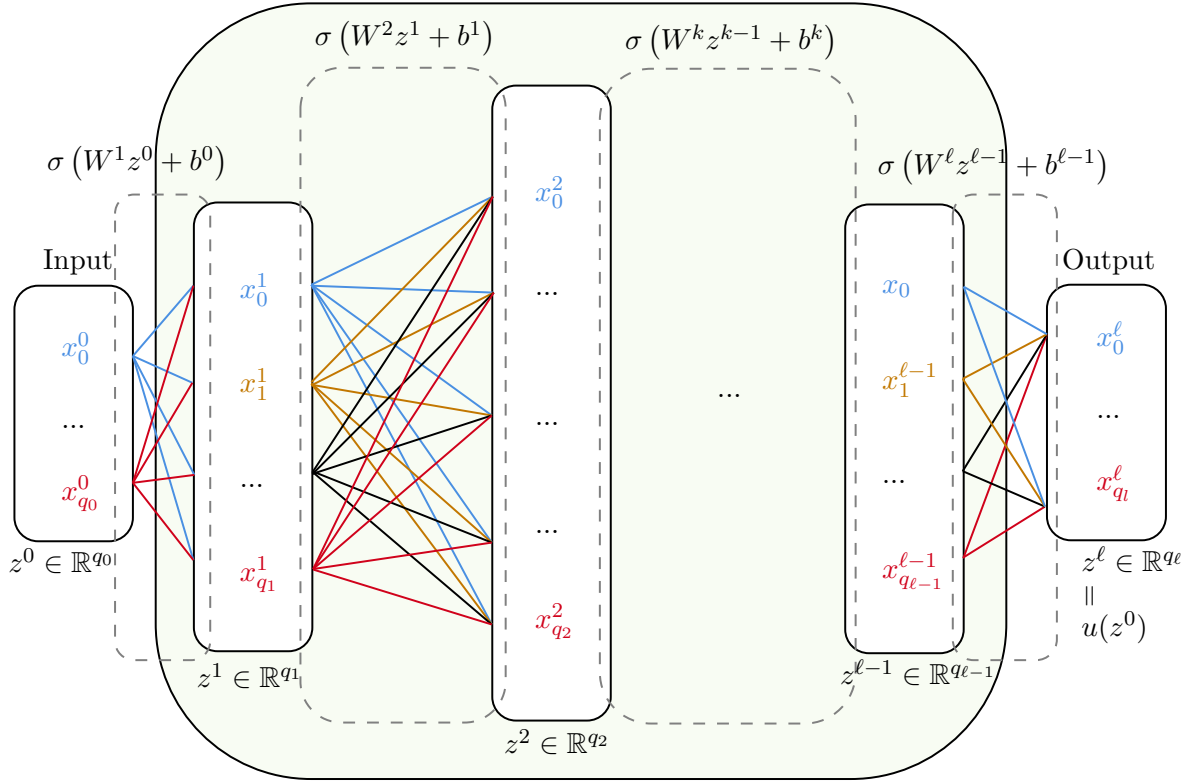
\begin{figure}[!ht]
    \begin{center}

        \tikzset{every picture/.style={line width=0.75pt}} 

        \begin{tikzpicture}[x=0.75pt,y=0.75pt,yscale=-1,xscale=1]
        
        \draw  [fill={rgb, 255:red, 184; green, 233; blue, 134 }  ,fill opacity=0.1 ] (111.25,80.8) .. controls (111.25,37.56) and (146.31,2.5) .. (189.55,2.5) -- (461.95,2.5) .. controls (505.19,2.5) and (540.25,37.56) .. (540.25,80.8) -- (540.25,315.7) .. controls (540.25,358.94) and (505.19,394) .. (461.95,394) -- (189.55,394) .. controls (146.31,394) and (111.25,358.94) .. (111.25,315.7) -- cycle ;
        \draw  [fill={rgb, 255:red, 255; green, 255; blue, 255 }  ,fill opacity=1 ] (130.5,114.95) .. controls (130.5,108.35) and (135.85,103) .. (142.45,103) -- (178.3,103) .. controls (184.9,103) and (190.25,108.35) .. (190.25,114.95) -- (190.25,317.55) .. controls (190.25,324.15) and (184.9,329.5) .. (178.3,329.5) -- (142.45,329.5) .. controls (135.85,329.5) and (130.5,324.15) .. (130.5,317.55) -- cycle ;
        \draw   (40,156.95) .. controls (40,150.35) and (45.35,145) .. (51.95,145) -- (87.8,145) .. controls (94.4,145) and (99.75,150.35) .. (99.75,156.95) -- (99.75,262.05) .. controls (99.75,268.65) and (94.4,274) .. (87.8,274) -- (51.95,274) .. controls (45.35,274) and (40,268.65) .. (40,262.05) -- cycle ;
        \draw  [fill={rgb, 255:red, 255; green, 255; blue, 255 }  ,fill opacity=1 ] (560.75,156.5) .. controls (560.75,149.87) and (566.12,144.5) .. (572.75,144.5) -- (608.75,144.5) .. controls (615.38,144.5) and (620.75,149.87) .. (620.75,156.5) -- (620.75,261.5) .. controls (620.75,268.13) and (615.38,273.5) .. (608.75,273.5) -- (572.75,273.5) .. controls (566.12,273.5) and (560.75,268.13) .. (560.75,261.5) -- cycle ;
        \draw  [fill={rgb, 255:red, 255; green, 255; blue, 255 }  ,fill opacity=1 ] (459,115.95) .. controls (459,109.35) and (464.35,104) .. (470.95,104) -- (506.8,104) .. controls (513.4,104) and (518.75,109.35) .. (518.75,115.95) -- (518.75,318.55) .. controls (518.75,325.15) and (513.4,330.5) .. (506.8,330.5) -- (470.95,330.5) .. controls (464.35,330.5) and (459,325.15) .. (459,318.55) -- cycle ;
        \draw  [fill={rgb, 255:red, 255; green, 255; blue, 255 }  ,fill opacity=1 ] (281,55.95) .. controls (281,49.35) and (286.35,44) .. (292.95,44) -- (328.8,44) .. controls (335.4,44) and (340.75,49.35) .. (340.75,55.95) -- (340.75,352.55) .. controls (340.75,359.15) and (335.4,364.5) .. (328.8,364.5) -- (292.95,364.5) .. controls (286.35,364.5) and (281,359.15) .. (281,352.55) -- cycle ;
        \draw [color={rgb, 255:red, 74; green, 144; blue, 226 }  ,draw opacity=1 ]   (99.75,180.5) -- (130.25,145) ;
        \draw [color={rgb, 255:red, 74; green, 144; blue, 226 }  ,draw opacity=1 ]   (99.75,180.5) -- (130.25,283.5) ;
        \draw [color={rgb, 255:red, 74; green, 144; blue, 226 }  ,draw opacity=1 ]   (99.75,180.5) -- (129.75,194) ;
        \draw [color={rgb, 255:red, 74; green, 144; blue, 226 }  ,draw opacity=1 ]   (99.75,180.5) -- (130.75,240.5) ;
        \draw [color={rgb, 255:red, 208; green, 2; blue, 27 }  ,draw opacity=1 ][fill={rgb, 255:red, 74; green, 144; blue, 226 }  ,fill opacity=1 ]   (99.75,244.5) -- (130.25,145) ;
        \draw [color={rgb, 255:red, 208; green, 2; blue, 27 }  ,draw opacity=1 ][fill={rgb, 255:red, 74; green, 144; blue, 226 }  ,fill opacity=1 ]   (99.75,244.5) -- (130.25,283.5) ;
        \draw [color={rgb, 255:red, 208; green, 2; blue, 27 }  ,draw opacity=1 ][fill={rgb, 255:red, 74; green, 144; blue, 226 }  ,fill opacity=1 ]   (99.75,244.5) -- (129.75,194) ;
        \draw [color={rgb, 255:red, 208; green, 2; blue, 27 }  ,draw opacity=1 ][fill={rgb, 255:red, 74; green, 144; blue, 226 }  ,fill opacity=1 ]   (99.75,244.5) -- (130.75,240.5) ;
        \draw [color={rgb, 255:red, 74; green, 144; blue, 226 }  ,draw opacity=1 ]   (189.75,144.5) -- (281.25,100) ;
        \draw [color={rgb, 255:red, 74; green, 144; blue, 226 }  ,draw opacity=1 ]   (189.75,144.5) -- (280.75,274.66) ;
        \draw [color={rgb, 255:red, 74; green, 144; blue, 226 }  ,draw opacity=1 ]   (189.75,145) -- (281.75,148.49) ;
        \draw [color={rgb, 255:red, 74; green, 144; blue, 226 }  ,draw opacity=1 ]   (189.75,144.5) -- (280.75,210.5) ;
        \draw [color={rgb, 255:red, 74; green, 144; blue, 226 }  ,draw opacity=1 ]   (189.75,144.5) -- (280.75,316) ;
        \draw [color={rgb, 255:red, 194; green, 119; blue, 0 }  ,draw opacity=1 ]   (189.75,193.75) -- (281.25,100) ;
        \draw [color={rgb, 255:red, 194; green, 119; blue, 0 }  ,draw opacity=1 ]   (189.75,194.3) -- (280.75,274.66) ;
        \draw [color={rgb, 255:red, 194; green, 119; blue, 0 }  ,draw opacity=1 ]   (189.75,195.32) -- (281.75,148.49) ;
        \draw [color={rgb, 255:red, 194; green, 119; blue, 0 }  ,draw opacity=1 ]   (189.75,195) -- (280.75,210.5) ;
        \draw [color={rgb, 255:red, 194; green, 119; blue, 0 }  ,draw opacity=1 ]   (189.75,194.42) -- (280.75,316) ;
        \draw [color={rgb, 255:red, 0; green, 0; blue, 0 }  ,draw opacity=1 ]   (190.21,239.35) -- (281.25,100) ;
        \draw [color={rgb, 255:red, 0; green, 0; blue, 0 }  ,draw opacity=1 ]   (190.21,238.67) -- (280.75,274.66) ;
        \draw [color={rgb, 255:red, 0; green, 0; blue, 0 }  ,draw opacity=1 ]   (190.21,240.17) -- (281.75,148.49) ;
        \draw [color={rgb, 255:red, 0; green, 0; blue, 0 }  ,draw opacity=1 ]   (190.21,240) -- (280.75,210.5) ;
        \draw [color={rgb, 255:red, 0; green, 0; blue, 0 }  ,draw opacity=1 ]   (190.21,238.9) -- (280.75,316) ;
        \draw [color={rgb, 255:red, 208; green, 2; blue, 27 }  ,draw opacity=1 ]   (190.67,288) -- (281.25,100) ;
        \draw [color={rgb, 255:red, 208; green, 2; blue, 27 }  ,draw opacity=1 ]   (190.67,288) -- (280.75,274.66) ;
        \draw [color={rgb, 255:red, 208; green, 2; blue, 27 }  ,draw opacity=1 ]   (190.67,288) -- (281.75,148.49) ;
        \draw [color={rgb, 255:red, 208; green, 2; blue, 27 }  ,draw opacity=1 ]   (190.66,288) -- (280.75,210.5) ;
        \draw [color={rgb, 255:red, 208; green, 2; blue, 27 }  ,draw opacity=1 ]   (190.66,288.26) -- (280.75,316) ;
        \draw [color={rgb, 255:red, 208; green, 2; blue, 27 }  ,draw opacity=1 ]   (518.25,282.39) -- (559.75,256.5) ;
        \draw [color={rgb, 255:red, 208; green, 2; blue, 27 }  ,draw opacity=1 ]   (518.75,282) -- (560.75,169.5) ;
        \draw [color={rgb, 255:red, 0; green, 0; blue, 0 }  ,draw opacity=1 ]   (518.25,238.61) -- (560.75,169.5) ;
        \draw [color={rgb, 255:red, 0; green, 0; blue, 0 }  ,draw opacity=1 ]   (518.25,238.61) -- (559.75,256.5) ;
        \draw [color={rgb, 255:red, 194; green, 119; blue, 0 }  ,draw opacity=1 ]   (518.25,192) -- (560.75,169.5) ;
        \draw [color={rgb, 255:red, 194; green, 119; blue, 0 }  ,draw opacity=1 ]   (518.25,192) -- (559.75,256.5) ;
        \draw [color={rgb, 255:red, 74; green, 144; blue, 226 }  ,draw opacity=1 ]   (518.75,145.5) -- (560.75,169.5) ;
        \draw [color={rgb, 255:red, 74; green, 144; blue, 226 }  ,draw opacity=1 ]   (518.75,145.5) -- (559.75,256.5) ;
        \draw  [color={rgb, 255:red, 128; green, 128; blue, 128 }  ,draw opacity=1 ][dash pattern={on 4.5pt off 4.5pt}] (184.25,55.7) .. controls (184.25,44.27) and (193.52,35) .. (204.95,35) -- (267.05,35) .. controls (278.48,35) and (287.75,44.27) .. (287.75,55.7) -- (287.75,344.8) .. controls (287.75,356.23) and (278.48,365.5) .. (267.05,365.5) -- (204.95,365.5) .. controls (193.52,365.5) and (184.25,356.23) .. (184.25,344.8) -- cycle ;
        \draw  [color={rgb, 255:red, 128; green, 128; blue, 128 }  ,draw opacity=1 ][dash pattern={on 4.5pt off 4.5pt}] (90.75,108.6) .. controls (90.75,103.3) and (95.05,99) .. (100.35,99) -- (129.15,99) .. controls (134.45,99) and (138.75,103.3) .. (138.75,108.6) -- (138.75,324.4) .. controls (138.75,329.7) and (134.45,334) .. (129.15,334) -- (100.35,334) .. controls (95.05,334) and (90.75,329.7) .. (90.75,324.4) -- cycle ;
        \draw  [color={rgb, 255:red, 128; green, 128; blue, 128 }  ,draw opacity=1 ][dash pattern={on 4.5pt off 4.5pt}] (513.25,109.5) .. controls (513.25,103.7) and (517.95,99) .. (523.75,99) -- (555.25,99) .. controls (561.05,99) and (565.75,103.7) .. (565.75,109.5) -- (565.75,323.5) .. controls (565.75,329.3) and (561.05,334) .. (555.25,334) -- (523.75,334) .. controls (517.95,334) and (513.25,329.3) .. (513.25,323.5) -- cycle ;
        \draw  [color={rgb, 255:red, 128; green, 128; blue, 128 }  ,draw opacity=1 ][dash pattern={on 4.5pt off 4.5pt}] (335.25,61.3) .. controls (335.25,47.05) and (346.8,35.5) .. (361.05,35.5) -- (438.45,35.5) .. controls (452.7,35.5) and (464.25,47.05) .. (464.25,61.3) -- (464.25,340.2) .. controls (464.25,354.45) and (452.7,366) .. (438.45,366) -- (361.05,366) .. controls (346.8,366) and (335.25,354.45) .. (335.25,340.2) -- cycle ;
        
        \draw (53,125) node [anchor=north west][inner sep=0.75pt]  [xscale=1,yscale=1] [align=left] {Input};
        \draw (567,125) node [anchor=north west][inner sep=0.75pt]  [xscale=1,yscale=1] [align=left] {Output};
        \draw (55.5,166.4) node [anchor=north west][inner sep=0.75pt]  [xscale=1,yscale=1]  {$ \begin{array}{l}
        \textcolor[rgb]{0.29,0.56,0.89}{x_{0}^{0}}\\
        \\
        ...\\
        \\
        \textcolor[rgb]{0.82,0.01,0.11}{x_{q_0}^{0}}
        \end{array}$};
        \draw (577.5,161.4) node [anchor=north west][inner sep=0.75pt]  [xscale=1,yscale=1]  {$ \begin{array}{l}
        \textcolor[rgb]{0.29,0.56,0.89}{x_{0}^{\ell}}\\
        \\
        ...\\
        \\
        \textcolor[rgb]{0.82,0.01,0.11}{x_{q_l}^{\ell}}
        \end{array}$};
        \draw (145.5,136.9) node [anchor=north west][inner sep=0.75pt]  [xscale=1,yscale=1]  {$ \begin{array}{l}
        \textcolor[rgb]{0.29,0.56,0.89}{x_{0}^{1}}\\
        \\[0.3 cm]
        \textcolor[rgb]{0.76,0.47,0}{x_{1}^{1}}\\
        \\[0.3 cm]
        ...\\
        \\[0.3 cm]
        \textcolor[rgb]{0.82,0.01,0.11}{x_{q_1}^{1}}
        \end{array}$};
        \draw (469.5,134.4) node [anchor=north west][inner sep=0.75pt]  [xscale=1,yscale=1]  {$ \begin{array}{l}
        \textcolor[rgb]{0.29,0.56,0.89}{x_{0}}\\
        \\[0.3 cm]
        \textcolor[rgb]{0.76,0.47,0}{x_{1}^{\ell-1}}\\
        \\[0.3 cm]
        ...\\
        \\[0.3 cm]
        \textcolor[rgb]{0.82,0.01,0.11}{x_{q_{\ell-1}}^{\ell-1}}
        \end{array}$};
        \draw (294.5,87.9) node [anchor=north west][inner sep=0.75pt]  [xscale=1,yscale=1]  {$ \begin{array}{l}
        \textcolor[rgb]{0.29,0.56,0.89}{x_{0}^{2}}\\
        \\[0.3 cm]
        ...\\
        \\
        \\[0.3 cm]
        ...\\
        \\
        \\
        ...\\[0.3 cm]
        \\
        \textcolor[rgb]{0.82,0.01,0.11}{x_{q_{2}}^{2}}
        \end{array}$};
        \draw (35.5,275.9) node [anchor=north west][inner sep=0.75pt]  [xscale=1,yscale=1]  {$z^{0}\in\mathbb{R}^{q_0}$};
        \draw (570,274.9) node [anchor=north west][inner sep=0.75pt]  [xscale=1,yscale=1]  {\makecell{$\begin{array}{lll}
            z^{\ell}\in\mathbb{R}^{q_{\ell}} \\[0.em]
            \,\rotatebox{90}{=} \\[-0.25em]
            u(z^0)
        \end{array}$}};
        \draw (130.5,333.9) node [anchor=north west][inner sep=0.75pt]  [xscale=1,yscale=1]  {$z^{1}\in\mathbb{R}^{q_1}$};
        \draw (285.5,367.4) node [anchor=north west][inner sep=0.75pt]  [xscale=1,yscale=1]  {$z^{2}\in\mathbb{R}^{q_2}$};
        \draw (455,331.9) node [anchor=north west][inner sep=0.75pt]  [xscale=1,yscale=1]  {$z^{\ell-1}\in\mathbb{R}^{q_{\ell-1}}$};
        \draw (55,72.9) node [anchor=north west][inner sep=0.75pt]  [xscale=1,yscale=1]  {$\sigma \left( W^{1} z^{0} +b^{0}\right)$};
        \draw (474,74.4) node [anchor=north west][inner sep=0.75pt]  [xscale=1,yscale=1]  {$\sigma \left( W^{\ell} z^{\ell-1} +b^{\ell-1}\right)$};
        \draw (190,9.9) node [anchor=north west][inner sep=0.75pt]  [xscale=1,yscale=1]  {$\sigma \left( W^{2} z^{1} +b^{1}\right)$};
        \draw (393,198) node [anchor=north west][inner sep=0.75pt]  [xscale=1,yscale=1] [align=left] {...};
        \draw (346.45,11.4) node [anchor=north west][inner sep=0.75pt]  [xscale=1,yscale=1]  {$\sigma \left( W^{k} z^{k-1} +b^{k}\right)$};

        \end{tikzpicture}
        \caption{%
            Diagram of a fully connected {\rb neural network},
            where the output of each layer is a vector $z_k \in \mathbb{R}^{q_k}$,
            with $q_k$ the number of neurons of the layer $k$. Computing $z_k$ involves a weight matrix $W^k \in \mathcal M_{q_k, q_{k-1}}(\mathbb R)$, a bias vector $b^k\in \mathbb{R}^{q_k}$
            and a nonlinear activation function $\sigma:\mathbb{R}\to\mathbb{R}$
            applied componentwise.
        }
        \label{fig:figure_pinns}
    \end{center}
\end{figure}

Consider again the problem \eqref{eq:poisson}, defined on an open set $\Upomega$. The main idea behind fully-connected neural networks is to represent the solution $u\in H^1_0(\Upomega)$ of the forward problem \eqref{eq:poisson} by $u_\theta\in C^\infty(\mathbb{R}^n)$, a composition of nonlinear parametric functions that takes $x\in\Upomega$ as input and returns an approximation of $u(x)$ as output, see \cite{raissi2019physics}. The network and associated notation are depicted in \cref{fig:figure_pinns}. The parameters $\theta = \lbrace W^k,b^k\rbrace_{k=1}^\ell$ of the {\rb neural network}, called trainable weights, are then optimized by minimizing a loss function.


To obtain a solution respecting the boundary conditions, we introduce two functions $\alpha, \beta \in C^1(\mathbb{R}^n, \mathbb{R})$ such that $\alpha$ vanishes on $\partial\Upomega$, and $\beta$ is equal, on $\partial\Upomega$, to the boundary condition of the Poisson problem \eqref{eq:poisson}.
%
%
%
Of course, for homogeneous Dirichlet boundary conditions, $\beta = 0$ is suitable. Then, for all $x \in \mathbb{R}^n$,
\begin{equation}\label{eq:bc_net}
    v_\theta(x) = \alpha(x) \, u_\theta(x) + \beta(x)
\end{equation}
satisfies the approximation of the PDE given by the PINN \cite{LagLikFot1998,franck2023approximately}.
For example, if~$v_\theta$ has to satisfy homogeneous Dirichlet conditions on $\mathbb S^1$, the unit sphere of $\mathbb R^2$, it is sufficient to define, for $x = (x_1,x_2)\in\mathbb B^2$,
\begin{equation*}
    \begin{cases}
        \alpha(x)=1-x_1^2-x_2^2, \\
        \beta(x)=0.
    \end{cases}
\end{equation*}
This can of course be generalized for other domains and boundary conditions.

The network we train is still $u_\theta$, but the solution is represented by~$v_\theta$, which is forced by construction to satisfy the boundary conditions on $\partial \Upomega$. Imposed that way, the boundary conditions are easier to implement, and they are less reliant on parameter adjustment than if we had added a penalization term to the loss function.

\subsubsection{The Dirichlet energy, an appropriate loss function}

To obtain the optimal network parameters $\theta^*$, we minimize a loss function during the training process. Let us design this loss function. Given that the neural network will be trained without any data, the loss function must incorporate the Dirichlet energy. This is because minimizing the Dirichlet energy yields the solution {\rb of} the underlying PDE, as indicated by \eqref{fv}.

As the computational domain for the Poisson problem \eqref{eq:poisson} can have a complex topology, meshing it can be cumbersome and time-consuming. Therefore, the integral in the loss function $\mathcal{J}$ is approximated with the Monte-Carlo method~\cite{caflisch1998monte}, rather than with a classical quadrature method.
This leads us to defining $\mathcal{J}_N$, a discrete quadrature of $\mathcal{J}$, given by
\begin{equation} \label{eq:loss}
    \mathcal{J}_N\left(\theta; \lbrace x^i \rbrace_{i=1}^{N}\right) =
    \frac{V_0}{N} \sum_{i=1}^{N}
    \left \lbrace \frac12 \lvert \nabla v_\theta(x^i)\rvert^2
    - f(x^i)v_\theta(x^i) \right \rbrace,
\end{equation}
where $\lbrace x^i \rbrace_{i=1}^{N}\in\Upomega^N$ are $N$ collocation points on which to evaluate the loss function~\eqref{eq:loss}, and where~$u_\theta$ and~$v_\theta$ are defined in \cref{sub:fcnn}.
Recall from~\cite{caflisch1998monte} that the convergence speed of the Monte-Carlo algorithm does not depend on the dimension of the problem, but only on $N$.
Since $\smash{u^f_\Upomega}$ solving the Poisson problem~\eqref{eq:poisson} is the unique minimizer of~$\mathcal{J}$, it is expected that $v_\theta$ is a fair approximation of the minimizer of $\mathcal{J}_N$.
The gradient of the loss function \eqref{eq:loss} is subsequently computed using the PyTorch library~\cite{AnsYan2024}.

Starting with randomly initialized parameters $\theta$, we use a gradient descent method to find optimal parameters $\theta^*$, i.e., ones that minimize \eqref{eq:loss}. At iteration $k$, the parameters are updated with the Adam Optimizer \cite{kingma2014adam}, a stochastic gradient descent algorithm.
Once the training procedure is complete, the trained {\rb neural network} $v_{\theta^*}$ provides an approximation of the solution to the forward problem \eqref{eq:poisson}.

\subsubsection{Parametric problems}

Note that the source term of the Poisson problem can be selected as a parametric function. Upon completing the training procedure, and with a slightly increased but comparable computation time, the Poisson problem can be solved for various parameter values within a given set.

Let $n_\mu$ be the number of parameters in the problem. We denote by $\mathbb{M} \subset \mathbb{R}^{n_\mu}$ the parameter space. The parametric neural network then takes as input both a position $x$ within the computational domain $\Upomega$ and parameters $\mu$ from the parameter space $\mathbb{M}$. Let $f^p:\Upomega\times \mathbb{M} \to \mathbb{R}$ represent the parametric source term.

The solution $u^p:\Upomega\times \mathbb{M} \to \mathbb{R}$ of the parametric Poisson problem satisfies
\begin{equation}\label{eq:parametric_poisson}
    \begin{cases}
        -\Updelta u^p(x;\mu) = f^p(x;\mu), \, & \text{for} \, (x,\mu)\in \Upomega\times\mathbb{M};         \\
        u^p(x;\mu) = 0, \,                    & \text{for} \, (x,\mu)\in \partial\Upomega\times\mathbb{M}. \\
    \end{cases}
\end{equation}
As an example, the space of parameters could be $\mathbb{M}  = \mathbb R \times \mathbb{R}$, and $f^p$ could be given by
\begin{equation*}
    f^p:\big(x=(x_1,x_2);\mu = (\mu_1, \mu_2)\big)
    \mapsto
    \exp \left(1-\left(\frac{x_1}{\mu_1}\right)^2 - \left(\frac{x_2}{\mu_2}\right)^2\right).
\end{equation*}

At the end of the training procedure, the {\rb neural network} will be such that $(x,\mu)\in \Upomega \times \mathbb{M} \mapsto v_\theta(x;\mu)$ approximates the parametric Poisson problem \eqref{eq:parametric_poisson} a.e. in $\Upomega\times\mathbb{M}$.
The parametric loss function $\mathcal{J}_N^p$ then becomes
\begin{equation*}
    \label{eq:loss_parametric}
    \mathcal{J}_N^p\left(\theta; \lbrace x^i, \mu^i \rbrace_{i=1}^{N}\right) =
    \frac{V_0}{N} \sum_{i=1}^{N}
    \left \lbrace \frac12 \lvert \nabla v_\theta^p(x^i; \mu^i)\rvert^2
    - f^p(x^i;\mu^i) \, v_\theta^p(x^i; \mu^i) \right \rbrace,
\end{equation*}
where $\lbrace x^i, \mu^i \rbrace_{i=1}^{N}$ are $N$ collocation points in $\Upomega\times\mathbb{M}$.
As defined in \cref{sub:fcnn}, we set $v_\theta^p = u_\theta^p \alpha^p + \beta^p$,
where the functions
$\alpha^p, \beta^p \in C^1(\mathbb{R}^2\times\mathbb{M}, \mathbb{R})$
are such that $\alpha^p$ vanishes on $\partial\Upomega$,
and $\beta^p$ is equal, on $\partial\Upomega$,
to the boundary condition of the Poisson problem \eqref{eq:poisson}.
Note that $\mathcal{J}_N^p\big(\theta; \lbrace x^i, \mu^i \rbrace_{i=1}^{N}\big)$,
for $\lbrace x^i, \mu^i \rbrace_{i=1}^{N}\in \Upomega\times \mathbb{M}$,
is a quadrature of the limit functional $\mathcal{J}^p$,
defined for $v^p=u^p \alpha^p + \beta^p$,
with $u^p:\Upomega\times \mathbb{M}\mapsto\mathbb R$, by
\begin{equation*}
    \mathcal{J}^p(u^p) =
    \int_{\Upomega\times\mathbb{M}} \left(
    \frac{1}{2}|\nabla v^p(x;\mu)|^2
    -
    f^p(x;\mu)v^p (x;\mu)
    \right) \mathrm{dx} \, \mathrm{d}\upmu.
\end{equation*}

\subsection{SympNets}

In the previous section, we addressed solving a PDE in a given domain $\Upomega$. Recall that our goal is to minimize the Dirichlet energy with respect to both the solution $u$ and the domain $\Upomega$. This also means minimizing~$\mathcal{J}(\Upomega, u)$ with respect to 
$\Upomega$, among the open sets in $\mathbb R^n$ with volume $V_0$. To achieve this, we must parameterize the set of open subsets of $\mathbb R^n$. For practical and simplicity reasons, we assume the optimal set is connected and enjoys several regularity properties. This will lead us to introduce appropriate invertible differentiable transformations of $\mathbb S^{n-1}$.


\subsubsection{Using symplectic maps for shape optimization}

Our objective is to develop a neural network representation for a shape. To achieve this, we seek an optimal volume-preserving diffeomorphism that maps $\mathbb S^{n-1}$ to a shape that minimizes the shape criterion. By focusing on finding the optimal diffeomorphism, this method aligns with ``geometric optimization'' as it allows us to refine the shape's geometry while preserving its topological characteristics (in particular its genus).

Rather than considering generic diffeomorphisms, we focus on the special case of symplectic maps\footnote{For those unfamiliar with symplectic maps, they are used to model the evolution of dynamical systems while preserving certain geometric properties. For instance, the flow of a Hamiltonian system is a symplectic map, preserving volume in the phase space of the dynamical system.
    For more details on symplectic maps, the reader is referred to \cite{Arnold,nakahara2018geometry,hairer2006structure}. In a nutschell, in $\mathbb{R}^n$ with $n = 2d$, a $C^1$-diffeo\-morphism~$\mathcal S$ is called symplectic if its Jacobian matrix $D\mathcal{S}$ satisfies $(D\mathcal{S})^\intercal J (D\mathcal{S}) = J$, where $J$ is a block matrix called the standard symplectic form on $\mathbb R^{n}$, with $J_{11} = J_{22} = 0_d$ and $J_{12} = - J_{21} = I_d$. This definition is enough to see that $\lvert \det D\mathcal{S}\rvert = 1$, and therefore that $\mathcal{S}$ is volume-preserving.
}, which are $C^1$-diffeo\-morphisms which also preserve volume.

This property is particularly relevant in our context since the shape optimization problem involves a volume constraint. However, symplectic maps are only defined for even-dimensional spaces. This is consistent with our focus of optimizing shapes in $\mathbb R^2$. Note that, in the specific case of $\mathbb R^2$, the symplectic form coincides with the volume form. This implies that, in $\mathbb R^2$, a diffeomorphism is volume-preserving if and only if it is a symplectic map.


The following section is dedicated to presenting
the architecture of the {\rb neural network} that represents the shape of the domain and detailing the training process to approximate any symplectic map..

\subsubsection{\texorpdfstring{Symplectic maps and {\rb neural networks}}{Symplectic maps and neural networks}}
\label{sec:sympnets}

A SympNet is a {\rb neural network} with a symplectic structure capable of approximating any symplectic map. SympNets were first introduced in \cite{JIN2020166} to approximate the flows of Hamiltonian systems.

Our objective is now to train a SympNet to replicate any continuous and differentiable transformation $\mathcal{T}$ of $\mathbb{R}^{2d}$. This section does not cover shape optimization; instead, it focuses on defining SympNets. Specifically, we aim to explain how to train a SympNet to learn the transformation of the ($2d-1$)-dimensional sphere into a given shape in $\mathbb{R}^{2d}$. generated by a given symplectic map $\mathcal{T}$.
\cref{sec:shape_optim_with_sympnets} is dedicated to solving a PDE in a domain generated with a symplectic map.

To properly define SympNets, we first need to define shear maps,
which can be seen as building blocks of symplectic maps.

\begin{definition}[Shear maps \cite{Arnold}]
    One of the simplest families of symplectic transformations from $\mathbb{R}^{2d}$ into $\mathbb{R}^{2d}$ is called ``shear maps'', and is defined by
    \begin{equation*}
        f_\text{\rm up}
        \begin{pmatrix}
            x_1 \\
            x_2
        \end{pmatrix}
        =
        \begin{pmatrix}
            x_1 + \nabla V_\text{\rm up}(x_2) \\
            x_2
        \end{pmatrix}
        \quad\text{and}\quad
        f_\text{\rm down}
        \begin{pmatrix}
            x_1 \\
            x_2
        \end{pmatrix}
        =
        \begin{pmatrix}
            x_1 \\
            x_2 + \nabla V_\text{\rm down}(x_1)
        \end{pmatrix}
    \end{equation*}
    for $x = (x_1, x_2)^\intercal$ with $x_1 \in \mathbb{R}^d$ and $x_2 \in \mathbb{R}^d$, where $V_\text{\rm up/down} \in C^2(\mathbb{R}^d,\mathbb{R}^d)$, and $\nabla V_\text{\rm up/down}: \mathbb{R}^d \to \mathbb{R}^d$ is the gradient of $V_\text{\rm up/down}$.
\end{definition}

Using the previous definition,
the architecture of SympNets is based on the following lemma.
\begin{lemma}[\cite{JIN2020166}]
    Any symplectic map can be approximated by composing of several shear maps, and the composition of several symplectic maps still remains symplectic.
\end{lemma}

In practice, we could directly approximate the smooth functions $V_{\text{up}/\text{down}}$ using a fully connected {\rb neural network}. However, this approach would require the computation of the gradient of $V_{\text{up}/\text{down}}$ at each iteration. SympNets are designed to bypass this computational constraint, with the advantage of requiring very few trainable parameters compared to fully connected {\rb neural networks}. To this aim, \cite{JIN2020166} introduces gradient modules, the principles of which are summarized below.


Let $q>0$ be the depth of the {\rb neural network}.
In practice, we set $q>2d$.
We define $\hat{\sigma}_{K, a, b}$
the approximation of $\nabla V_{\text{up}/\text{down}}$ in terms of
an activation function~$\sigma: \mathbb{R} \to \mathbb{R}$,
two vectors $a, b \in \mathbb{R}^q$,
a matrix $K \in \mathcal M _{q,d}(\mathbb R)$,
and $\mathrm{diag}(a) = (a_i \delta_{ij})_{1 \leq i, j \leq q }$,
as follows, for all $x \in \mathbb{R}^d$,
\begin{equation*}
    \hat{\sigma}_{K, a, b}(x)
    =
    K^\intercal \mathrm{diag}(a) \sigma(K x + b),
\end{equation*}
where it is understood that $\sigma$ is applied to a vector componentwise.
Then, gradient modules $\mathcal G_\text{up}$ and~$\mathcal G_\text{down}$ are defined to approximate $f_\text{up}$ and $f_\text{down}$, by
\begin{equation*}
    \mathcal{G}_\text{up}
    \begin{pmatrix}
        x_1 \\
        x_2
    \end{pmatrix}
    =
    \begin{pmatrix}
        x_1 + \hat{\sigma}_{K, a, b}(x_2) \\
        x_2
    \end{pmatrix}
    \quad\text{and}\quad
    \mathcal{G}_\text{down}
    \begin{pmatrix}
        x_1 \\
        x_2
    \end{pmatrix}
    =
    \begin{pmatrix}
        x_1                               \\
        x_2 + \hat{\sigma}_{K, a, b}(x_1) \\
    \end{pmatrix}.
\end{equation*}
These functions are called gradient modules because $\hat{\sigma}_{K, a, b}$ is able to approximate any Jacobian of smooth vector fields, see~\cite[Appendix~A]{JIN2020166}.

Moreover, it is also possible to define a parametric gradient module to
approximate a family of symplectic maps
indexed by $n_\mu$ parameters
$\mu \in \mathbb M \subset \mathbb{R}^{n_\mu}$.
To that end, we define a second matrix
$K_\mu \in \mathcal M _{q,n_\mu}(\mathbb R)$,
and replace $\hat{\sigma}_{K, a, b}(x)$ with
\begin{equation*}
    \tilde{\sigma}_{K, K_\mu, a, b}(x;\mu)
    =
    K^\intercal \sigma(K x + b + K_\mu \mu).
\end{equation*}
One can show that $\tilde{\sigma}_{K, K_\mu, a, b}$ is a gradient module \cite{9716789},
and that the whole network remains symplectic
with respect to $x\in\mathbb{R}^{2d}$,
for each parameter $\mu \in \mathbb M$.

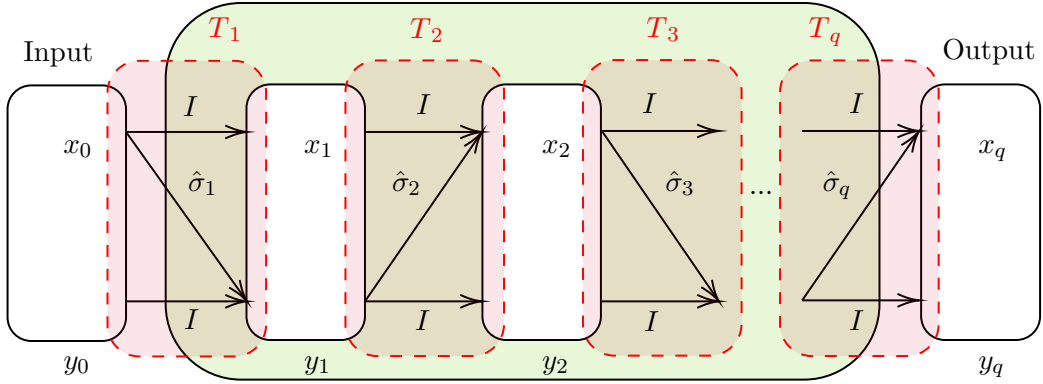
\begin{figure}[!ht]
    \begin{center}

        \tikzset{every picture/.style={line width=0.75pt}} 

        \begin{tikzpicture}[x=0.75pt,y=0.75pt,yscale=-1,xscale=1]
        
        \draw  [fill={rgb, 255:red, 184; green, 233; blue, 134 }  ,fill opacity=0.33 ] (149.75,41.42) .. controls (149.75,20.42) and (166.77,3.4) .. (187.77,3.4) -- (471.73,3.4) .. controls (492.73,3.4) and (509.75,20.42) .. (509.75,41.42) -- (509.75,155.48) .. controls (509.75,176.48) and (492.73,193.5) .. (471.73,193.5) -- (187.77,193.5) .. controls (166.77,193.5) and (149.75,176.48) .. (149.75,155.48) -- cycle ;
        \draw   (70,56.95) .. controls (70,50.35) and (75.35,45) .. (81.95,45) -- (117.8,45) .. controls (124.4,45) and (129.75,50.35) .. (129.75,56.95) -- (129.75,162.05) .. controls (129.75,168.65) and (124.4,174) .. (117.8,174) -- (81.95,174) .. controls (75.35,174) and (70,168.65) .. (70,162.05) -- cycle ;
        \draw  [fill={rgb, 255:red, 255; green, 255; blue, 255 }  ,fill opacity=1 ] (190.5,56.45) .. controls (190.5,49.85) and (195.85,44.5) .. (202.45,44.5) -- (238.3,44.5) .. controls (244.9,44.5) and (250.25,49.85) .. (250.25,56.45) -- (250.25,161.55) .. controls (250.25,168.15) and (244.9,173.5) .. (238.3,173.5) -- (202.45,173.5) .. controls (195.85,173.5) and (190.5,168.15) .. (190.5,161.55) -- cycle ;
        \draw  [fill={rgb, 255:red, 255; green, 255; blue, 255 }  ,fill opacity=1 ] (309.5,56.45) .. controls (309.5,49.85) and (314.85,44.5) .. (321.45,44.5) -- (357.3,44.5) .. controls (363.9,44.5) and (369.25,49.85) .. (369.25,56.45) -- (369.25,161.55) .. controls (369.25,168.15) and (363.9,173.5) .. (357.3,173.5) -- (321.45,173.5) .. controls (314.85,173.5) and (309.5,168.15) .. (309.5,161.55) -- cycle ;
        \draw  [fill={rgb, 255:red, 255; green, 255; blue, 255 }  ,fill opacity=1 ] (530.75,56.5) .. controls (530.75,49.87) and (536.12,44.5) .. (542.75,44.5) -- (578.75,44.5) .. controls (585.38,44.5) and (590.75,49.87) .. (590.75,56.5) -- (590.75,161.5) .. controls (590.75,168.13) and (585.38,173.5) .. (578.75,173.5) -- (542.75,173.5) .. controls (536.12,173.5) and (530.75,168.13) .. (530.75,161.5) -- cycle ;
        \draw    (129.75,68.5) -- (188.75,68.5) ;
        \draw [shift={(190.75,68.5)}, rotate = 180] [color={rgb, 255:red, 0; green, 0; blue, 0 }  ][line width=0.75]    (10.93,-3.29) .. controls (6.95,-1.4) and (3.31,-0.3) .. (0,0) .. controls (3.31,0.3) and (6.95,1.4) .. (10.93,3.29)   ;
        \draw    (129.75,154) -- (188.75,154) ;
        \draw [shift={(190.75,154)}, rotate = 180] [color={rgb, 255:red, 0; green, 0; blue, 0 }  ][line width=0.75]    (10.93,-3.29) .. controls (6.95,-1.4) and (3.31,-0.3) .. (0,0) .. controls (3.31,0.3) and (6.95,1.4) .. (10.93,3.29)   ;
        \draw    (250.25,154) -- (307.25,154) ;
        \draw [shift={(309.25,154)}, rotate = 180] [color={rgb, 255:red, 0; green, 0; blue, 0 }  ][line width=0.75]    (10.93,-3.29) .. controls (6.95,-1.4) and (3.31,-0.3) .. (0,0) .. controls (3.31,0.3) and (6.95,1.4) .. (10.93,3.29)   ;
        \draw    (250.25,68.5) -- (307.25,68.5) ;
        \draw [shift={(309.25,68.5)}, rotate = 180] [color={rgb, 255:red, 0; green, 0; blue, 0 }  ][line width=0.75]    (10.93,-3.29) .. controls (6.95,-1.4) and (3.31,-0.3) .. (0,0) .. controls (3.31,0.3) and (6.95,1.4) .. (10.93,3.29)   ;
        \draw    (470.75,68) -- (527.75,68) ;
        \draw [shift={(529.75,68)}, rotate = 180] [color={rgb, 255:red, 0; green, 0; blue, 0 }  ][line width=0.75]    (10.93,-3.29) .. controls (6.95,-1.4) and (3.31,-0.3) .. (0,0) .. controls (3.31,0.3) and (6.95,1.4) .. (10.93,3.29)   ;
        \draw    (470.75,153.5) -- (527.75,153.5) ;
        \draw [shift={(529.75,153.5)}, rotate = 180] [color={rgb, 255:red, 0; green, 0; blue, 0 }  ][line width=0.75]    (10.93,-3.29) .. controls (6.95,-1.4) and (3.31,-0.3) .. (0,0) .. controls (3.31,0.3) and (6.95,1.4) .. (10.93,3.29)   ;
        \draw    (369.75,154) -- (426.75,154) ;
        \draw [shift={(428.75,154)}, rotate = 180] [color={rgb, 255:red, 0; green, 0; blue, 0 }  ][line width=0.75]    (10.93,-3.29) .. controls (6.95,-1.4) and (3.31,-0.3) .. (0,0) .. controls (3.31,0.3) and (6.95,1.4) .. (10.93,3.29)   ;
        \draw    (369.75,68) -- (426.75,68) ;
        \draw [shift={(428.75,68)}, rotate = 180] [color={rgb, 255:red, 0; green, 0; blue, 0 }  ][line width=0.75]    (10.93,-3.29) .. controls (6.95,-1.4) and (3.31,-0.3) .. (0,0) .. controls (3.31,0.3) and (6.95,1.4) .. (10.93,3.29)   ;
        \draw [color={rgb, 255:red, 0; green, 0; blue, 0 }  ,draw opacity=1 ][fill={rgb, 255:red, 189; green, 16; blue, 224 }  ,fill opacity=1 ]   (129.75,68.5) -- (189.59,152.37) ;
        \draw [shift={(190.75,154)}, rotate = 234.49] [color={rgb, 255:red, 0; green, 0; blue, 0 }  ,draw opacity=1 ][line width=0.75]    (10.93,-3.29) .. controls (6.95,-1.4) and (3.31,-0.3) .. (0,0) .. controls (3.31,0.3) and (6.95,1.4) .. (10.93,3.29)   ;
        \draw    (250.25,154) -- (308.11,70.15) ;
        \draw [shift={(309.25,68.5)}, rotate = 124.61] [color={rgb, 255:red, 0; green, 0; blue, 0 }  ][line width=0.75]    (10.93,-3.29) .. controls (6.95,-1.4) and (3.31,-0.3) .. (0,0) .. controls (3.31,0.3) and (6.95,1.4) .. (10.93,3.29)   ;
        \draw    (369.75,68) -- (427.62,152.35) ;
        \draw [shift={(428.75,154)}, rotate = 235.55] [color={rgb, 255:red, 0; green, 0; blue, 0 }  ][line width=0.75]    (10.93,-3.29) .. controls (6.95,-1.4) and (3.31,-0.3) .. (0,0) .. controls (3.31,0.3) and (6.95,1.4) .. (10.93,3.29)   ;
        \draw    (470.75,153.5) -- (528.61,69.65) ;
        \draw [shift={(529.75,68)}, rotate = 124.61] [color={rgb, 255:red, 0; green, 0; blue, 0 }  ][line width=0.75]    (10.93,-3.29) .. controls (6.95,-1.4) and (3.31,-0.3) .. (0,0) .. controls (3.31,0.3) and (6.95,1.4) .. (10.93,3.29)   ;
        \draw  [color={rgb, 255:red, 255; green, 0; blue, 0 }  ,draw opacity=1 ][fill={rgb, 255:red, 208; green, 2; blue, 27 }  ,fill opacity=0.1 ][dash pattern={on 4.5pt off 4.5pt}] (120.4,48.6) .. controls (120.4,39.76) and (127.56,32.6) .. (136.4,32.6) -- (184.4,32.6) .. controls (193.24,32.6) and (200.4,39.76) .. (200.4,48.6) -- (200.4,165.8) .. controls (200.4,174.64) and (193.24,181.8) .. (184.4,181.8) -- (136.4,181.8) .. controls (127.56,181.8) and (120.4,174.64) .. (120.4,165.8) -- cycle ;
        \draw  [color={rgb, 255:red, 255; green, 0; blue, 0 }  ,draw opacity=1 ][fill={rgb, 255:red, 208; green, 2; blue, 27 }  ,fill opacity=0.1 ][dash pattern={on 4.5pt off 4.5pt}] (240.4,48.4) .. controls (240.4,39.56) and (247.56,32.4) .. (256.4,32.4) -- (304.4,32.4) .. controls (313.24,32.4) and (320.4,39.56) .. (320.4,48.4) -- (320.4,165.6) .. controls (320.4,174.44) and (313.24,181.6) .. (304.4,181.6) -- (256.4,181.6) .. controls (247.56,181.6) and (240.4,174.44) .. (240.4,165.6) -- cycle ;
        \draw  [color={rgb, 255:red, 255; green, 0; blue, 0 }  ,draw opacity=1 ][fill={rgb, 255:red, 208; green, 2; blue, 27 }  ,fill opacity=0.1 ][dash pattern={on 4.5pt off 4.5pt}] (360.2,48.2) .. controls (360.2,39.36) and (367.36,32.2) .. (376.2,32.2) -- (424.2,32.2) .. controls (433.04,32.2) and (440.2,39.36) .. (440.2,48.2) -- (440.2,165.4) .. controls (440.2,174.24) and (433.04,181.4) .. (424.2,181.4) -- (376.2,181.4) .. controls (367.36,181.4) and (360.2,174.24) .. (360.2,165.4) -- cycle ;
        \draw  [color={rgb, 255:red, 255; green, 0; blue, 0 }  ,draw opacity=1 ][fill={rgb, 255:red, 208; green, 2; blue, 27 }  ,fill opacity=0.1 ][dash pattern={on 4.5pt off 4.5pt}] (459.6,48.4) .. controls (459.6,39.56) and (466.76,32.4) .. (475.6,32.4) -- (523.6,32.4) .. controls (532.44,32.4) and (539.6,39.56) .. (539.6,48.4) -- (539.6,165.6) .. controls (539.6,174.44) and (532.44,181.6) .. (523.6,181.6) -- (475.6,181.6) .. controls (466.76,181.6) and (459.6,174.44) .. (459.6,165.6) -- cycle ;
        
        \draw (76.5,21) node [anchor=north west][inner sep=0.75pt]  [xscale=1,yscale=1] [align=left] {Input};
        \draw (540,20) node [anchor=north west][inner sep=0.75pt]  [xscale=1,yscale=1] [align=left] {Output};
        \draw (90,65) node [anchor=north west][inner sep=0.75pt]  [xscale=1,yscale=1]  {$ \begin{array}{l}
        x_{0}\\
        \\
        \\
        \\
        \\
        \\
        y_{0}
        \end{array}$};
        \draw (211.5,65) node [anchor=north west][inner sep=0.75pt]  [xscale=1,yscale=1]  {$ \begin{array}{l}
        x_{1}\\
        \\
        \\
        \\
        \\
        \\
        y_{1}
        \end{array}$};
        \draw (331.5,65) node [anchor=north west][inner sep=0.75pt]  [xscale=1,yscale=1]  {$ \begin{array}{l}
        x_{2}\\
        \\
        \\
        \\
        \\
        \\
        y_{2}
        \end{array}$};
        \draw (551.5,65) node [anchor=north west][inner sep=0.75pt]  [xscale=1,yscale=1]  {$ \begin{array}{l}
        x_{q}\\
        \\
        \\
        \\
        \\
        \\
        y_{q}
        \end{array}$};
        \draw (442.5,97) node [anchor=north west][inner sep=0.75pt]  [xscale=1,yscale=1] [align=left] {...};
        \draw (157,48.9) node [anchor=north west][inner sep=0.75pt]  [xscale=1,yscale=1]  {$I$};
        \draw (493,47.4) node [anchor=north west][inner sep=0.75pt]  [xscale=1,yscale=1]  {$I$};
        \draw (493,156.9) node [anchor=north west][inner sep=0.75pt]  [xscale=1,yscale=1]  {$I$};
        \draw (388.5,47.4) node [anchor=north west][inner sep=0.75pt]  [xscale=1,yscale=1]  {$I$};
        \draw (389.5,157.4) node [anchor=north west][inner sep=0.75pt]  [xscale=1,yscale=1]  {$I$};
        \draw (274.5,47.9) node [anchor=north west][inner sep=0.75pt]  [xscale=1,yscale=1]  {$I$};
        \draw (274,156.9) node [anchor=north west][inner sep=0.75pt]  [xscale=1,yscale=1]  {$I$};
        \draw (157.5,156.4) node [anchor=north west][inner sep=0.75pt]  [xscale=1,yscale=1]  {$I$};
        \draw (169.9,8.9) node [anchor=north west][inner sep=0.75pt]  [xscale=1,yscale=1]  {${\textcolor[rgb]{1,0,0}{T}}\textcolor[rgb]{1,0,0}{_{1}}$};
        \draw (472.7,9.3) node [anchor=north west][inner sep=0.75pt]  [xscale=1,yscale=1]  {${\textcolor[rgb]{1,0,0}{T}}\textcolor[rgb]{1,0,0}{_{q}}$};
        \draw (390.7,9) node [anchor=north west][inner sep=0.75pt]  [xscale=1,yscale=1]  {${\textcolor[rgb]{1,0,0}{T}}\textcolor[rgb]{1,0,0}{_{3}}$};
        \draw (271.6,9.5) node [anchor=north west][inner sep=0.75pt]  [xscale=1,yscale=1]  {${\textcolor[rgb]{1,0,0}{T}}\textcolor[rgb]{1,0,0}{_{2}}$};
        \draw (159.6,86.6) node [anchor=north west][inner sep=0.75pt]  [xscale=1,yscale=1]  {$\hat{\sigma} _{1}$};
        \draw (262.4,86.6) node [anchor=north west][inner sep=0.75pt]  [xscale=1,yscale=1]  {$\hat{\sigma} _{2}$};
        \draw (400.4,86.6) node [anchor=north west][inner sep=0.75pt]  [xscale=1,yscale=1]  {$\hat{\sigma} _{3}$};
        \draw (479.6,86.6) node [anchor=north west][inner sep=0.75pt]  [xscale=1,yscale=1]  {$\hat{\sigma} _{q}$};

        \end{tikzpicture}
    \end{center}
    \caption{SympNet architecture.
        Here the $(T_i)_{1\leq i \leq q}$ are shear symplectic maps, while the $(\hat\sigma_i)_{1\leq i \leq q}$ denote their associated gradient modules.}
    \label{fig:figure_symp_net}
\end{figure}

Finally, the architecture of SympNets is made
by composing several gradient modules,
as shown in \cref{fig:figure_symp_net}.

Once the SympNet has been constructed,
a loss function has to be designed to optimize
all the trainable parameters contained in gradient modules.
To that end, the starting space is assumed to be included in a compact set $K$, here the sphere~$\mathbb S^1$. First, we draw points $\lbrace x^i \rbrace_{i=1}^{N_b}$ on the boundary of $K$, and evaluate their image by~$\mathcal{T}$, a given symplectic map that the SympNet $T_\omega$ with trainable weights $\omega$ will be trained to learn.
    {\rb If one wants to approximate a certain symplectic map $\mathcal T$ in terms of SympNets, one would have to minimize the following loss function
        \begin{equation*}
            \mathcal{J}_{S}\left(\omega; \lbrace x^i \rbrace_{i=1}^{N}\right) =
            \sum_{i=1}^{N} \lvert T_{\omega}(x^i) - \mathcal{T}(x^i)\rvert^2,
        \end{equation*}
        with $T_{\omega}$ the SympNet approximation of $\mathcal T$ and $|\cdot|$ the euclidean norm on $\mathbb{R}^2$.}

Equipped with PINNs (or rather DeepRitz networks)
introduced in \cref{sec:pinns}
as well as SympNets,
we now combine them in the next section to solve
the shape optimization problem \eqref{eq:optim_ener}.



\subsection{Solving a PDE with PINNs in a domain generated with a symplectic map}
\label{sec:shape_optim_with_sympnets}

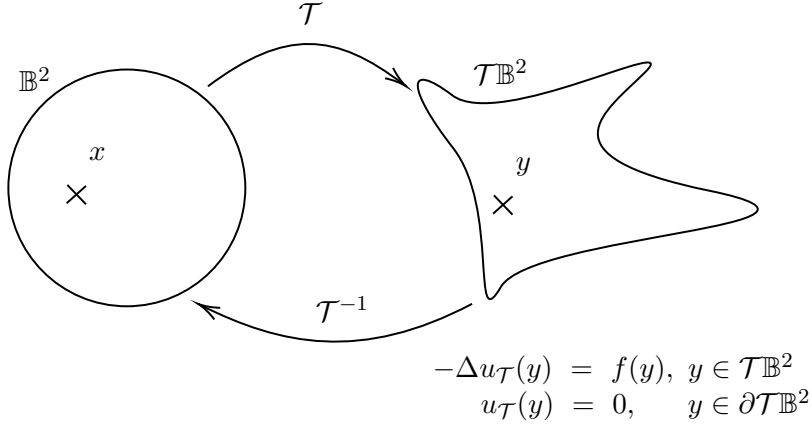
\begin{figure}[!ht]
    \centering
        \tikzset{every picture/.style={line width=0.75pt}} 

        \begin{tikzpicture}[x=0.75pt,y=0.75pt,yscale=-1,xscale=1]
        
        \draw   (163,93.75) .. controls (163,60.75) and (189.75,34) .. (222.75,34) .. controls (255.75,34) and (282.5,60.75) .. (282.5,93.75) .. controls (282.5,126.75) and (255.75,153.5) .. (222.75,153.5) .. controls (189.75,153.5) and (163,126.75) .. (163,93.75) -- cycle ;
        \draw    (192.75,92.25) -- (202.13,102) ;
        \draw    (192.38,101.5) -- (202.13,92.75) ;
        \draw   (386.88,47) .. controls (410.63,65) and (504.25,15.25) .. (484.25,35.25) .. controls (464.25,55.25) and (429.5,77.25) .. (516.88,95.25) .. controls (604.25,113.25) and (426.81,117.18) .. (411.38,143.75) .. controls (395.94,170.32) and (410.7,103.77) .. (387.38,73.75) .. controls (364.05,43.73) and (363.13,29) .. (386.88,47) -- cycle ;
        \draw    (407.88,97.5) -- (417.25,107.25) ;
        \draw    (407.5,106.75) -- (417.25,98) ;
        \draw    (263.63,42.5) .. controls (303.03,12.95) and (328.11,15.42) .. (362.07,41.3) ;
        \draw [shift={(363.63,42.5)}, rotate = 217.85] [color={rgb, 255:red, 0; green, 0; blue, 0 }  ][line width=0.75]    (10.93,-3.29) .. controls (6.95,-1.4) and (3.31,-0.3) .. (0,0) .. controls (3.31,0.3) and (6.95,1.4) .. (10.93,3.29)   ;
        \draw    (396.88,152.25) .. controls (341.06,181.7) and (300.32,173.18) .. (259.85,152.87) ;
        \draw [shift={(258.63,152.25)}, rotate = 26.91] [color={rgb, 255:red, 0; green, 0; blue, 0 }  ][line width=0.75]    (10.93,-3.29) .. controls (6.95,-1.4) and (3.31,-0.3) .. (0,0) .. controls (3.31,0.3) and (6.95,1.4) .. (10.93,3.29)   ;
        
        \draw (167,31.4) node [anchor=north west][inner sep=0.75pt]    {$\mathbb B^2$};
        \draw (202.25,71.65) node [anchor=north west][inner sep=0.75pt]    {$x$};
        \draw (417.38,76.9) node [anchor=north west][inner sep=0.75pt]    {$y$};
        \draw (395.88,27.4) node [anchor=north west][inner sep=0.75pt]    {$\mathcal T \mathbb B^2$};
        \draw (307.88,-1.1) node [anchor=north west][inner sep=0.75pt]    {$\mathcal{T}$};
        \draw (316.63,147.15) node [anchor=north west][inner sep=0.75pt]    {$\mathcal{T}^{-1}$};
        \draw (329.93,175.25) node [anchor=north west][inner sep=0.75pt]    {$ \begin{array}{l}
        \ \ \ \ \ \ \ \ -\Delta u_{\mathcal{T}}( y) \ =\ f( y) ,\ y\in \mathcal T \mathbb B^2\\
        \ \ \ \ \ \ \ \ \ \ \ \ \ u_{\mathcal{T}}( y) \ =\ 0,\ \ \ \ \ y\in \partial \mathcal T \mathbb B^2
        \end{array}$};

        \end{tikzpicture}    
    \caption{Transformation of the ball $\mathbb B^2$ by a symplectic map $\mathcal{T}$.}\label{fig:transfo}
\end{figure}

Now that we know how to learn a symplectic transformation with a SympNet, we now wish to combine this approach with PINNs to solve the Poisson problem \eqref{eq:poisson} in this shape.
For simplicity, we consider the 2D case,
i.e., we take $d=1$ and transform the unit sphere $\mathbb B^2$
into a shape $\mathcal{T} \mathbb B^2$ homeomorphic to the ball with volume $V_0$.
However, directly adapting the approach described in \cref{sec:pinns}
leads to difficulties in properly defining the boundary conditions.
Indeed, the solution of the PDE
should be expressed as in \eqref{eq:bc_net},
where the trained network $u_\theta$ is
multiplied by a function $\alpha$ and
summed to a function~$\beta$.
In this case, we would have to devise new expressions of
$\alpha$ and $\beta$ in $\mathcal {T}\mathbb B^2$.
Faced with this difficulty, we choose to bypass that problem,
by solving the PDE satisfied by
$w:\mathbb B^2 \to \mathbb{R}$,
defined, for a.e. $x\in \mathbb B^2$ and $y\in \mathcal T \mathbb B^2$ such that $y=\mathcal T x$, by
\begin{equation*}
    w(x) = (u_{\mathcal{T}} \circ \mathcal T)(x) = u_{\mathcal{T}}(y),
\end{equation*}
with $u_\mathcal{T}$ the solution of the Poisson problem \eqref{eq:poisson} in $\mathcal T \mathbb B^2$.
Since the inputs of the function $v$ lie in the unit sphere, it is again very easy to find appropriate functions $\alpha$ and $\beta$. \cref{fig:transfo} recaps this notation.

{\rb Given $w$, we now determine the PDE it satisfies}. The answer lies in the following lemma (that remains true in higher dimensions).

\begin{lemma}
    \label{lem:formulation_with_metric_tensor}
    Let $\mathcal{T}$ be a symplectic map on $\mathbb R^2$ and $u_\mathcal{T}\in H^1_0(\mathcal T\mathbb B^2)$, such that $w=u_\mathcal{T}\circ\mathcal{T}\in H^1_0(\mathbb B^2)$. If $u_\mathcal{T}\in H^1_0(\mathcal T \mathbb B^2)$ is the solution of the Poisson problem~\eqref{eq:poisson}, then $w\in H^1_0(\mathbb B^2)$ is solution of
    \begin{equation} \label{eq:diff_form}
        \begin{cases}
            -\mathrm{div}\big( A \nabla w \big) = \tilde{f}, & \text{ in } \mathbb B^2; \\
            w = 0,                                           & \text{ on } \mathbb S^1,
        \end{cases}
    \end{equation}
    with $A:\mathbb B^2 \to \mathcal{M}_2(\mathbb{R})$ is defined by
    $A = J_\mathcal{T}^{-1}\cdot J_\mathcal{T}^{-\intercal}$, with $\tilde{f} = f \circ \mathcal{T}$ and with $J_\mathcal{T}$ the jacobian matrix of $\mathcal T$, {\rb and $J_\mathcal{T}^{-\intercal} = \left(J_\mathcal{T}^{-1}\right)^\intercal = \left(J_\mathcal{T}^\intercal\right)^{-1}$}.
    Moreover, the following problem can be formulated in a weaker sense, as an optimization problem
    \begin{equation}
        \label{eq:optim_form}
        \inf_{}
        \bigg\lbrace\frac12 \int_{\mathbb B^2}
        A \nabla w\cdot \nabla w
        - \int_{\mathbb B^2}\tilde{f}w,
        \quad \exists u_\mathcal{T}\in H^1_0(\mathcal{T}\mathbb B^1),
        \quad w=u_\mathcal{T}\circ\mathcal{T}\in H^1_0(\mathbb B^1)
        \bigg\rbrace.
    \end{equation}
\end{lemma}
\begin{proof}
    Assume $u_\mathcal{T}$ solves the Poisson problem \eqref{eq:poisson}.
    Then, for all $\varphi \in H^1_0(\mathcal{T}\mathbb B^1)$, the following variational formulation holds:
    \begin{equation*}
        \int_{\mathcal{T}\mathbb B^2} \nabla u_\mathcal{T}(y)\cdot\nabla\varphi(y)\mathrm{dy}
        =
        \int_{\mathcal{T}\mathbb B^2}f(y)\varphi(y)\mathrm{dy}.
    \end{equation*}
    Proceeding with the change of variable $y =\mathcal{T} x$, we get
    \begin{equation*}
        \int_{\mathbb B^2}
        |J_\mathcal{T}(x)|\nabla u_\mathcal{T}(\mathcal{T}x)\cdot\nabla\varphi(\mathcal{T}x)\mathrm{dx}
        =
        \int_{\mathbb B^2} |J_\mathcal{T}(x)|f(\mathcal{T}x)\varphi(\mathcal{T}x)\mathrm{dx},
    \end{equation*}
    with $|J_\mathcal{T}|$ the determinant of the Jacobian matrix of $\mathcal T$.
    We now introduce $w=u_\mathcal{T}\circ\mathcal{T}$, $\tilde{f}=f\circ\mathcal{T}$ and $\tilde{\varphi}=\varphi\circ\mathcal{T}$. By using the chain rule, we get for a.e. $x\in \mathbb B^2$,
    $\nabla u_\mathcal{T}(\mathcal{T}x)=J_\mathcal{T}^{-\intercal}(x)\nabla v(x)$ and
    $\nabla \varphi(\mathcal{T}x) = J_\mathcal{T}^{-\intercal}(x)\nabla \tilde \varphi(x)$.
    Since $\mathcal{T}$ is a symplectic map and therefore is volume-preserving, $|J_\mathcal{T}(x)|=1$,
    and we obtain
    \begin{equation}\label{eq:avant_derniere_ligne}
        \int_{\mathbb B^2}
        \bigg[ J_\mathcal{T}^{-\intercal}(x)\nabla w(x) \bigg]
        \cdot
        \bigg[ J_\mathcal{T}^{-\intercal}(x)\nabla \tilde \varphi(x) \bigg]
        \mathrm{dx}
        =
        \int_{\mathbb B^2} \tilde f(x)\tilde \varphi(x)\mathrm{dx}.
    \end{equation}
    Then, \eqref{eq:avant_derniere_ligne} directly leads us to the variational problem of finding $v\in H^1_0(\mathbb B^2)$, such that, for all $\tilde{\varphi} \in H^1_0(\mathbb B^2)$,
    \begin{equation}
        \label{eq:fv_diff_forme}
        \int_{\mathbb B^2}
        \bigg[ J_\mathcal{T}^{-1}(x)J_\mathcal{T}^{-\intercal}(x)\bigg]
        \nabla w(x) \cdot \nabla \tilde{\varphi}(x)
        \mathrm{dx}
        =
        \int_{\mathbb B^2}\tilde{f}(x)\tilde{\varphi}(x)\mathrm{dx}.
    \end{equation}
    The proof is thus concluded.
\end{proof}

\begin{remark}
    Note that the proof of \cref{lem:formulation_with_metric_tensor}
    follows a similar reasoning as in the proof of~\cite[Theorem 5.3.2]{henrot-pierre}.
    Indeed, \cref{lem:formulation_with_metric_tensor} can formally be seen
    as computing a shape derivative, but here the shape derivative is simply represented
    through the derivative of the symplectic map.
\end{remark}

\subsection{Shape optimization with PINNs and SympNets}
\label{sec:strategy_shapo}

To propose a shape optimization algorithm,
we now just have to minimize a loss function depending on both PINNs and SympNets.
It is given by
\begin{equation}
    \label{eq:loss_optim_forme}
    \mathcal{J}_{P/S}\big(\theta, \omega; \lbrace x^i, \mu^i \rbrace_{i=1}^{N}\big) =
    \frac{V_0}{N} \sum_{i=1}^{N}
    \left\lbrace
    \frac12 A_{\omega} \nabla v_{\theta,\omega} \cdot \nabla v_{\theta,\omega}
    - \tilde{f}_{\omega}v_{\theta,\omega} \right\rbrace (x^i; \mu^i),
\end{equation}
with
\begin{itemize}
    \item $\theta$ the trainable weights of the PINN, $\omega$ the trainable weights of the SympNet;
    \item $v_{\theta,\omega}:x \in \mathbb B^2 \mapsto \alpha(x)u_\theta(T_\omega x)+\beta(x) \in \mathbb{R}$ the solution of the Poisson problem set in $T_\omega\mathbb B^2$;
    \item $T_{\omega}:\mathbb{R}^{2d} \to \mathbb{R}^{2d}$ the SympNet;
    \item $\mathrm{A}_\omega$ the diffusion matrix defined by $\mathrm{A}_\omega = \mathrm{J}_{T_\omega}^{-1}\mathrm{J}_{T_\omega}^{-\intercal}$;
    \item $u_{\theta}:{T}_{\omega}\mathbb B^2\to \mathbb{R}$ the PINN;
    \item $\alpha: \mathbb B^2 \mapsto \mathbb R$ a $C^\infty$ function that vanishes on $\mathbb S^1$;
    \item $\beta: \mathbb B^2 \mapsto \mathbb R$ a $C^\infty$ function that satisfies the boundary condition of the Poisson problem (here, $\beta=0$ on $\mathbb S^1$);
    \item $\tilde{f}_\omega:x \in \mathbb B^2 \mapsto (f\circ T_\omega)(x) \in \mathbb R$;
    \item $N$ the number of collocation points;
    \item $\lbrace x^i, \mu^i \rbrace_{i=1}^{N}$ a set of random collocation points and values of the parameters.
\end{itemize}

\begin{remark}
    If we uniformly sample points $\lbrace x^i \rbrace_i$ in the unit ball $\mathbb B^2$ and apply  a symplectic transformation~$\mathcal{T}$ to obtain a shape $\mathcal{T}\mathbb B^2$, then the points $\lbrace \mathcal{T}x^i \rbrace_i$ are again uniformly distributed in $\mathcal{T}\mathbb B^2$. This is immediate, since $\mathcal{T}$ is volume-preserving.
\end{remark}

\begin{remark}
    One of the exciting aspects of this formulation of the Dirichlet energy optimization problem is that we iterate the gradient descent on the SympNet and PINN at the same time. There is no need to wait for the PDE to be solved before iterating on the shape. The gradient descent is now parallelizable.
\end{remark}

\begin{remark}
    The major advantage of this formulation is that the shape derivative can be calculated using automatic differentiation techniques.
\end{remark}

{\rb %
\begin{remark}
    Note that adapting our method to higher dimensions
    would require changing the SympNet architecture,
    firstly because they are restricted to even dimensions,
    and secondly because the equivalence between being a symplectomorphism
    and being a volume-preserving differentiable map is only true in $\mathbb R^2$.
    Despite these restrictions only due to the network representing the shape,
    the global methodology is applicable to higher dimensions.
    For instance, the DeepRitz method works in any dimension,
    and the symplecticity hypothesis in \cref{lem:formulation_with_metric_tensor}
    can be relaxed to the volume-preserving property
    of the $C^1$-diffeomorphism $\mathcal T$.
\end{remark}%
}

\begin{remark}
    At each iteration, contrary to classical shape derivative methods, the SympNet representing the shape and the PINN representing the PDE are updated, even if the PINN gives a poor approximation of the PDE solution within the shape. The PINN is expected to solve Poisson's equation \eqref{eq:poisson} in the shape only when the algorithm has converged.
    This allows us to parallelize the entire loop and go beyond conventional methods, which are intrinsically sequential.
\end{remark}

{\rb %
\begin{remark}
    A natural question one can ask is whether
    the neural networks could be replaced with classical methods.
    On the one hand, replacing DeepRitz with a finite element method
    would alter our method's mesh-less nature (since (re)meshing would be needed)
    and its intrinsic parallelizability
    (since joint optimization would no longer be possible).
    On the other hand, coupling PINNs (to solve the PDE)
    with a classical method (to represent and optimize the shape)
    has the same problems, and has already been documented in~\cite{odot2023real}.
\end{remark}%
}

In \cref{alg:cap}, we summarize the main steps
of the shape optimization algorithm.

\begin{algorithm}[!ht]
    \caption{Shape optimization algorithm}\label{alg:cap}
    \begin{algorithmic}
        \Require \
        \begin{itemize}
            \item {\it collocation domain:} $\mathbb B^2\times \mathbb M$;
            \item {\it PINN and SympNet hyperparameters:} see \cref{appendix:hyperparameters};
            \item {\it algorithm parameters:} gradient descent step, number of epochs, number of collocation points $N$, see \cref{appendix:hyperparameters};
        \end{itemize}
        \While{epoch $<$ total number of epochs}
        \begin{itemize}
            \item[$\leadsto$] randomly draw $N$ collocation points $\{x^i,\mu^i\}_{i=1}^N\in \left(B^2\times\mathbb M\right)^N$;
            \item[$\leadsto$] for each collocation point $(x^i,\mu^i)$,
                  \textit{this step is vectorizable and parallelizable}:
                  \begin{itemize}
                      \item determination of the transformation $T_\omega(x^i,\mu^i)$; 
                      \item knowing $T_\omega(x^i,\mu^i)$, computation of the diffusion matrix $A_\omega(x^i,\mu^i)$ and of the source term $\tilde{f}_\omega(x^i,\mu^i)$ involved in \eqref{eq:diff_form};
                      \item computation of $v_{\theta,\omega}$, the composition of the PINN by the SympNet at each point $(x^i, \mu^i)$;
                  \end{itemize}
            \item[$\leadsto$] calculation of the loss given by \eqref{eq:loss_optim_forme};
            \item[$\leadsto$] calculation by automatic differentiation of the loss derivative with respect to PINN and SympNet trainable weights $(\theta,\omega)$;
            \item[$\leadsto$] computation of a gradient method step and update of the PINN and SympNet trainable weights.
        \end{itemize}
        \EndWhile
    \end{algorithmic}
\end{algorithm}

\section{Numerical results with Dirichlet boundary conditions}
\label{sec:numerics}

Equipped the algorithm from \cref{sec:strategy_shapo},
we are now ready to provide a comprehensive numerical study.
As a sanity check, we first present the results of the PINN method in \cref{sec:poisson_with_PINN}
and of the SympNets in \cref{sec:sympnet_numerics}.
When no exact solution is available, the results of the code will be compared to results obtained via a fixed point algorithm from \cite[Algorithm 1]{chambolle:hal-04140177}, implemented in \texttt{FreeFem++} \cite{MR3043640}.
To test our SympNets, we introduce the symplectic map
$\mathcal T_\lambda = \mathcal S^1_\lambda \circ \mathcal S^2_\lambda: \mathbb{R}^2 \to \mathbb{R}^2$, with
\begin{equation}
    \label{eq:bizaroid}
    \begin{cases}
        \mathcal S^1_\lambda:
        (x_1,x_2)
        \mapsto
        (x_1 - \lambda x_2^2+ 0.3 \sin(\frac {x_2}\lambda) - 0.2  \sin(8x_2), \, x_2), \\
        \mathcal S^2_\lambda:
        (x_1,x_2)
        \mapsto
        (x_1, \, x_2 + 0.2\lambda x_1 + 0.12 \cos(x_1)).
    \end{cases}
\end{equation}
Note that $\mathcal T_\lambda$ is a nonlinear function of $\lambda$.
This map will be the support of several numerical experiments.
Then, we combine both approaches in \cref{sec:gesonn_numerics}
to solve the shape optimization problem \eqref{eq:optim_ener}.
Finally, we solve a more intricate problem, namely the exterior Bernoulli free-boundary problem, in \cref{sec:bernoulli}.
The hyperparameters of the networks are summarized in \cref{tab:hyperparameters}.

\subsection{Solving the Poisson problem with the PINN method}
\label{sec:poisson_with_PINN}

\begin{figure}[!ht]
    \centering
    \begin{tikzpicture}
        \node[xshift=-0.25\textwidth] (top_left) at (0,0)
        {\includegraphics[width=0.45\textwidth]{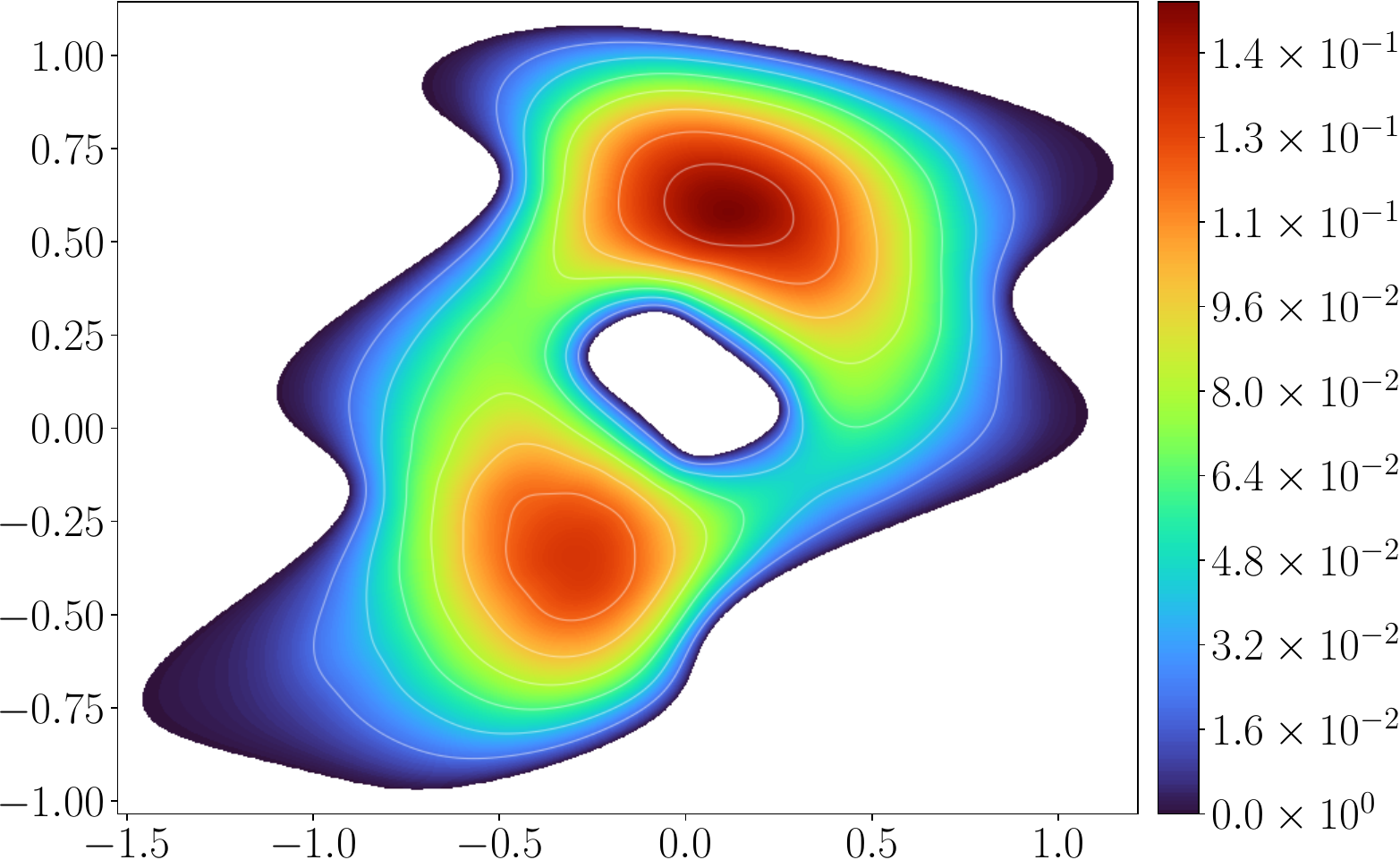}};
        \node[yshift=-0.02\textwidth] at (top_left.south) {(a) solution, $\mu=0.5$};

        \node[xshift=0.25\textwidth] (top_right) at (top_left.east)
        {\includegraphics[width=0.45\textwidth]{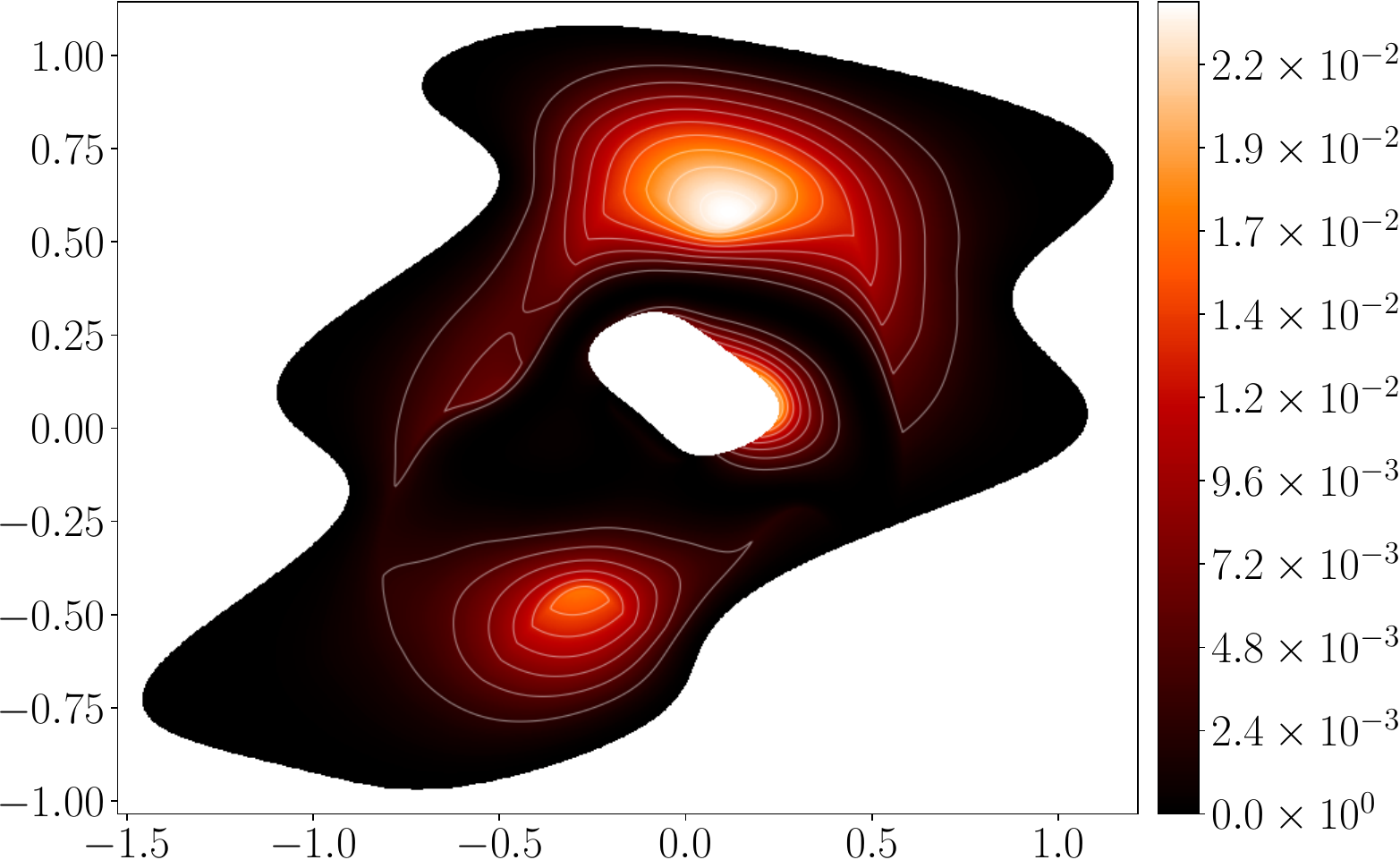}};
        \node[yshift=-0.02\textwidth] at (top_right.south) {(b) error, $\mu=0.5$};

        \node[yshift=-0.21\textwidth] (bottom_left) at (top_left.south)
        {\includegraphics[width=0.45\textwidth]{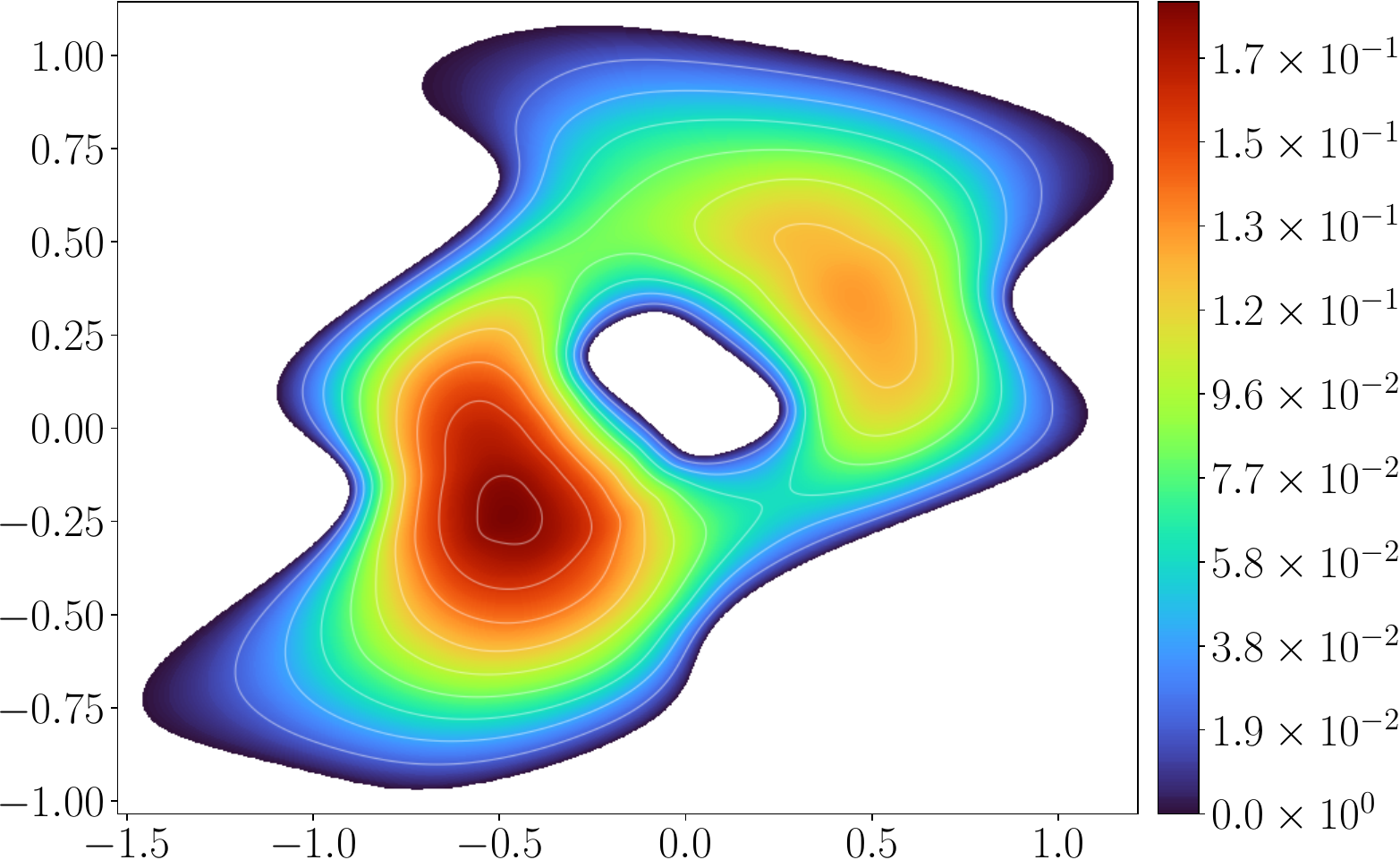}};
        \node[yshift=-0.02\textwidth] at (bottom_left.south) {(c) solution, $\mu=1.5$};

        \node[yshift=-0.21\textwidth] (bottom_right) at (top_right.south)
        {\includegraphics[width=0.45\textwidth]{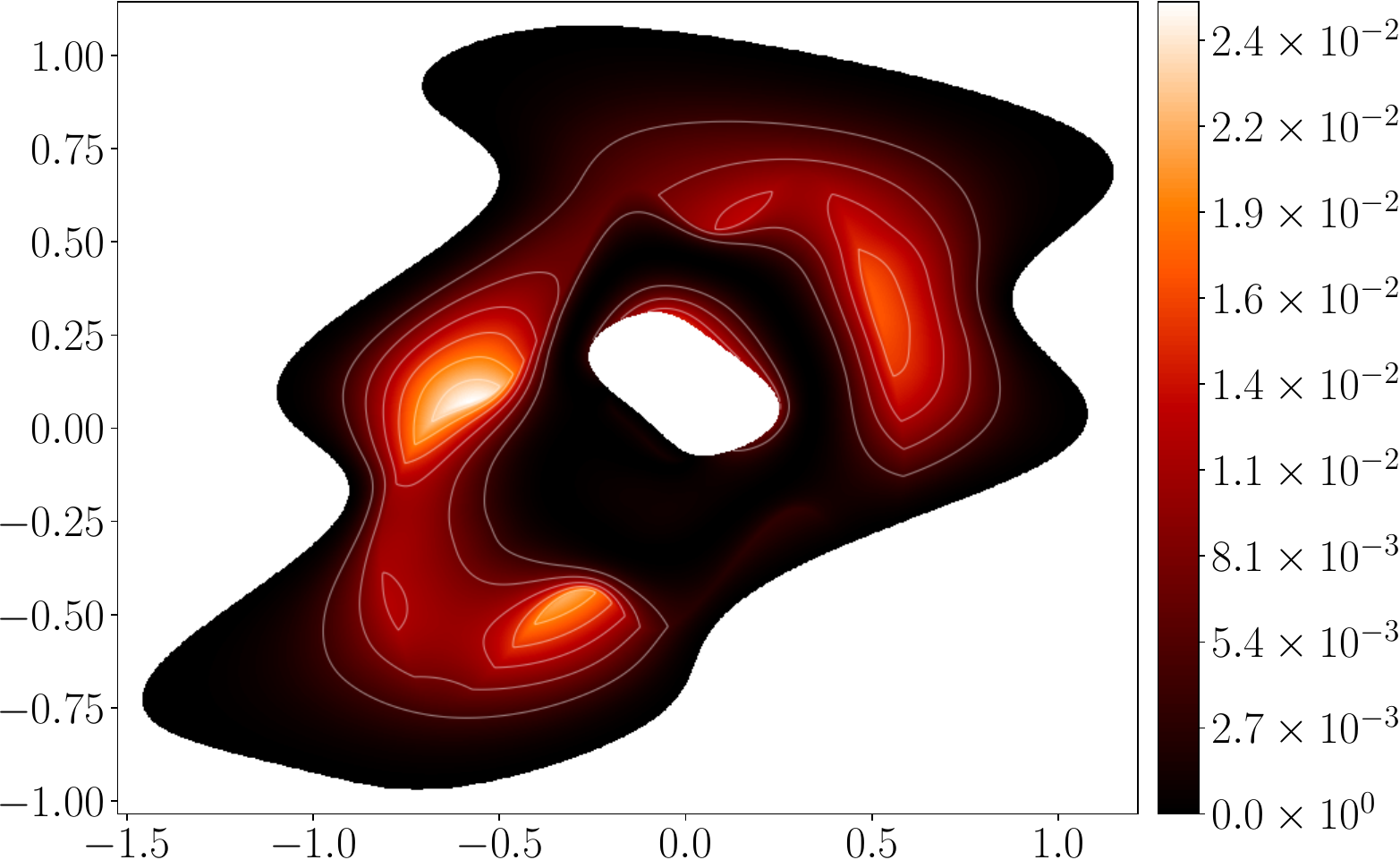}};
        \node[yshift=-0.02\textwidth] at (bottom_right.south) {(d) error, $\mu=1.5$};
    \end{tikzpicture}
    \caption{%
        Results of the PINN applied to the Poisson problem~\eqref{eq:poisson}
        with the source term \eqref{eq:parametric_source},
        on the domain given by applying the symplectic map~\eqref{eq:bizaroid}
        with $\lambda = 0.5$ to the annulus with inner radius $0.2$ and outer radius $1$.
            {\rall Top panels: results with $\mu = 0.5$;
                bottom panels: results with $\mu = 1.5$.
                Left panels: approximate solutions;
                right panels: errors with respect to a finite element simulation
                on a fine grid ($1000$ nodes per edge, around $10^6$ elements).}
    }
    \label{fig:PINN_approximation}
\end{figure}
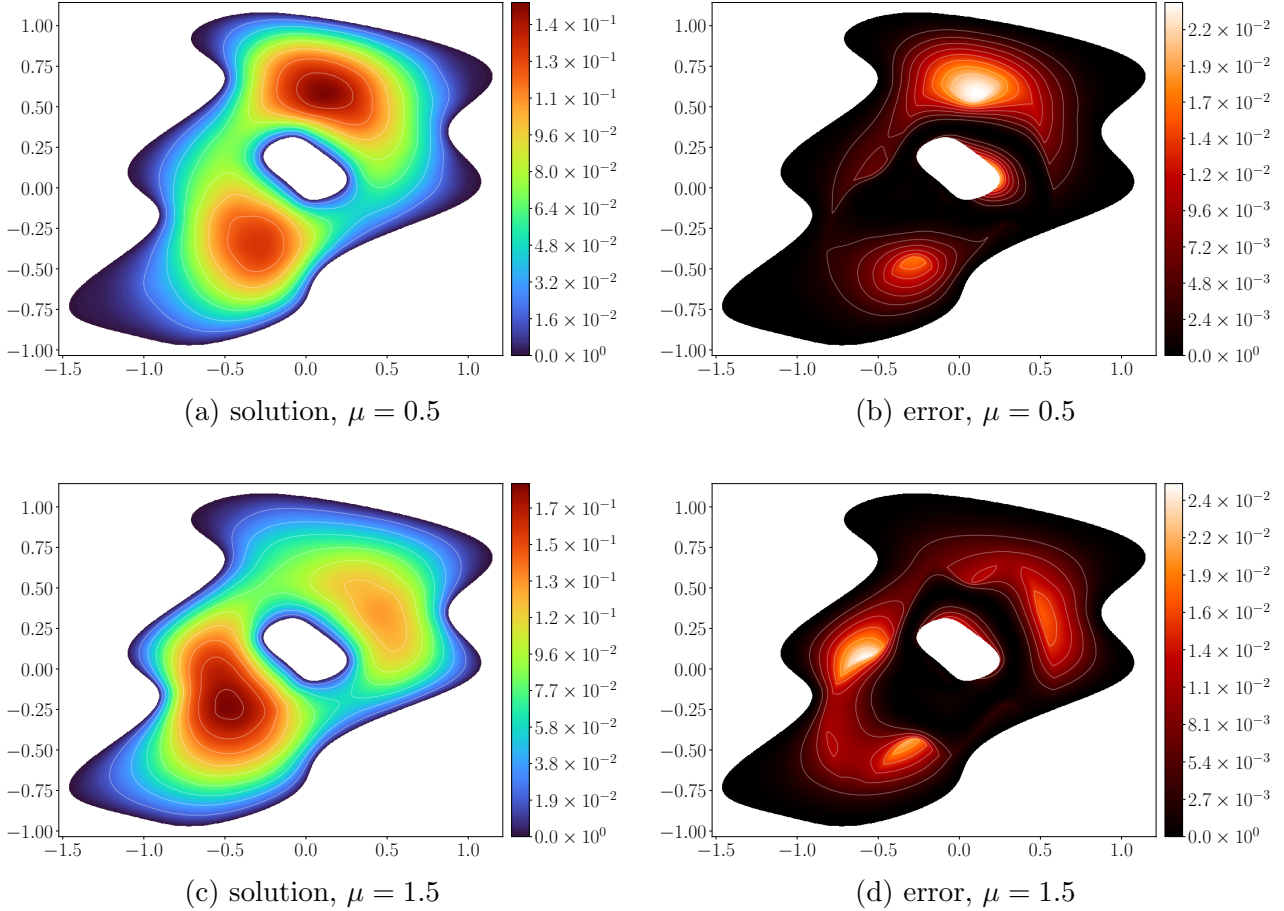

This section is dedicated to checking that the
energetic formulation \eqref{eq:optim_form}
introduced in \cref{lem:formulation_with_metric_tensor}
is indeed able to solve the Poisson problem \eqref{eq:poisson}
in a given shape $\Upomega$.
For this purpose, we choose a complex shape $\Upomega$
given by the symplectic map \eqref{eq:bizaroid}
with $\lambda = 0.5$
applied to the annulus with inner radius $0.2$ and outer radius $1$.
This leads to a shape with a hole.
More specifically, we highlight the capability of
our neural network-based approach to tackle a parametric
Poisson problem, with a parametric source term $f$, given by
\begin{equation}
    \label{eq:parametric_source}
    f(x_1, x_2; \mu) = \exp \left(
    1 - \left(\frac {x_1} \mu \right)^2 - \left(\mu x_2\right)^2
    \right),
\end{equation}
where the parameter $\mu \in \mathbb{M} = (0.5, 1.5)$
controls the elongation of the ellipse
described by the source term.
Therefore, the neural network is a function of the two space variables $x_1$ and $x_2$, as well as of the parameter~$\mu$.
The results are displayed on \cref{fig:PINN_approximation}.

To quantify the approximation quality,
we cannot directly use the Dirichlet energy~\eqref{eq:energy},
since we do not know its optimal value for this specific problem.
Instead, we elect to use the variational formulation \eqref{eq:fv_diff_forme}.
Indeed, we know that the solution to the Poisson problem
satisfies, for all $\tilde{\varphi} \in H^1_0(\mathbb B^2)$,
\begin{equation}
    \label{eq:fv_diff_forme_for_stats}
    \int_{0.5}^{1.5}
    \int_{\mathbb B^2}
    \bigg(
    \bigg[ J_\mathcal{T}^{-1}(x)J_\mathcal{T}^{-\intercal}(x)\bigg]
    \nabla w(x;\mu) \cdot \nabla \tilde{\varphi}(x)
    -
    \tilde{f}(x; \mu) \, \tilde{\varphi}(x)
    \bigg)
    \, \mathrm{dx} \, \mathrm{d}\upmu
    = 0.
\end{equation}
To that end, we generate a set of $10^3$ functions $\varphi \in H^1_0(\mathbb B^2)$,
defined by $\varphi(x_1, x_2) = \alpha(x_1, x_2) P(x_1, x_2)$,
where $\alpha$ is the previously defined function that vanishes on $\mathbb S^1$,
and where $P$ is a two-dimensional polynomial of degree $2$
whose six coefficients are uniformly sampled in $(0, 1)$.
Then, the integral in \eqref{eq:fv_diff_forme_for_stats} is approximated
using $20\,000$ collocation points for $x_1$, $x_2$ and $\mu$.
The results are reported in \cref{tab:stats_PINN},
where we observe that, on average, the error is of the order of $10^{-2}$.
Therefore, our approach indeed provides a good approximation of the solution
to the parametric Poisson problem on a complex shape.

\begin{table}[!ht]
    \centering
    \caption{%
        Statistics on the absolute value of the integral
        in \eqref{eq:fv_diff_forme_for_stats}
        for the approximation of the parametric Poisson problem
        in \cref{sec:poisson_with_PINN}.
    }
    \label{tab:stats_PINN}
    \begin{tabular}{cccc}
        \toprule
        Mean value            & Maximal value
                              & Minimal value         & Standard deviation    \\
        \cmidrule(lr){1-4}
        $8.43 \times 10^{-3}$ & $4.53 \times 10^{-2}$
                              & $1.10 \times 10^{-5}$ & $7.08 \times 10^{-3}$ \\
        \bottomrule
    \end{tabular}
\end{table}

\subsection{Learning a simply connected parameterized shape with a SympNet}
\label{sec:sympnet_numerics}

The goal of this section is to show that our parametric SympNet
introduced in \cref{sec:sympnets}
is indeed able to learn a given symplectic map
as a deformation of the unit sphere of $\mathbb{R}^2$.
To that end, we train a SympNet to learn the shape given by the symplectic map
\eqref{eq:bizaroid}, with $\lambda \in \mathbb{M} = (0.5, 2)$.
Therefore, the SympNet takes as input the parameter $\lambda$,
as well as the two space variables $x_1$ and $x_2$.
To measure the error made by the SympNet,
we use the Hausdorff distance, computed with the algorithm from~\cite{taha2015efficient},
between the true shape and the learned shape.

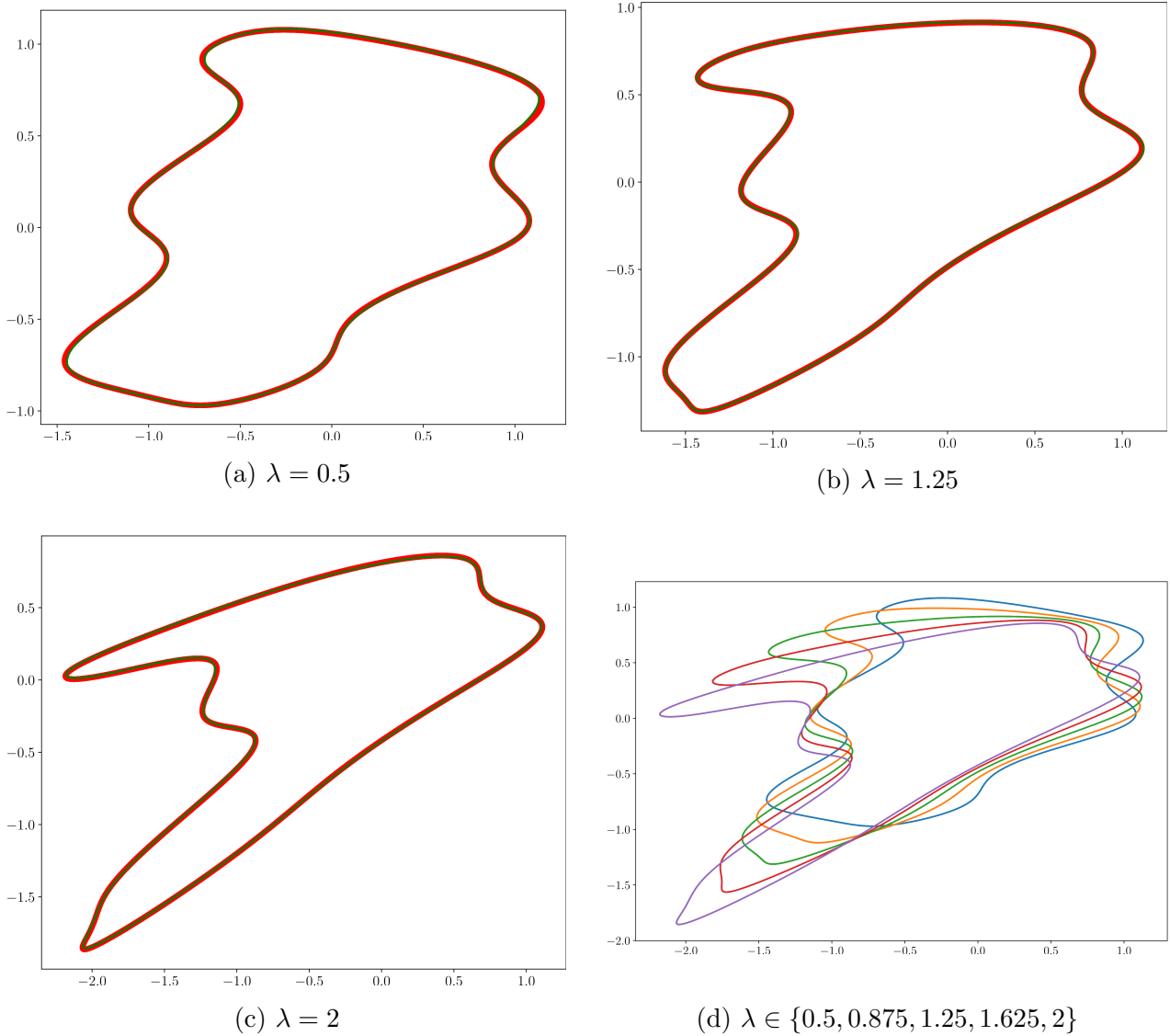
\begin{figure}[!ht]
    \centering
    \begin{tikzpicture}
        \node[xshift=-0.25\textwidth] (top_left) at (0,0)
        {\includegraphics[width=0.45\textwidth]{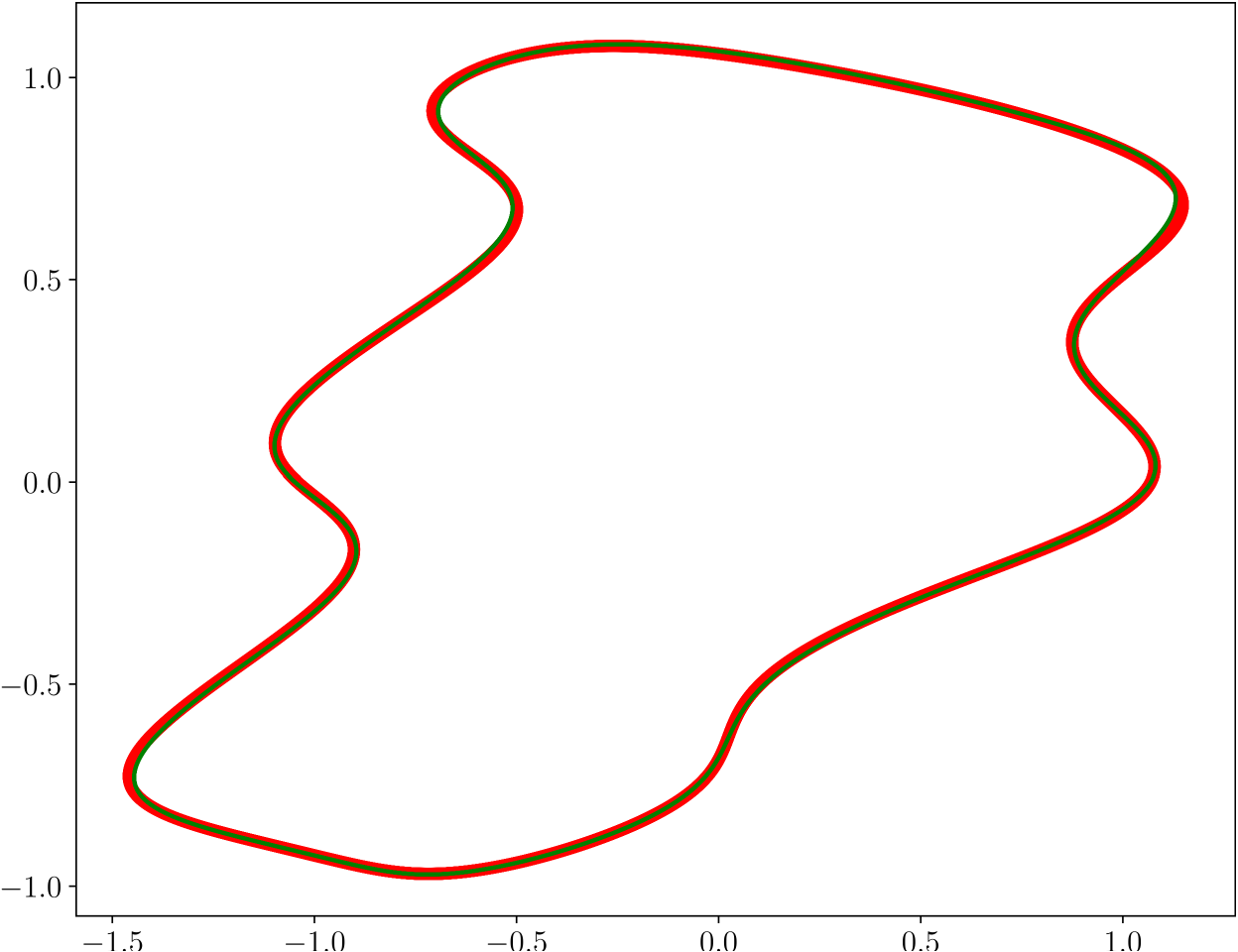}};
        \node[yshift=-0.02\textwidth] at (top_left.south) {(a) $\lambda=0.5$};

        \node[xshift=0.25\textwidth] (top_right) at (top_left.east)
        {\includegraphics[width=0.45\textwidth]{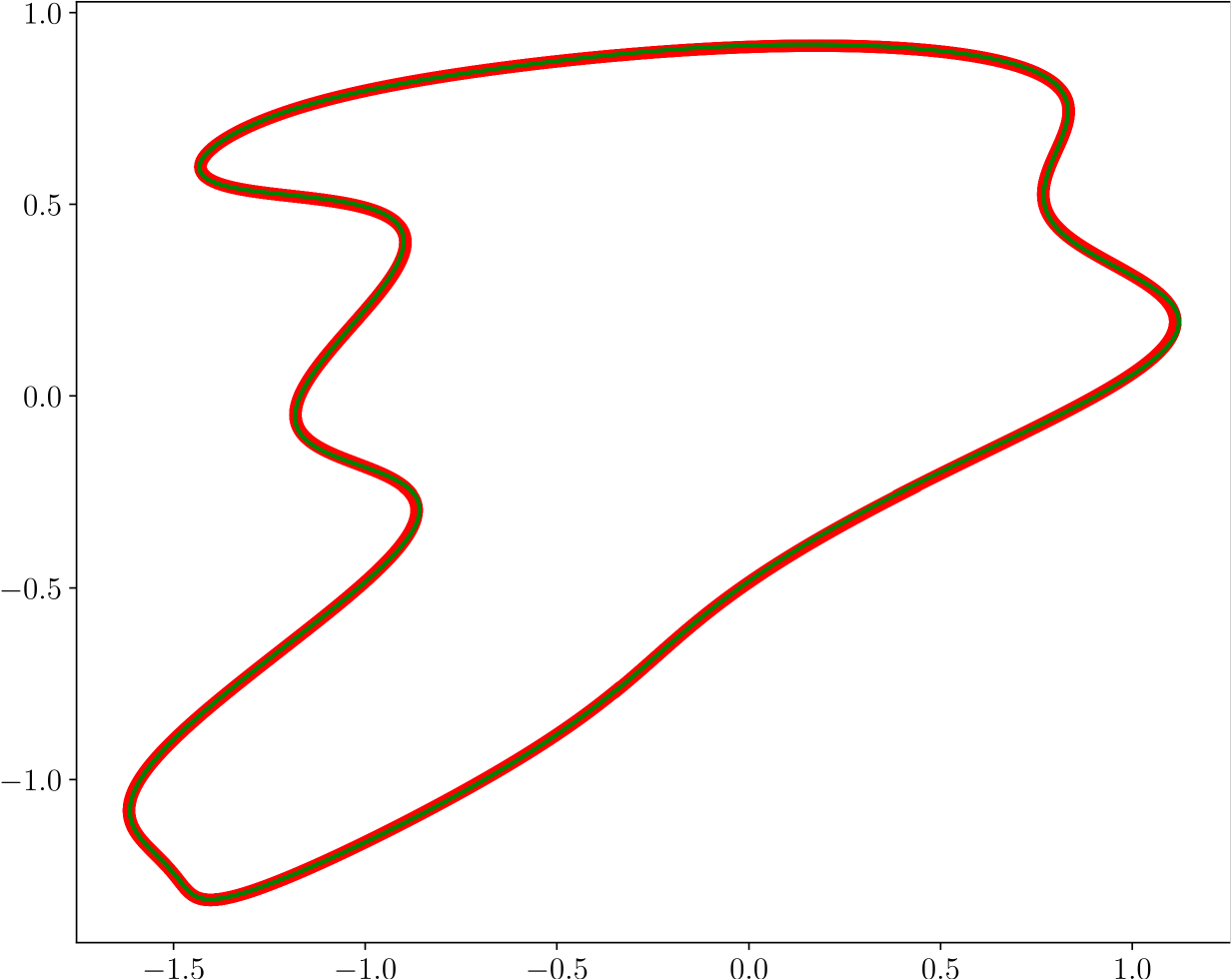}};
        \node[yshift=-0.02\textwidth] at (top_right.south) {(b) $\lambda=1.25$};

        \node[yshift=-0.25\textwidth] (bottom_left) at (top_left.south)
        {\includegraphics[width=0.45\textwidth]{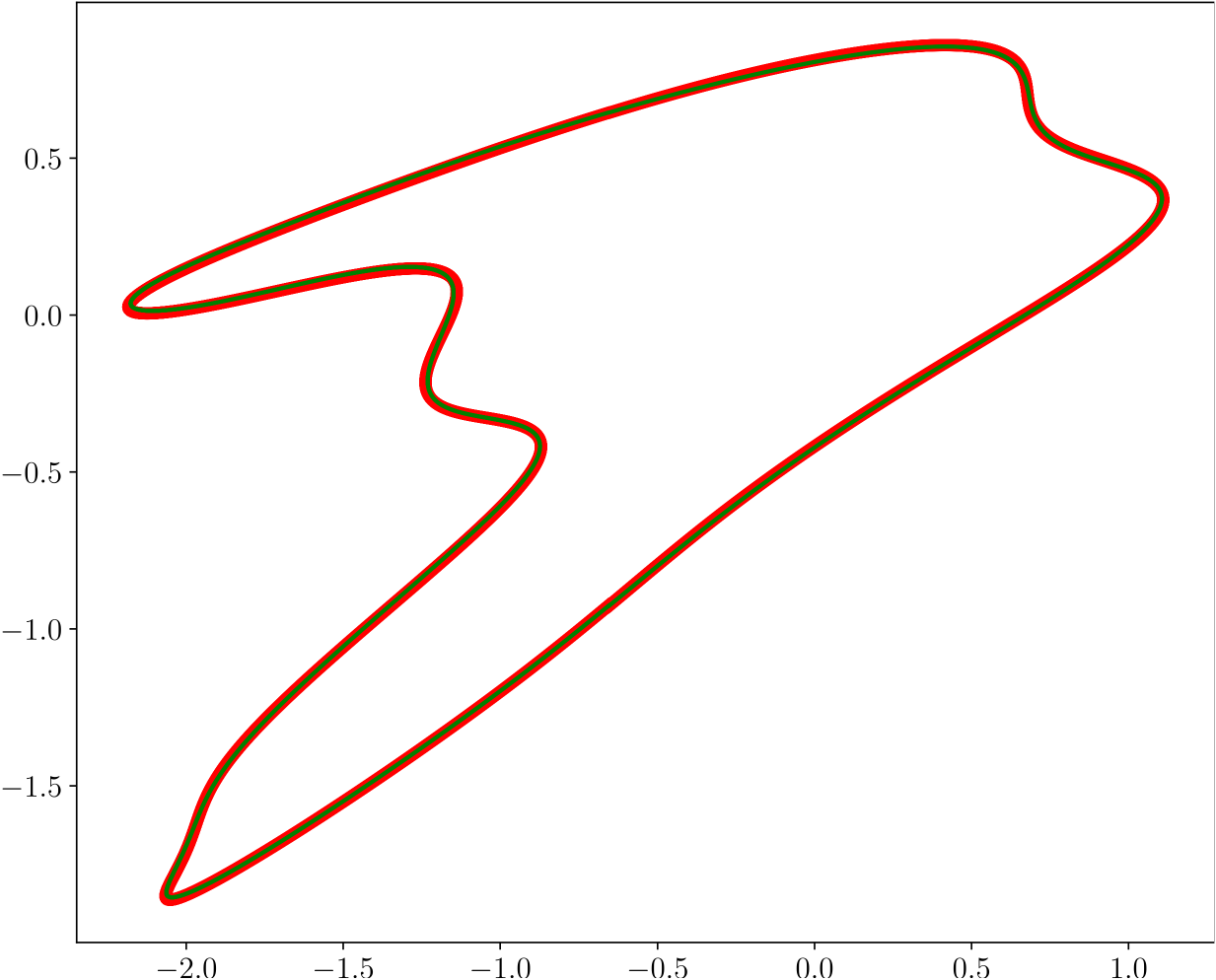}};
        \node[yshift=-0.02\textwidth] at (bottom_left.south) {(c) $\lambda=2$};

        \node[yshift=-0.25\textwidth] (bottom_right) at (top_right.south)
        {\includegraphics[width=0.45\textwidth]{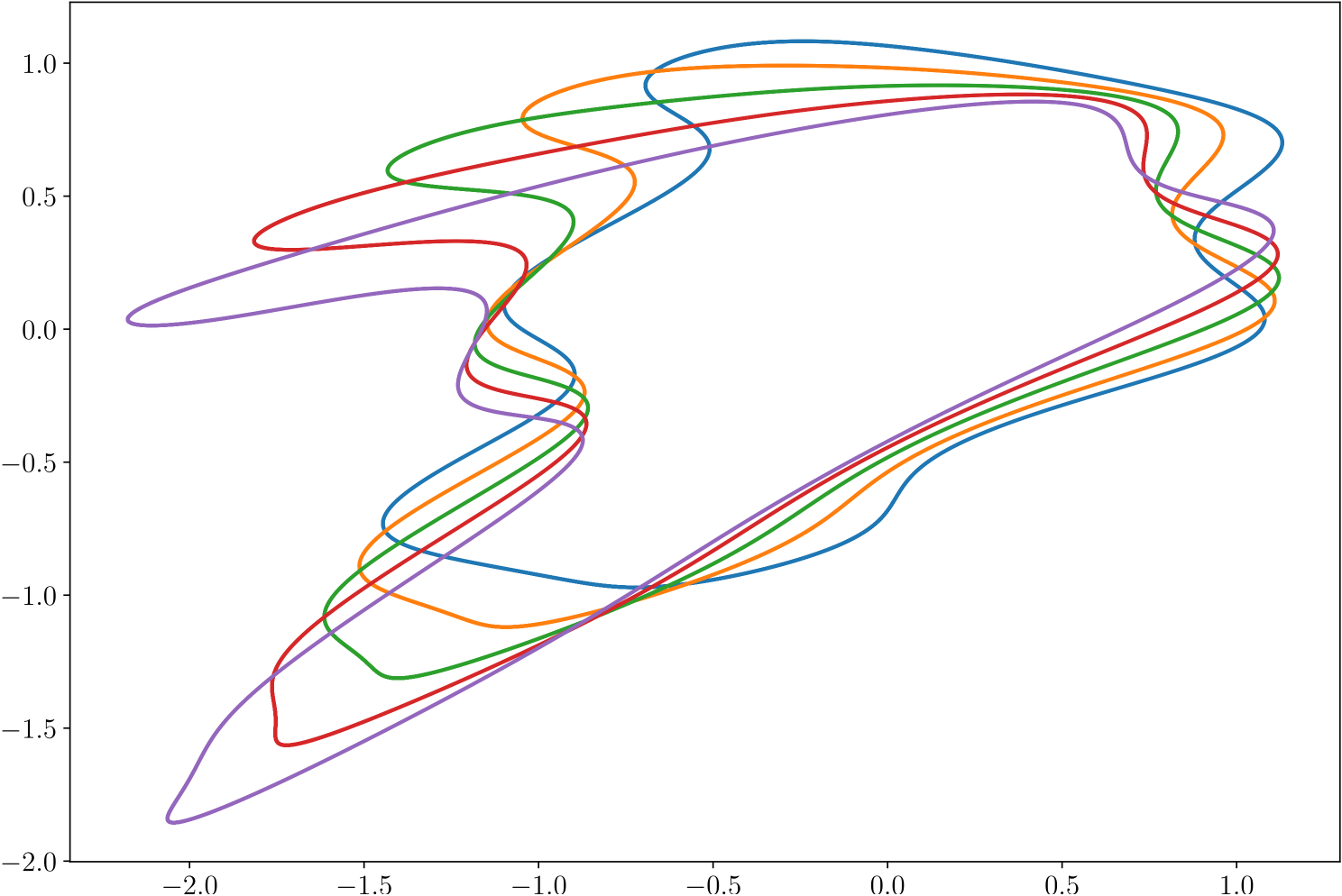}};
        \node[yshift=-0.045\textwidth] at (bottom_right.south) {(d) $\lambda \in \{0.5, 0.875,1.25, 1.625,2\}$};
    \end{tikzpicture}
    \caption{%
        SympNets for $\mathcal{T}_\lambda$,
        for $\lambda=0.5$ (upper left panel),
        $\lambda=1.25$ (upper right panel)
        and $\lambda=2$ (lower left panel).
        The red lines correspond to the reference shape,
        while the green ones correspond to the learned shape.
        The lower right picture is a superposition of learned shapes
        for $\lambda\in \{0.5, 0.875,1.25, 1.625,2\}$.%
    }
    \label{fig:sympnet_approximation}
\end{figure}

The results are displayed on \cref{fig:sympnet_approximation}.
We observe a good agreement between the reference shape and the learned shape,
for all considered values of $\lambda$.
Hence, our parametric SympNet is able to learn this rather complex parametric symplectic map.
    {\ra%
        We emphasize that, while \cref{fig:sympnet_approximation} features results with multiple parameters, we are interested in the quality of the approximation across the entire range $\lambda \in \mathbb{M} = (0.5, 2)$, rather than for any particular value of $\lambda$.
    }


To quantify the quality of the approximation,
we compute the Hausdorff distance between the reference shape and the learned shape.
To that end, we randomly sample $1000$ points in the parameter space $\mathbb M = (0.5, 2)$,
and compute some statistics on the Hausdorff distance.
They are collected in \cref{tab:stats_SympNet},
where we observe that the mean value of the Hausdorff distance is of the order of $10^{-2}$,
which makes up for a satisfactory approximation.

\begin{table}[!ht]
    \centering
    \caption{%
        Statistics on the Hausdorff distance between the reference shape and the learned shape,
        on the parameter set $\mathbb M = (0.5, 2)$.%
    }
    \label{tab:stats_SympNet}
    \begin{tabular}{cccc}
        \toprule
        Mean value            & Maximal value
                              & Minimal value         & Standard deviation    \\
        \cmidrule(lr){1-4}
        $9.45 \times 10^{-3}$ & $2.10 \times 10^{-2}$
                              & $6.05 \times 10^{-3}$ & $2.12 \times 10^{-3}$ \\
        \bottomrule
    \end{tabular}
\end{table}

\subsection{Solving the geometric Dirichlet energy problem with PINNs and SympNets}
\label{sec:gesonn_numerics}

The goal of this section is to show that our joint optimization procedure, involving both PINNs and SympNets
and based on the loss function \eqref{eq:loss_optim_forme},
is able to learn the optimal shape for the
parametric Poisson problem \eqref{eq:parametric_poisson}.
To that end, we first consider non-parametric problems to validate our approach,
before moving on to parametric ones.
According to \cref{thm:caracterisation_formes_optimales}, we will look at the first order optimality conditions and for a given shape $\Upomega$ and $u$ the associated solution of the PDE, we will consider the optimality error, an optimality criterion involving the standard deviation of \eqref{cion:surdet} on $\partial \Upomega$
\begin{equation*}
    \text{optimality error}:=\sqrt{\frac{1}{|\partial\Upomega|}\int_{\partial\Upomega}
        \left(\partial_n u(x) - \overline{\partial_n u}\right)^{\!2} \,
        \mathrm{d}\upsigma} = 0,
\end{equation*}
where $\overline{\partial_n u}$ denotes the average value
of $\partial_n u$ on $\partial\Upomega$.

\subsubsection{\texorpdfstring{Non-parametric problem: results with $f=1$}{Non-parametric problem: results with f=1}}
\label{sec:shapo_f_1}

\begin{figure}[!ht]
    \centering
    \begin{tikzpicture}
        \node[xshift=-0.25\textwidth] (top_left) at (0,0)
        {\includegraphics[width=0.352\textwidth]{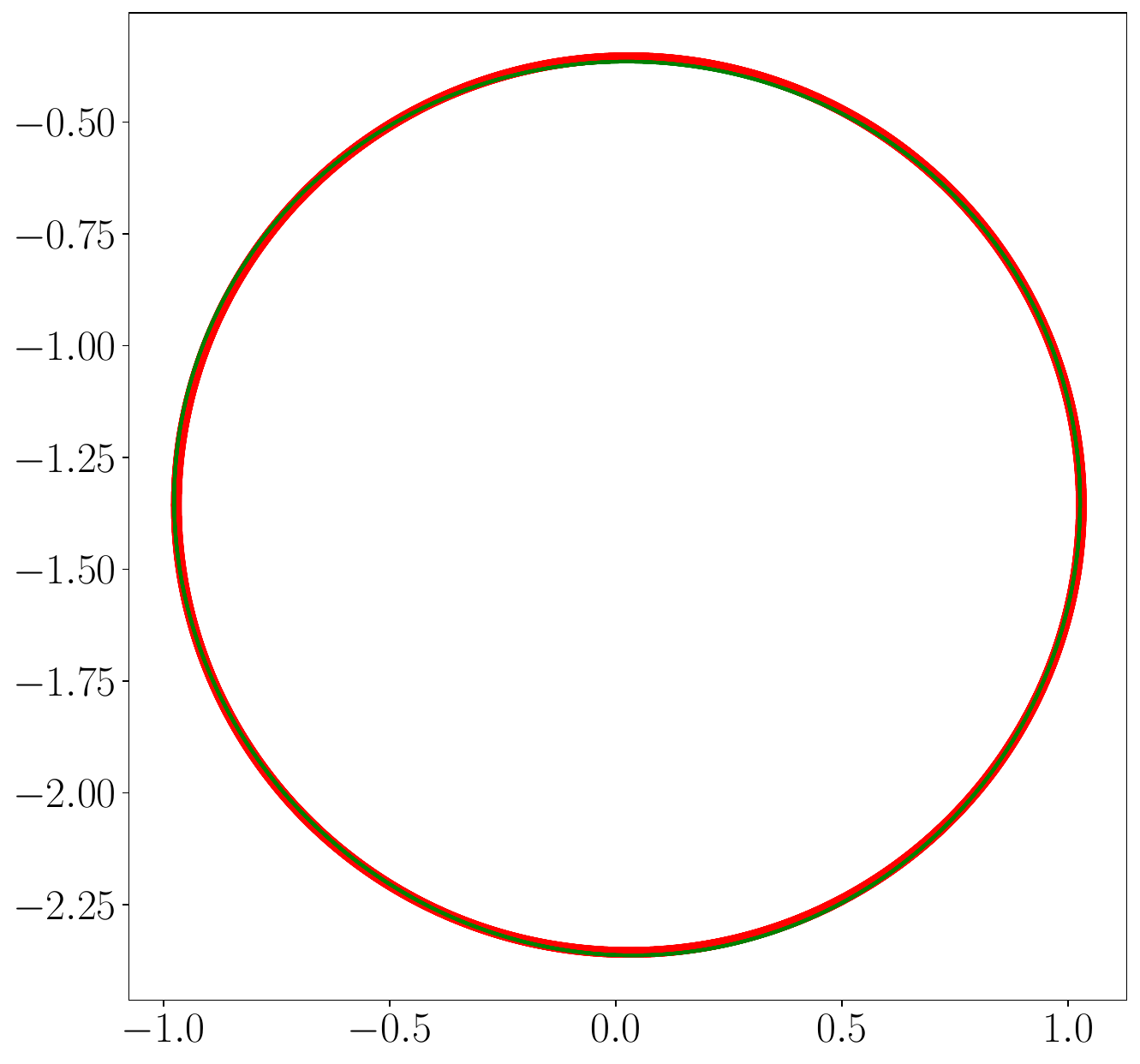}};

        \node[xshift=0.34\textwidth] (top_right) at (top_left.east)
        {\includegraphics[width=0.45\textwidth]{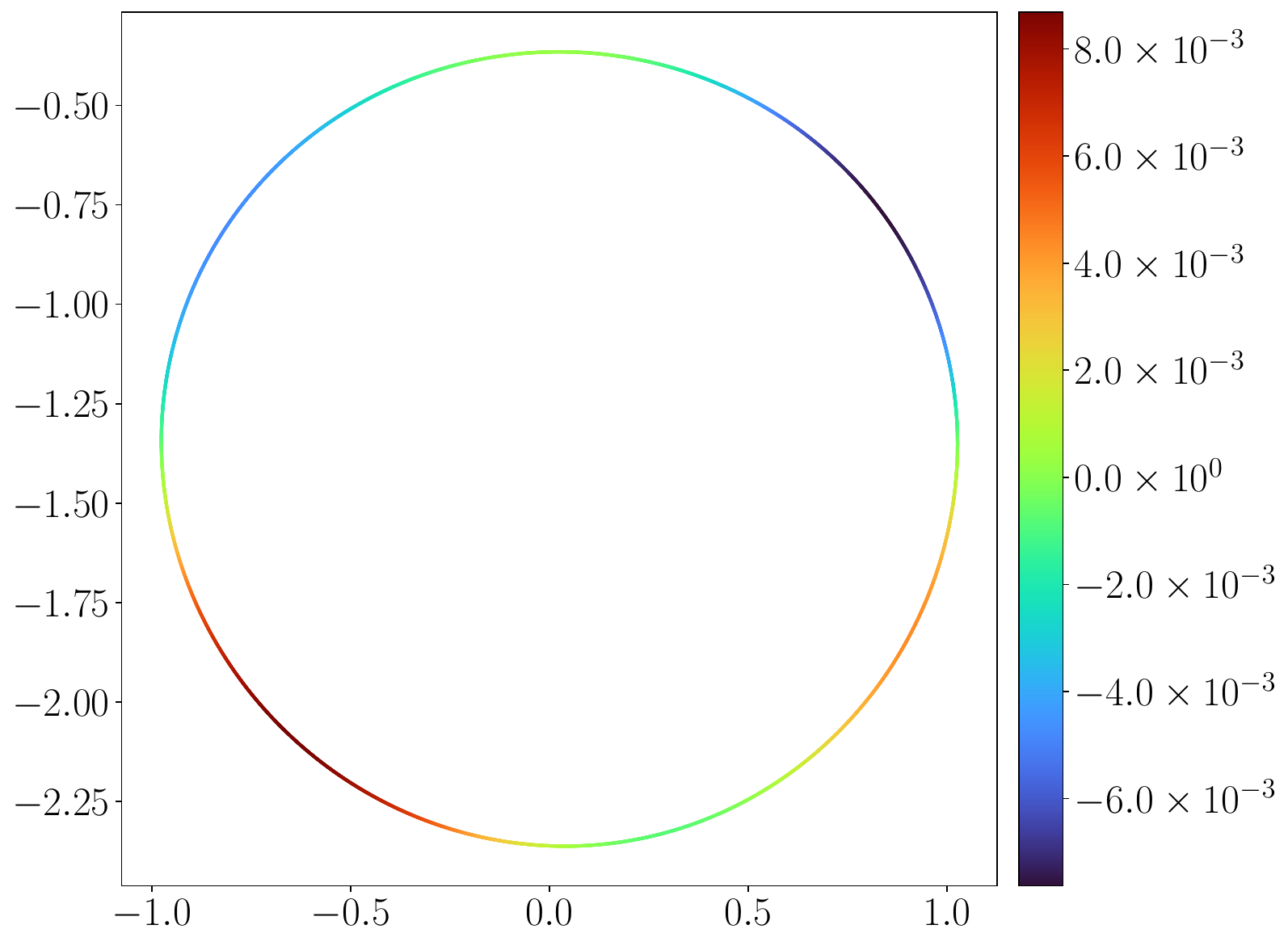}};

        \node[yshift=-0.2\textwidth] (bottom_right) at (top_right.south)
        {\includegraphics[width=0.45\textwidth]{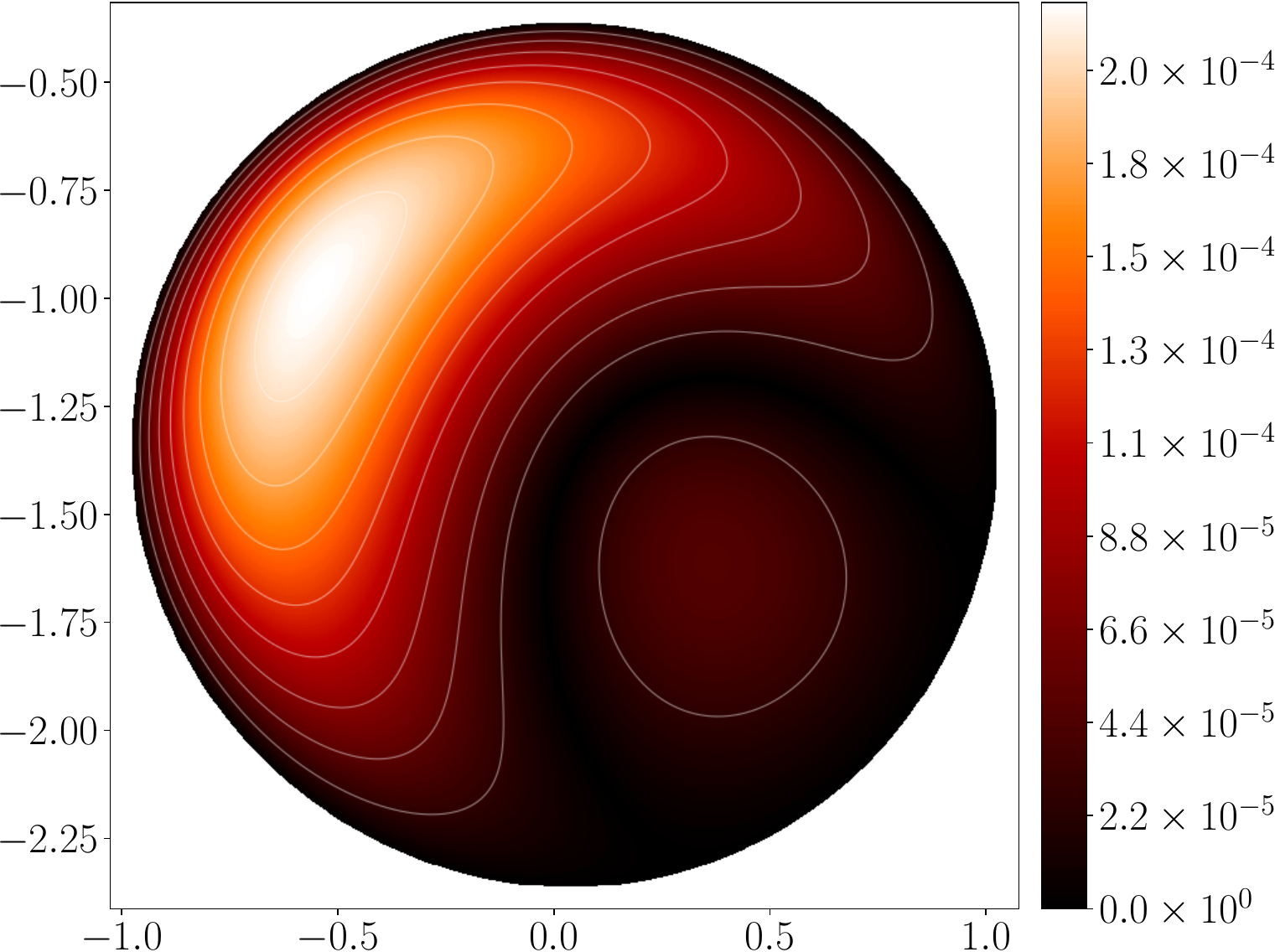}};

        \node[xshift=-0.25\textwidth] (bottom_left) at (bottom_right.west)
        {\includegraphics[width=0.45\textwidth]{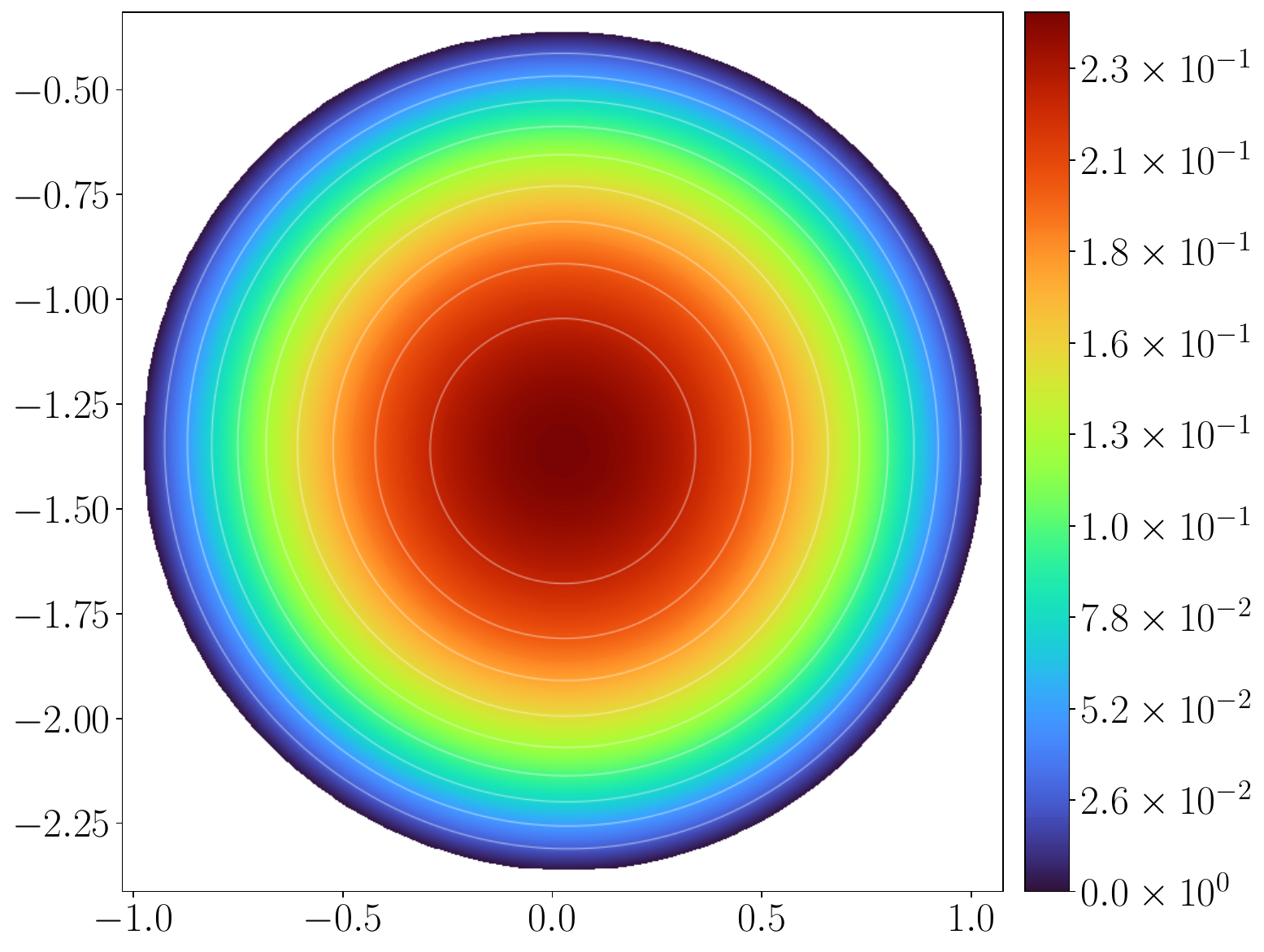}};

    \end{tikzpicture}

    \caption{%
        Shape optimization for the Dirichlet energy,
        for Poisson's equation with source term $f=1$.
        Top left panel: learned shape (green line)
        and reference shape (red line).
        Top right panel: deviation from the average of the optimality condition in the learned shape.
        Bottom left panel: approximate PDE solution.
        Bottom right panel: pointwise error
        between the true and approximate PDE solutions.
    }
    \label{fig:shapo_f_1}
\end{figure}

The results with $f=1$ are displayed on \cref{fig:shapo_f_1}. {\rb In this case, the optimal shape is known to be the unit sphere,
see~\cite{serrin1971symmetry}}.
However, we do not impose any constraint on the center of the unit sphere,
and therefore the optimal shape is not necessarily centered at the origin.
We observe a good agreement between the learned optimal shape
and the true one (top left panel),
as well as a good value of
the deviation from the average of the optimality condition (top right panel).
Moreover, the error between the approximate PDE solution
and the exact one is low (bottom right panel).
In addition, we report some metrics
in \cref{tab:stats_DeepShape_constant}
(namely, the Hausdorff distance between the learned and true shapes,
the optimality error,
and the~$L^2$ error between the exact and approximate solutions).
Remark that the Hausdorff distance reported in \cref{tab:stats_DeepShape_constant}
is of the same order as the one obtained in \cref{tab:stats_SympNet}
when using a SympNet to approximate a shape.
This is a good indicator that the SympNet in our approach
is fully working as expected.

\begin{table}[!ht]
    \centering
    \caption{%
        Relevant metrics related to the shape optimization
        of Poisson's equation with $f=1$.%
    }
    \label{tab:stats_DeepShape_constant}
    \begin{tabular}{ccc}
        \toprule
        Hausdorff distance
         & Optimality error
         & $L^2$ error           \\
        \cmidrule(lr){1-3}
        $3.09 \times 10^{-2}$
         & $4.34 \times 10^{-3}$
         & $9.15 \times 10^{-5}$ \\
        \bottomrule
    \end{tabular}
\end{table}

\begin{figure}[!ht]
    \centering
    \begin{tikzpicture}

        \node[xshift=-0.25\textwidth] (top_left) at (0,0)
        {\includegraphics[width=0.45\textwidth]{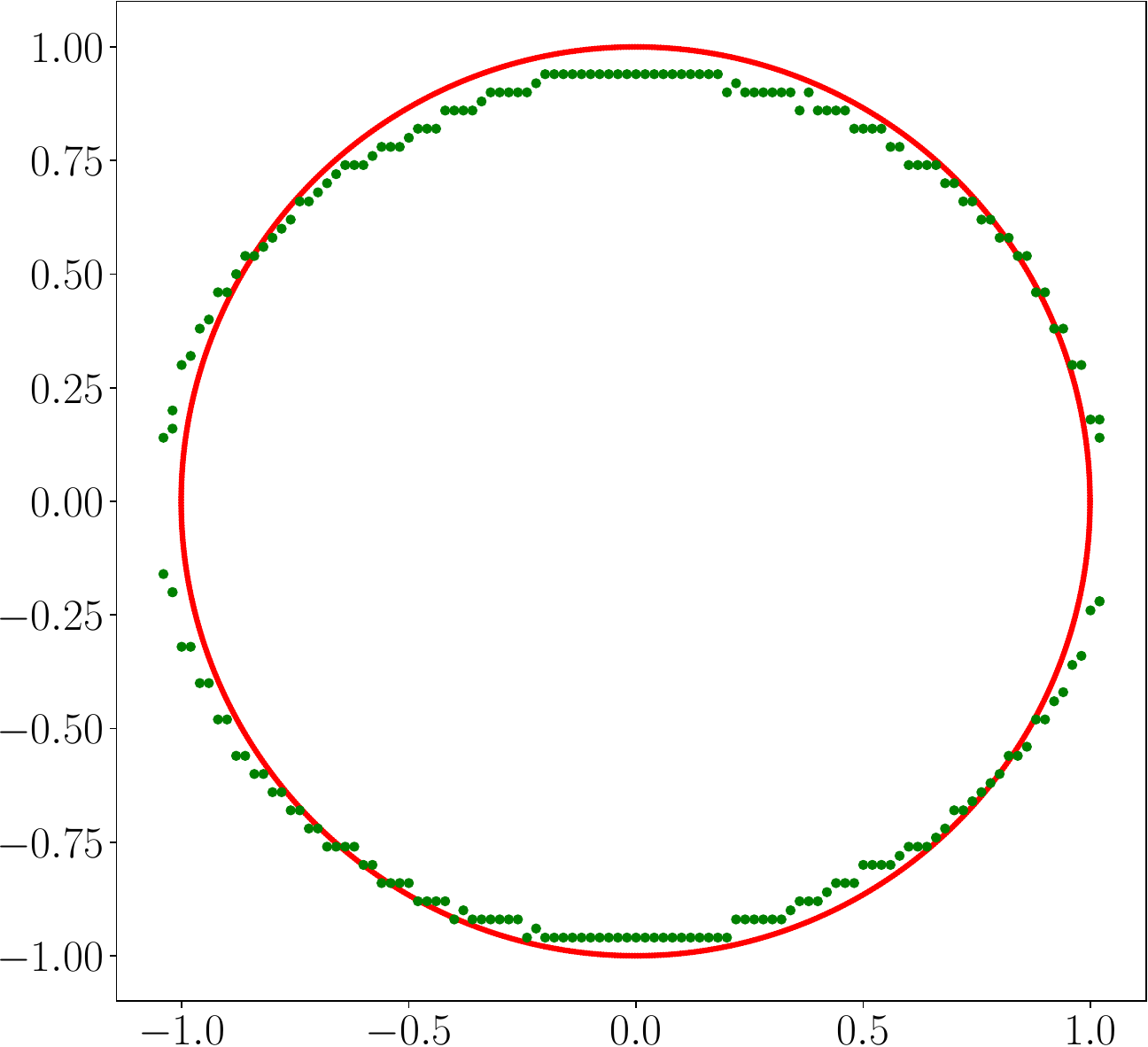}};
        \node[yshift=-0.02\textwidth] at (top_left.south) {(a) FEM, $R=100$};

        \node[xshift=0.25\textwidth] (top_right) at (top_left.east)
        {\includegraphics[width=0.45\textwidth]{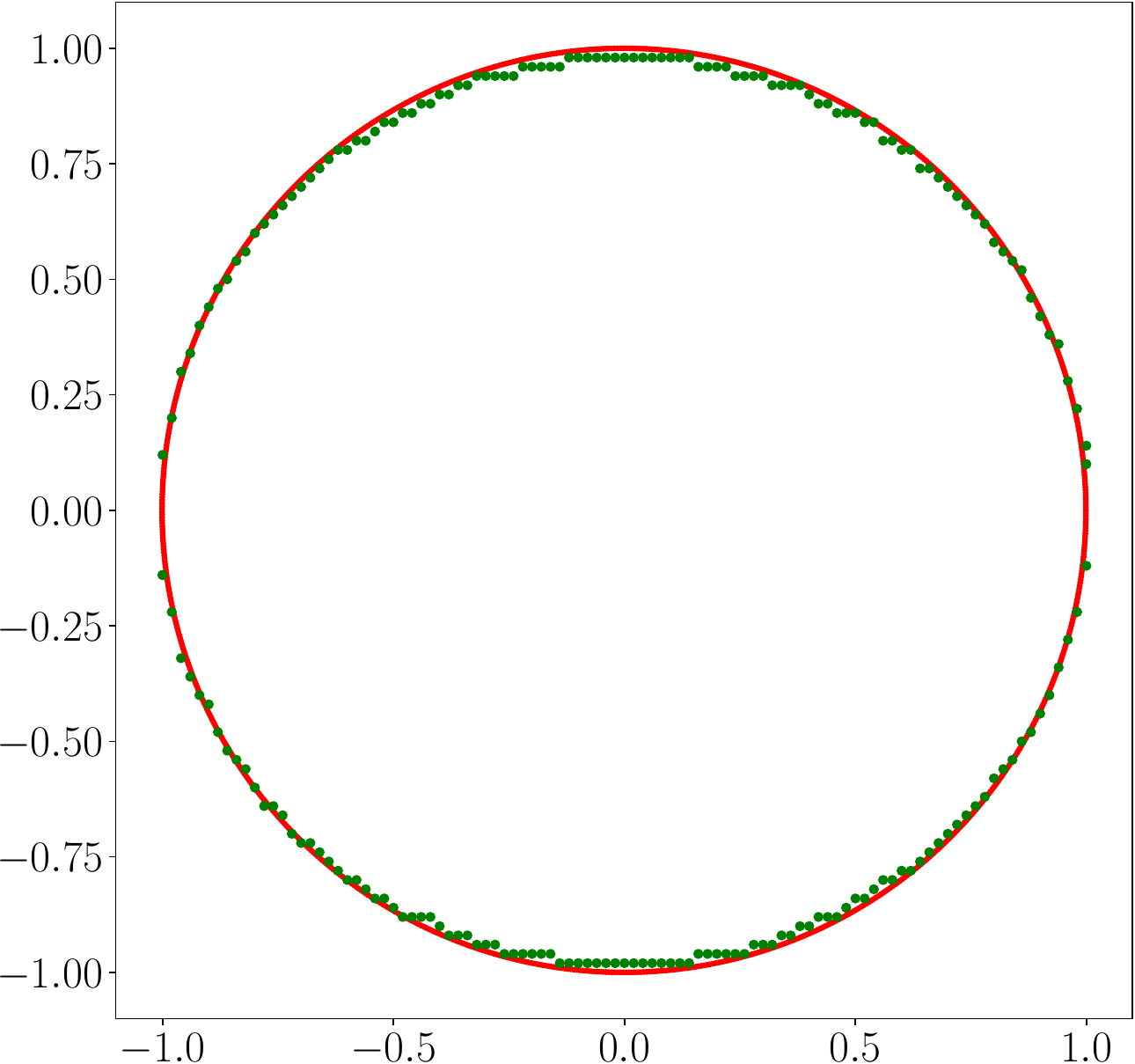}};
        \node[yshift=-0.02\textwidth] at (top_right.south) {(b) FEM, $R=250$};

        \node[yshift=-0.275\textwidth] (bottom_right) at (top_right.south)
        {\includegraphics[width=0.45\textwidth]{SIAM_one_backup_shape_error_error_cropped.pdf}};
        \node[yshift=-0.02\textwidth] at (bottom_right.south) {(d) GeSONN};

        \node[xshift=-0.25\textwidth] (bottom_left) at (bottom_right.west)
        {\includegraphics[width=0.45\textwidth]{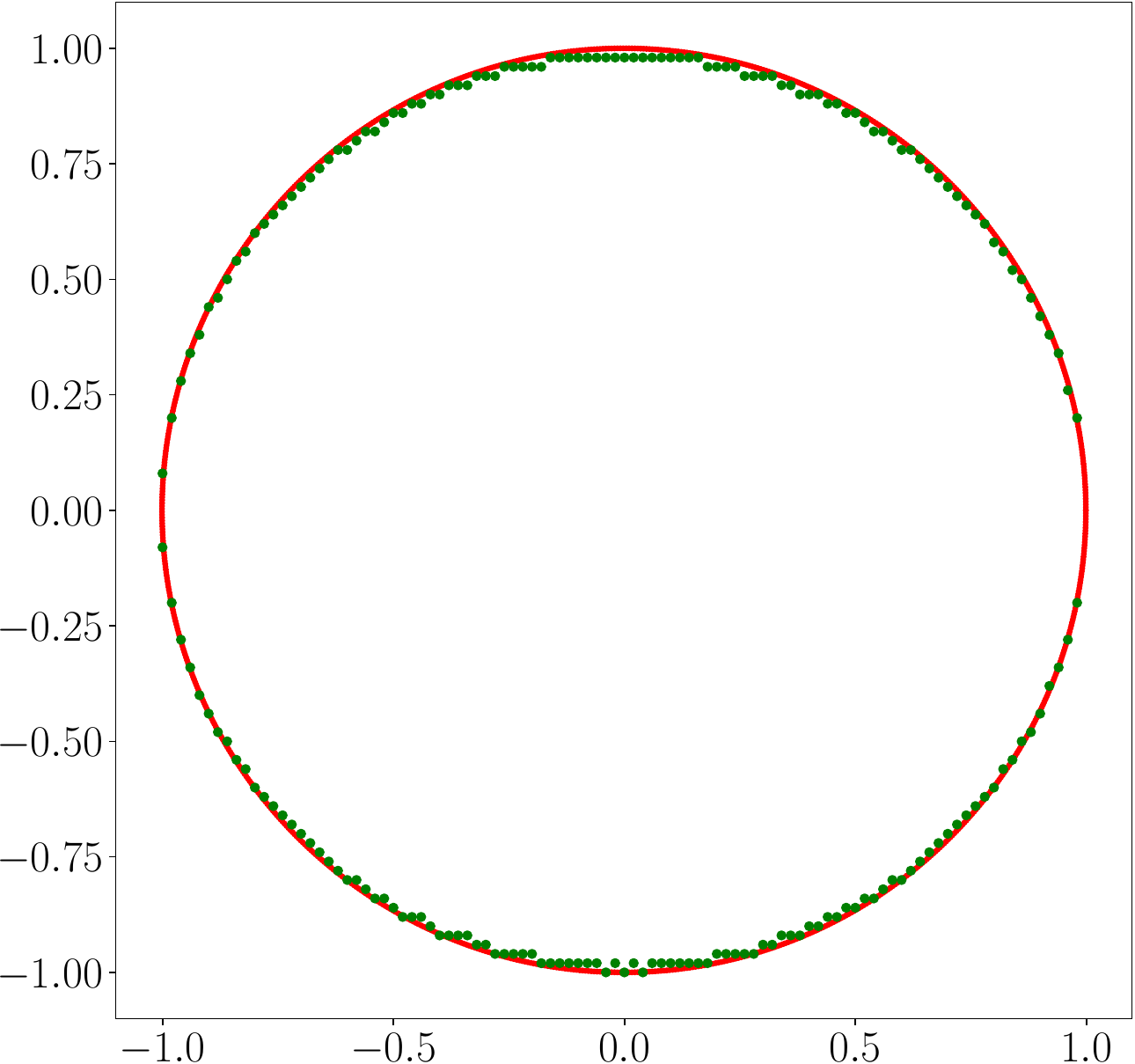}};
        \node[yshift=-0.02\textwidth] at (bottom_left.south) {(c) FEM, $R=500$};

    \end{tikzpicture}

    \caption{ \rb Shape optimization for the Dirichlet energy,
        for Poisson's equation with source term $f=1$.
        Sub-figures (a), (b) and (c) correspond to solutions obtained using the FEM-based procedure from \cite[Algorithm 1]{chambolle:hal-04140177}, with different mesh resolutions $R$ (each value of $R$ corresponds to a different number of points per edge; the total number of mesh elements is around $R^2$).
        Sub-figure (d) corresponds to the approximation with our method, see \cref{fig:shapo_f_1} for a depiction of other quantities.
        The exact solution is drawn in red, while the numerical approximations are drawn in green.
    }
    \label{fig:cv_algo_pt_fixe}
\end{figure}

{\rb%
\begin{remark}
   Since the problem under consideration has a known exact solution, we also compare the results of our method with those from an existing algorithm, whose convergence has been rigorously proven for a broad class of volume constraints in~\cite[Algorithm~1, Theorem~1]{chambolle:hal-04140177}. This algorithm is based on the finite element method (FEM), meaning that its computational time and accuracy are highly sensitive to the mesh resolution.
   We implemented a version of this algorithm using \texttt{FreeFem++}. For a meaningful comparison, both algorithms start with an initial shape defined as an ellipse with parameters $(2,0.5)$. While the computation time and accuracy of the FEM-based method depend on the mesh size, for our method they are influenced by the network parameters.
    The results are presented in \cref{tab:cv_algo_pt_fixe}, where the error is quantified as the $\ell^2$-norm between $N$ uniformly sampled points $(x_1^i, x_2^i)$ on the shape for our method, and the mesh points that define the optimal shape for the FEM-based method:
    \begin{equation*}
        \ell^2 \text{ error} \coloneqq \sqrt{\frac{1}{N}
            \sum_{i=1}^N \left( \left(x_1^i\right)^2 + \left(x_2^i\right)^2 - 1 \right)^2},
    \end{equation*}
    It is important to note that we did not compute the Hausdorff distance, as this would unfairly penalize the FEM-based method. The absence of sufficient points on the shape would artificially inflate this distance. In contrast, our method allows for the generation of as many points as needed on the shape, whereas the FEM-based method is limited to its predefined mesh points.
    In \cref{tab:cv_algo_pt_fixe},
    we observe that our method outperforms the FEM-based method in terms of both speed and accuracy;
    We also note that the error decreases with the mesh size for the FEM-based method.
    In \cref{fig:cv_algo_pt_fixe}, we present the optimal shapes obtained by both methods alongside the true optimal shape (a disk). This clearly demonstrates that our method produces an optimal shape that is closer to the true one.
    Finally, as explained in the next sections, our method is also capable of solving parametric shape optimization problems with a single training, incurring minimal additional computational cost. Given that the inference time is negligible (approximately $10^{-3}$ seconds), this represents a significant advantage of our approach.
\end{remark}%
}

\begin{table}[!ht]
    \centering
    \caption{%
        \rb	Comparison between the results of the FEM-based method
        \cite[Algorithm 1]{chambolle:hal-04140177} and our approach.
        For the FEM-based method, we report the results for three different
        mesh resolutions $R$, corresponding to the number of points per edge
        (roughly $R^2$ elements in the mesh).
        For fairness, both methods are run on the same CPU,
        even though GeSONN natively runs on GPUs thanks to the \texttt{pytorch} library.%
    }
    \label{tab:cv_algo_pt_fixe}
    \rb
    \begin{tabular}{ccccc}
        \toprule
        method
         & FEM ($R = 100$)
         & FEM ($R = 250$)
         & FEM ($R = 500$)
         & GeSONN                \\
        \cmidrule(lr){1-5}
        Computational time (s)
         & $53.3$
         & $509$
         & $3020$
         & $22.5$                \\
        $\ell^2 \text{ error}$
         & $7.20 \times 10^{-2}$
         & $2.83 \times 10^{-2}$
         & $2.21 \times 10^{-2}$
         & $1.99 \times 10^{-3}$ \\
        \bottomrule
    \end{tabular}
\end{table}

\subsubsection{\texorpdfstring{Non-parametric problem: results with $f$ given by \eqref{eq:parametric_source} with $\mu=0.5$}{Non-parametric problem: results with exponential source term}}

\begin{figure}[!ht]
    \centering
    \begin{tikzpicture}
        \node[xshift=-0.25\textwidth] (top_left) at (0,0)
        {\includegraphics[width=0.387\textwidth]{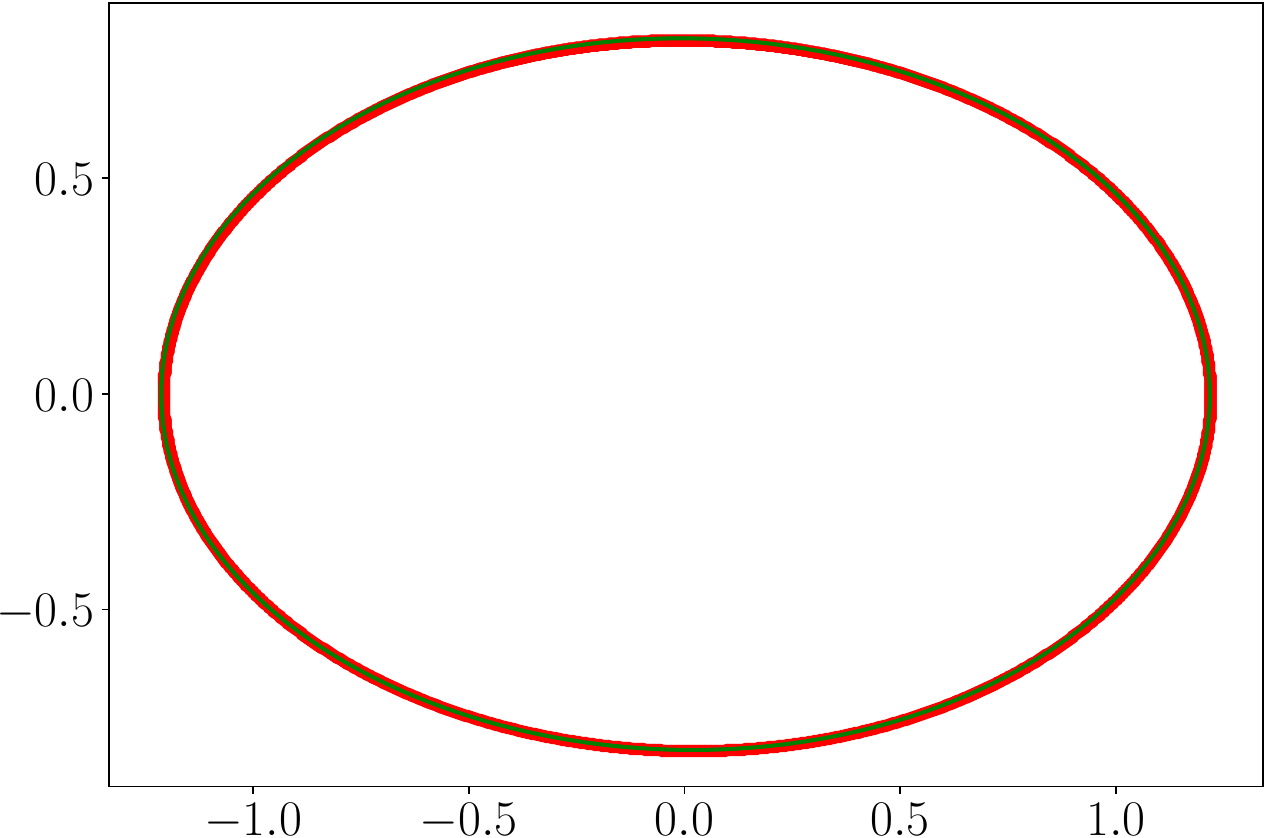}};

        \node[xshift=0.337\textwidth] (top_right) at (top_left.east)
        {\includegraphics[width=0.475\textwidth]{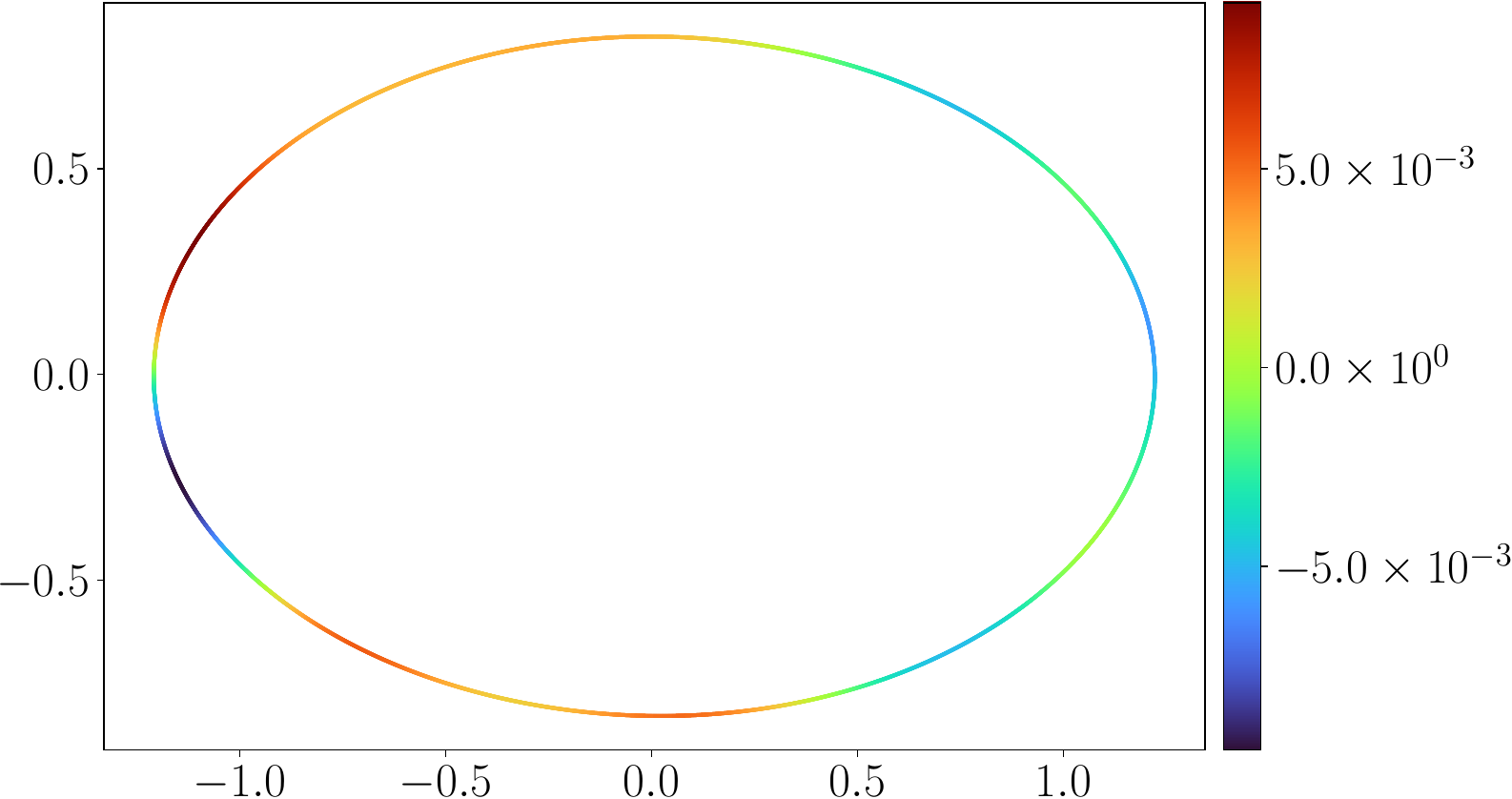}};

        \node[yshift=-0.15\textwidth] (bottom_right) at (top_right.south)
        {\includegraphics[width=0.475\textwidth]{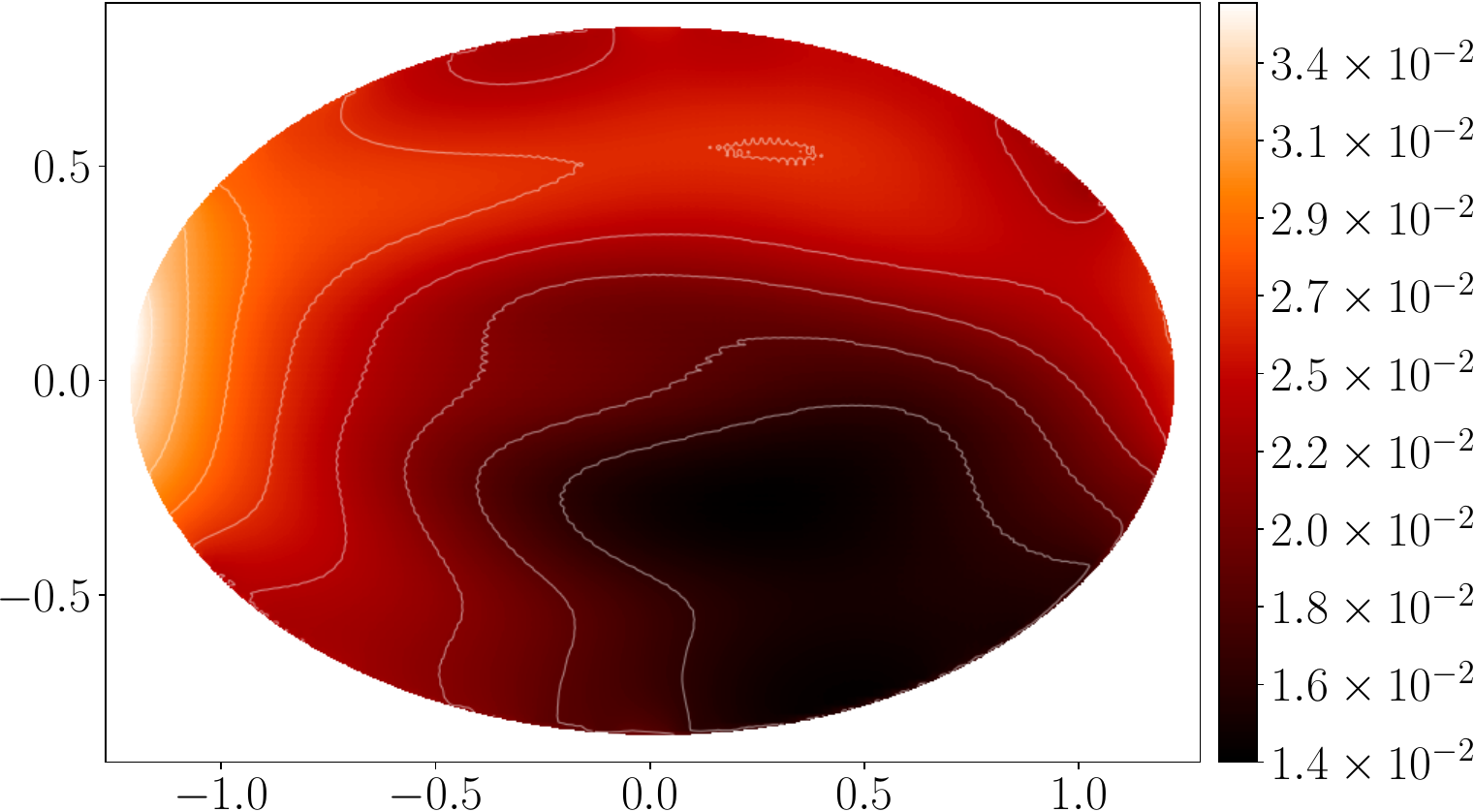}};

        \node[xshift=-0.25\textwidth] (bottom_left) at (bottom_right.west)
        {\includegraphics[width=0.475\textwidth]{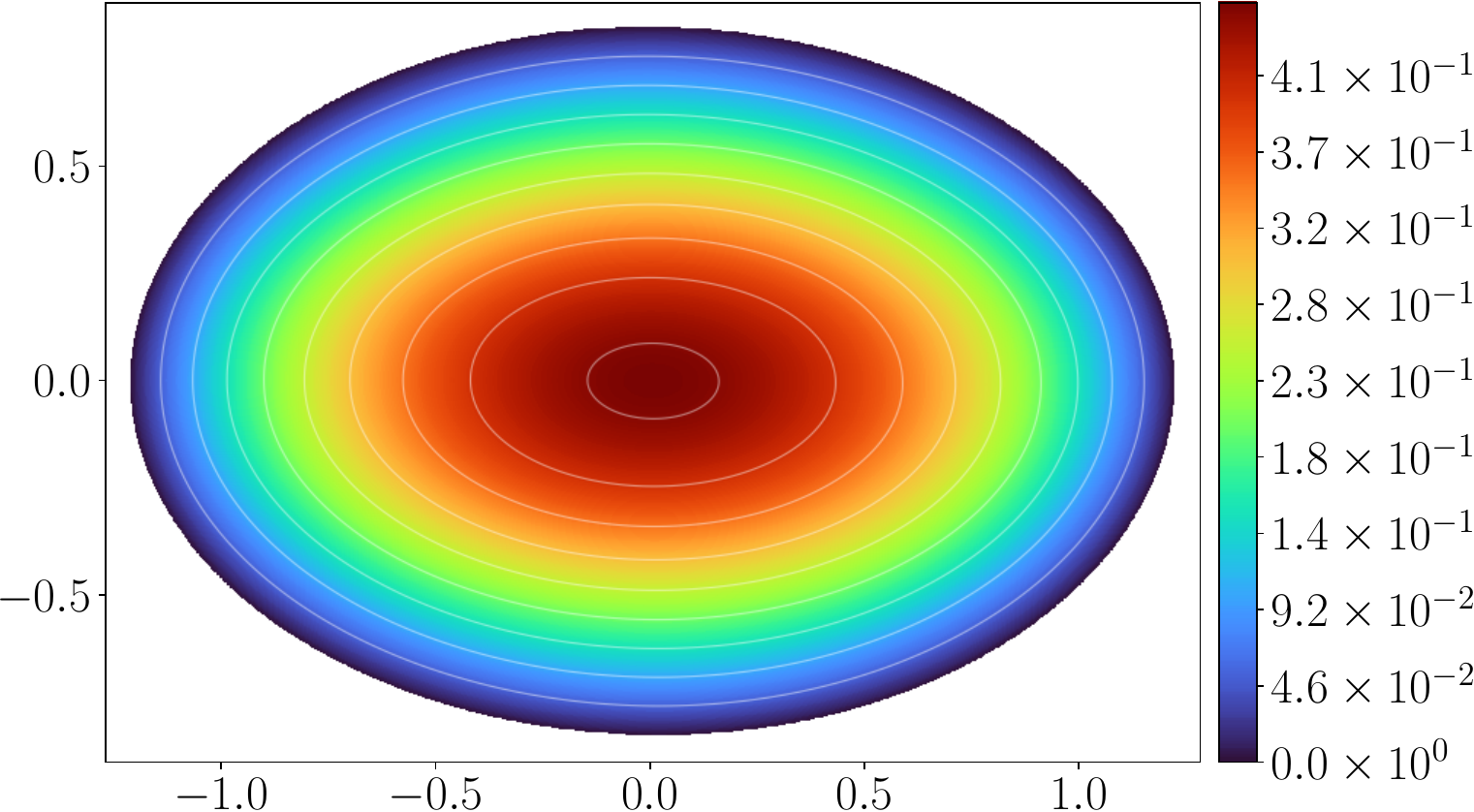}};

    \end{tikzpicture}

    \caption{%
        Shape optimization for the Dirichlet energy,
        for Poisson's equation with source term $f$
        given by \eqref{eq:parametric_source} with $\mu=0.5$. The results of a fixed-point algorithm \cite{chambolle:hal-04140177} is our reference for this numerical experiment.
        Top left panel: learned shape (green line)
        and reference shape (red line).
        Top right panel: deviation from the average of the optimality condition in the learned shape.
        Bottom left panel: approximate PDE solution.
        Bottom right panel: pointwise error
        between the reference and approximate PDE solutions.
    }
    \label{fig:shapo_f_exp}
\end{figure}

Next, we turn to a more complex, non-constant source term.
Namely, we take~$f$ given by \eqref{eq:parametric_source} with $\mu=0.5$,
and display the results on \cref{fig:shapo_f_exp}.
This time, the reference shape has been obtained by using \texttt{FreeFem++},
with a fixed point algorithm based on~\cite[Algorithm 1]{chambolle:hal-04140177}.
As a consequence, we only consider it as a point of reference rather than a ground truth.
We observe a good agreement between the learned shape
and the reference one, much like in \cref{sec:shapo_f_1}.
Moreover, to quantify the approximation quality,
we report in \cref{tab:stats_DeepShape_exp}
the same metrics as in \cref{sec:shapo_f_1}.

\begin{table}[!ht]
    \centering
    \caption{%
        Relevant metrics related to the shape optimization
        of Poisson's equation with source term $f$
        given by \eqref{eq:parametric_source} with $\mu=0.5$.%
    }
    \label{tab:stats_DeepShape_exp}
    \begin{tabular}{ccc}
        \toprule
        Hausdorff distance
         & Optimality error
         & $L^2$ error           \\
        \cmidrule(lr){1-3}
        $8.29 \times 10^{-3}$
         & $4.43 \times 10^{-3}$
         & $2.20 \times 10^{-2}$ \\
        \bottomrule
    \end{tabular}
\end{table}

\subsubsection{\texorpdfstring{Parametric problem: results with constant $f$}{Parametric problem: results with constant f}}
\label{sec:shapo_f_1_param}

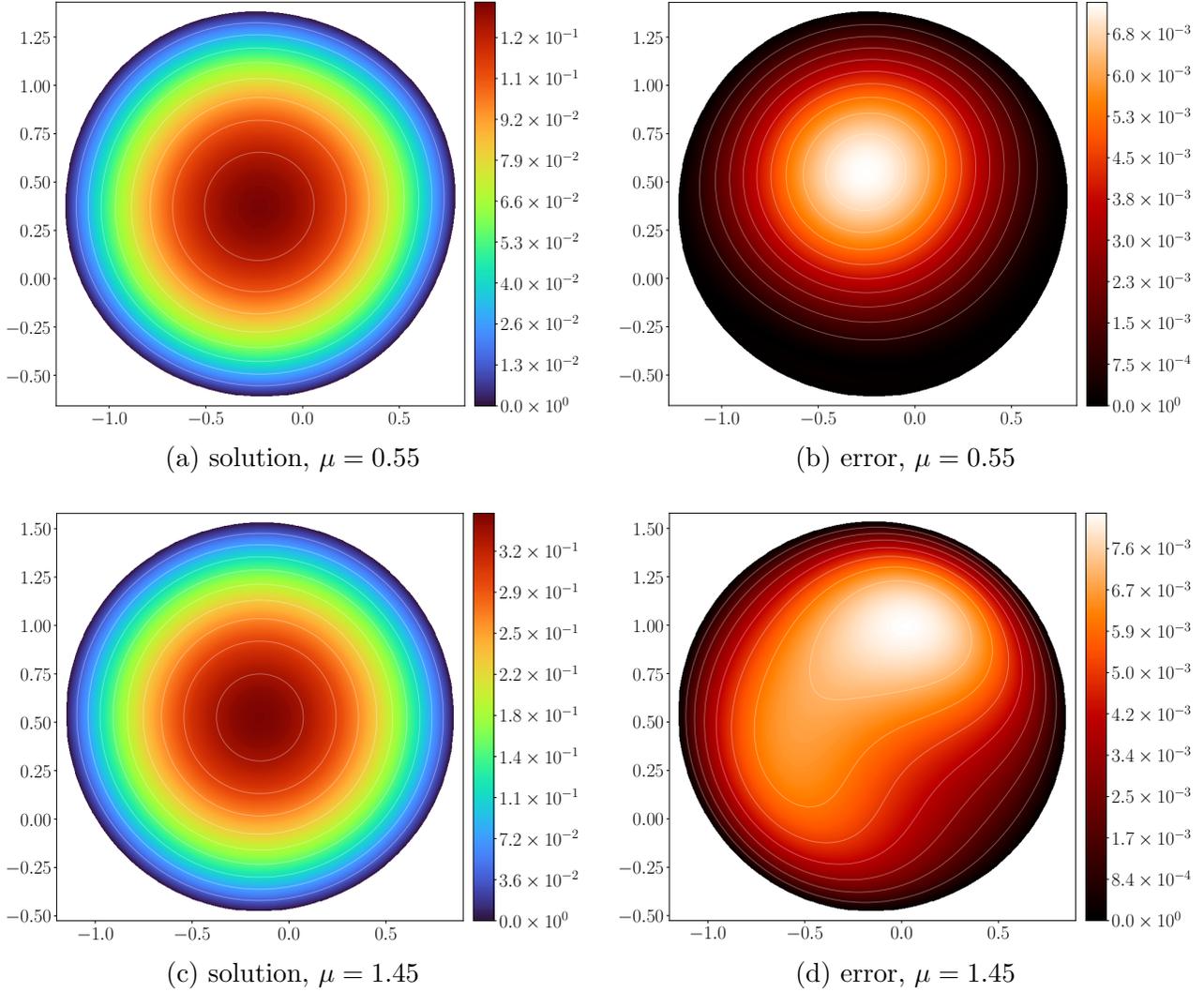
\begin{figure}[!ht]
    \centering
    \begin{tikzpicture}
        \node[xshift=-0.25\textwidth] (top_left) at (0,0)
        {\includegraphics[width=0.45\textwidth]{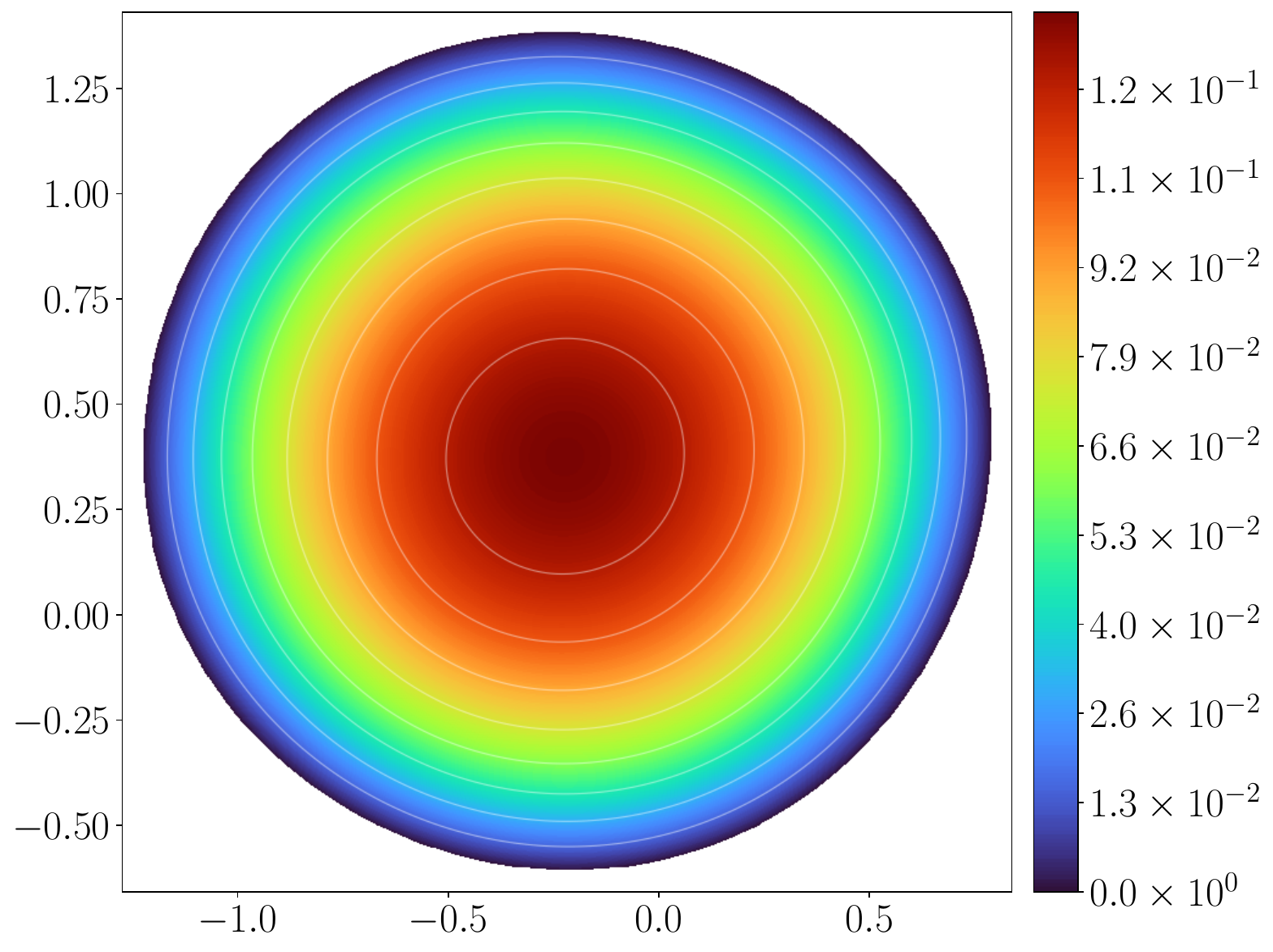}};
        \node[yshift=-0.02\textwidth] at (top_left.south) {(a) solution, $\mu=0.55$};

        \node[xshift=0.25\textwidth] (top_right) at (top_left.east)
        {\includegraphics[width=0.45\textwidth]{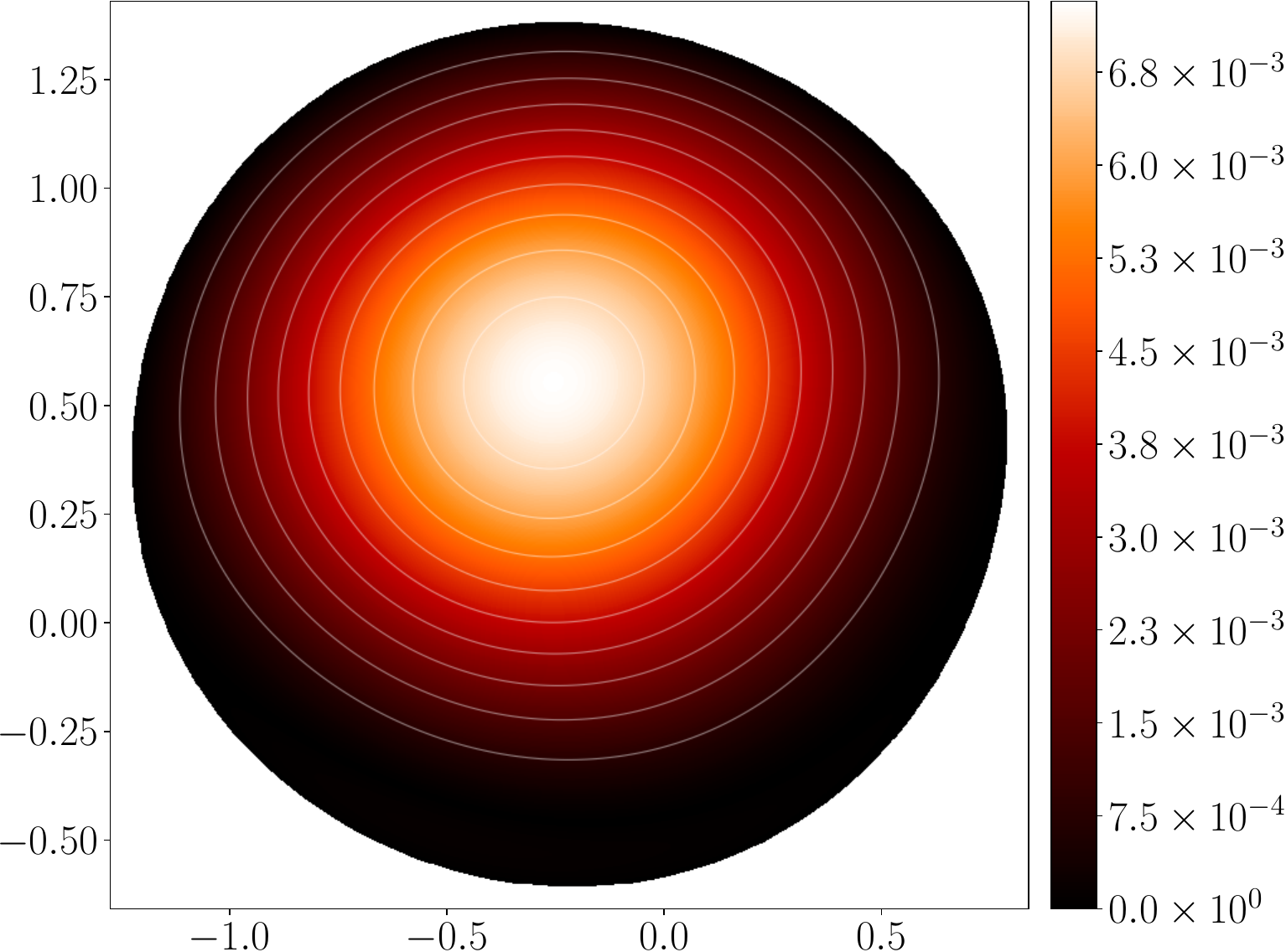}};
        \node[yshift=-0.02\textwidth] at (top_right.south) {(b) error, $\mu=0.55$};

        \node[yshift=-0.23\textwidth] (bottom_left) at (top_left.south)
        {\includegraphics[width=0.45\textwidth]{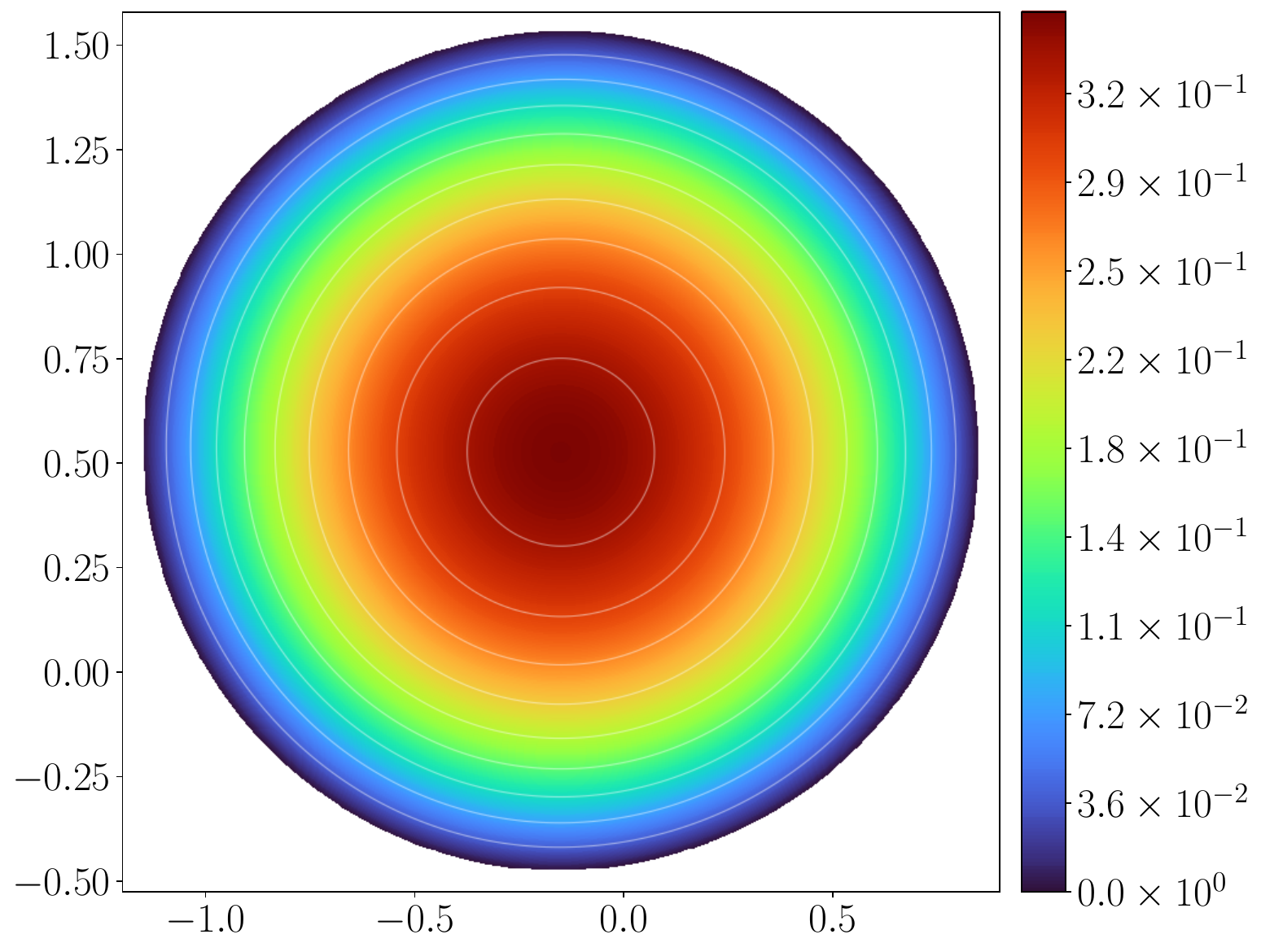}};
        \node[yshift=-0.02\textwidth] at (bottom_left.south) {(c) solution, $\mu=1.45$};

        \node[yshift=-0.23\textwidth] (bottom_right) at (top_right.south)
        {\includegraphics[width=0.45\textwidth]{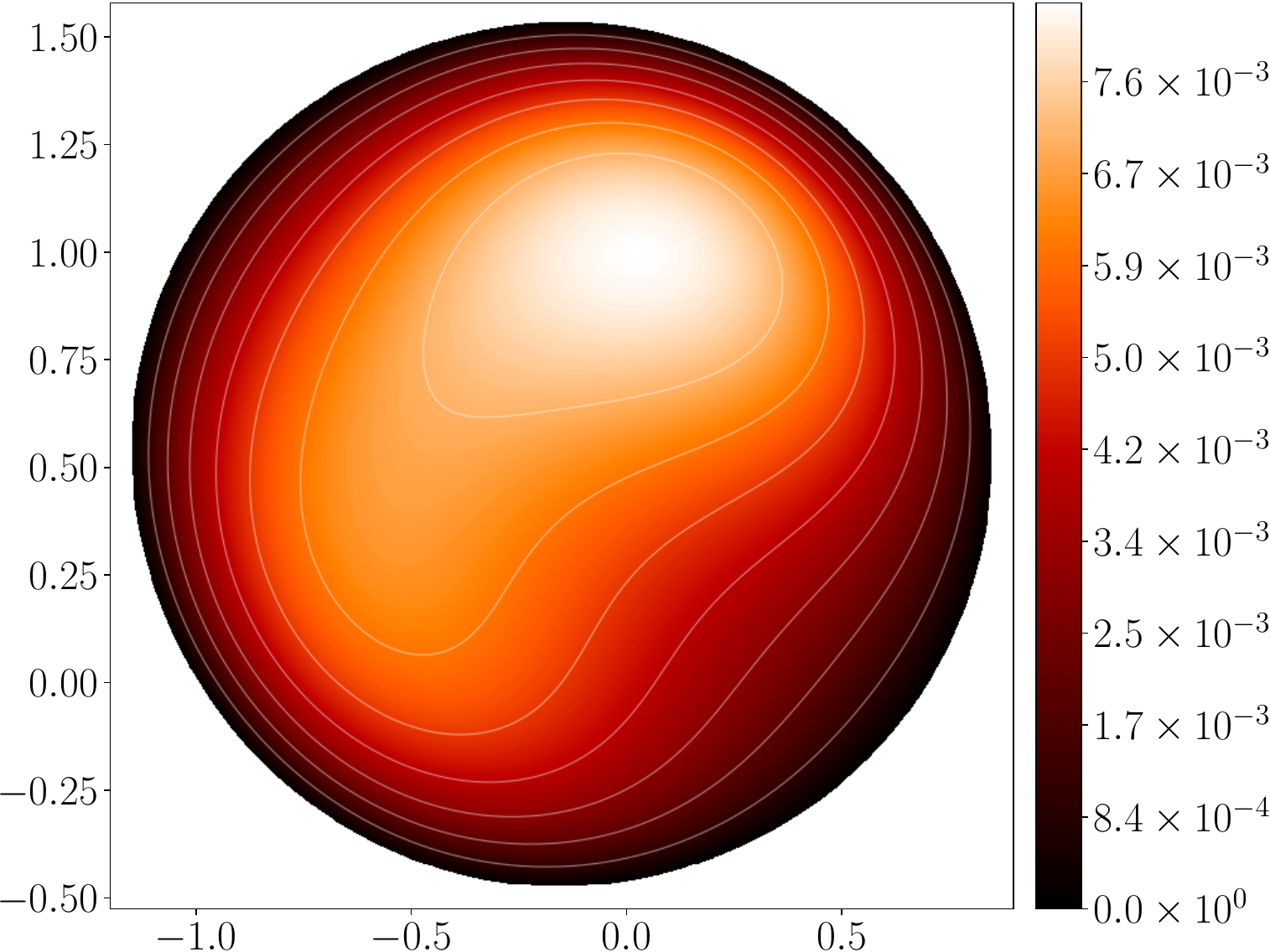}};
        \node[yshift=-0.02\textwidth] at (bottom_right.south) {(d) error, $\mu=1.45$};

    \end{tikzpicture}

    \caption{%
        Shape optimization for the Dirichlet energy,
        for Poisson's equation with source term $f=\mu$.
        For two values of $\mu$, we display the approximate solution (left column)
        and the pointwise error between the approximate and exact solutions (right column).
    }
    \label{fig:shapo_f_1_param}
\end{figure}

\begin{figure}[!ht]
    \centering

    \includegraphics[width=0.45\textwidth]{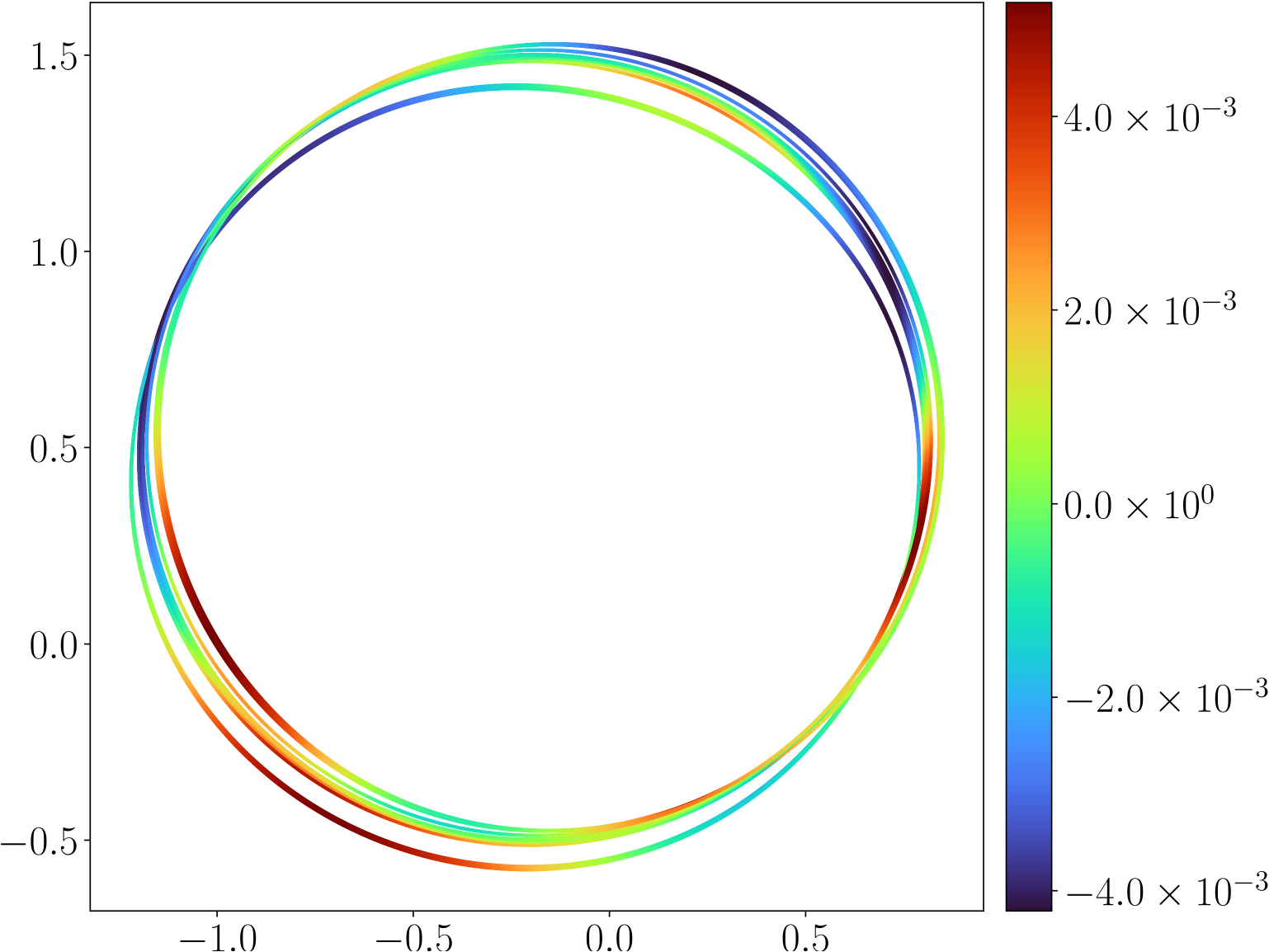}

    \caption{%
        Shape optimization for the Dirichlet energy,
        for Poisson's equation with source term $f=\mu$:
        optimal shapes and deviation from the average of the optimality condition for $10$ random values of~$\mu$.
    }
    \label{fig:shapo_f_1_param_optimality_conditions}
\end{figure}

We now turn to a parametric shape optimization problem,
where the source term is a constant $\mu \in \mathbb{M} = (0.5, 1.5)$.
The results and errors are displayed on \cref{fig:shapo_f_1_param},
where there is a good agreement between approximate and exact solution.
In addition, we display the optimal shapes and the optimality error for $10$ random values of $\mu$ on \cref{fig:shapo_f_1_param_optimality_conditions}.
We note that the first order optimality condition is indeed satisfied,
even though the shapes themselves are not necessarily centered at the origin.
Indeed, like in \cref{sec:shapo_f_1},
the optimal shape is the unit sphere, but the center of the sphere is not fixed.
Finally, we report some statistics on our three main metrics
in \cref{tab:stats_DeepShape_constant_param},
which confirm the relevance of our approach for a parametric problem.

\begin{table}[!ht]
    \centering
    \caption{%
        Statistics of relevant metrics in the case of
        a parametric constant source term $f=\mu$,
        obtained by computing each metric for $10^3$ values of $\mu$.%
    }
    \label{tab:stats_DeepShape_constant_param}
    \begin{tabular}{ccccc}
        \toprule
        Metric
         & Mean
         & Max
         & Min
         & Standard deviation    \\
        \cmidrule(lr){1-5}
        Hausdorff distance
         & $3.16 \times 10^{-2}$
         & $3.56 \times 10^{-2}$
         & $2.89 \times 10^{-2}$
         & $1.86 \times 10^{-3}$ \\
        optimality error
         & $4.23 \times 10^{-3}$
         & $7.55 \times 10^{-3}$
         & $2.11 \times 10^{-3}$
         & $3.67 \times 10^{-3}$ \\
        $L^2$ error
         & $4.59 \times 10^{-3}$
         & $7.38 \times 10^{-3}$
         & $1.74 \times 10^{-3}$
         & $1.81 \times 10^{-3}$ \\
        \bottomrule
    \end{tabular}
\end{table}

\subsubsection{\texorpdfstring{Parametric problem: results with $f$ given by \eqref{eq:parametric_source}}{Parametric problem: results with exponential source term}}
\label{sec:shapo_f_exp_param}

\begin{figure}[!ht]
    \centering
    \begin{tikzpicture}
        \node[xshift=-0.25\textwidth] (top_left) at (0,0)
        {\includegraphics[width=0.45\textwidth]{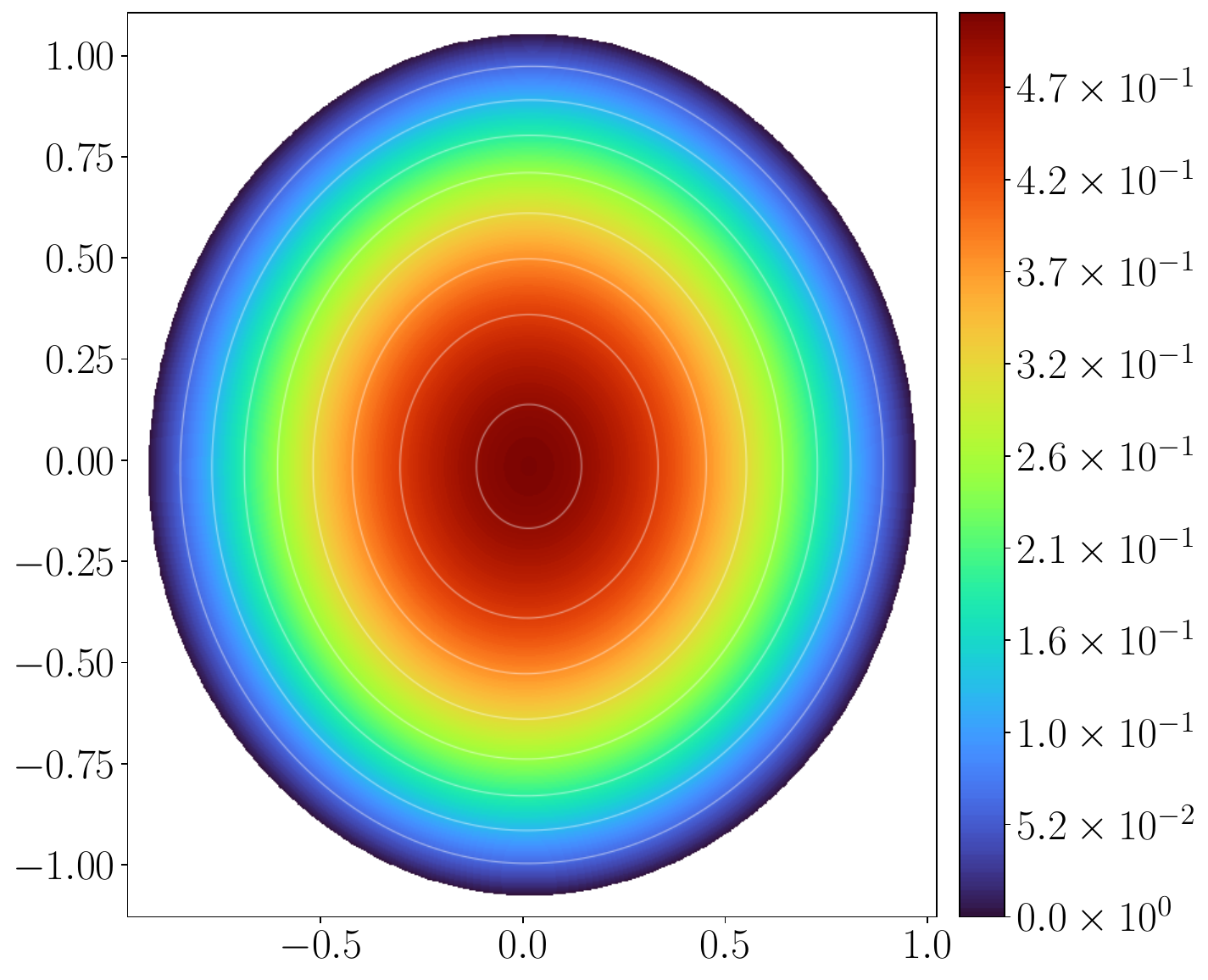}};
        \node[yshift=-0.02\textwidth] at (top_left.south) {(a) solution, $\mu=0.83$};

        \node[xshift=0.25\textwidth, yshift=0.0125\textwidth] (top_right) at (top_left.east)
        {\includegraphics[width=0.45\textwidth]{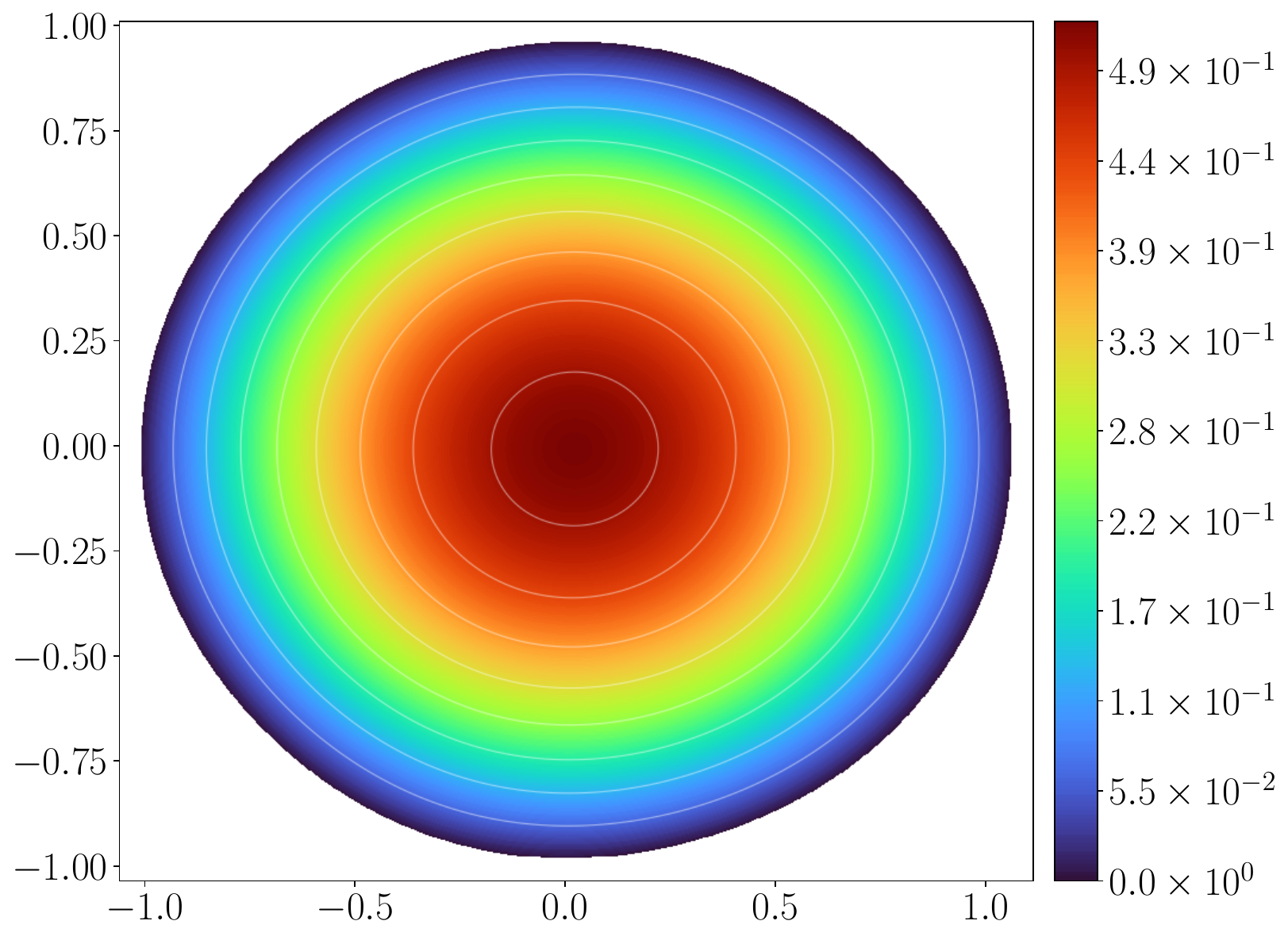}};
        \node[yshift=-0.02\textwidth] at (top_right.south) {(b) solution, $\mu=1.13$};

        \node[yshift=-0.23\textwidth] (bottom_left) at (top_left.south)
        {\includegraphics[width=0.45\textwidth]{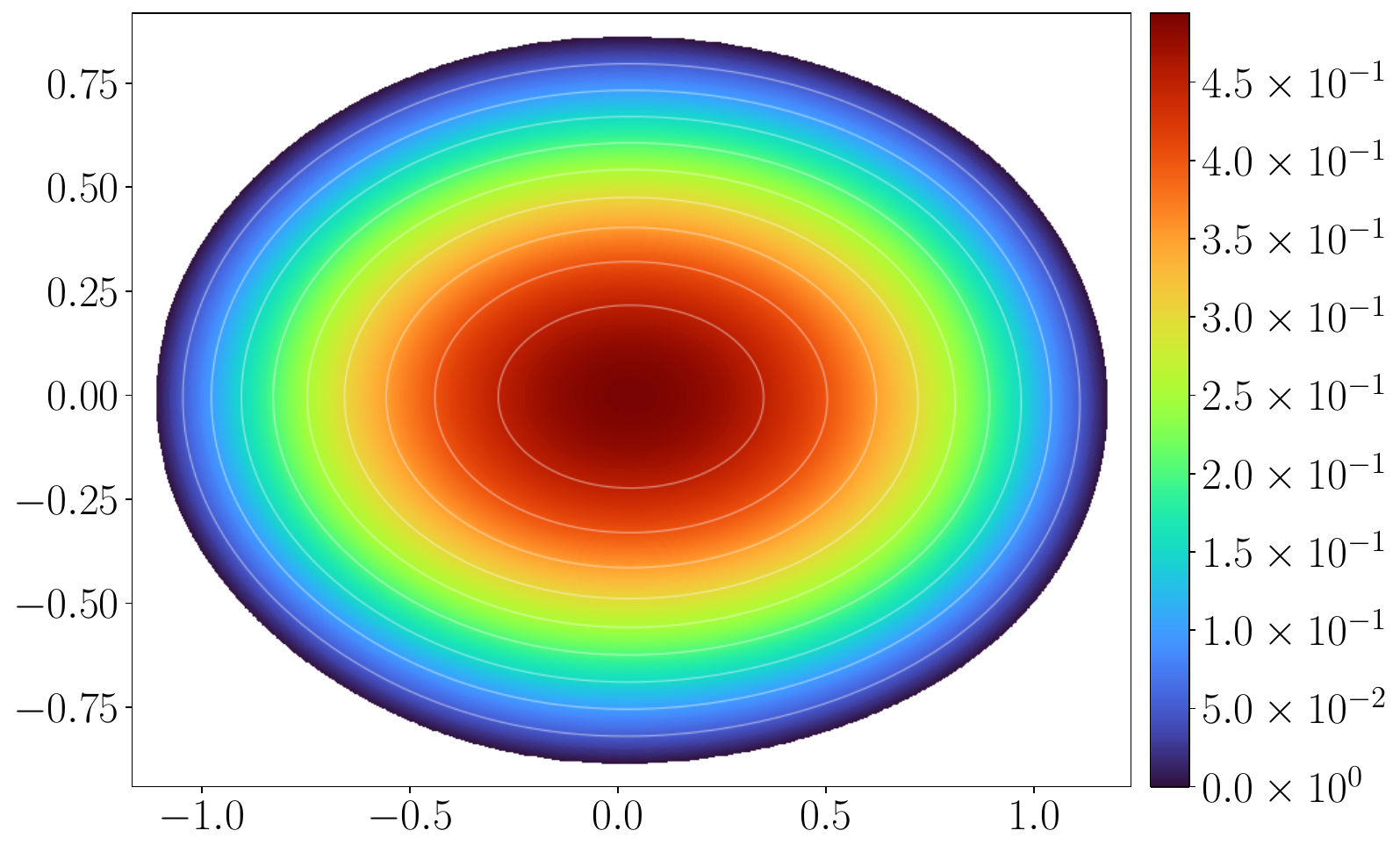}};
        \node[yshift=-0.02\textwidth] at (bottom_left.south) {(c) solution, $\mu=1.61$};

        \node[yshift=-0.23\textwidth] (bottom_right) at (top_right.south)
        {\includegraphics[width=0.45\textwidth]{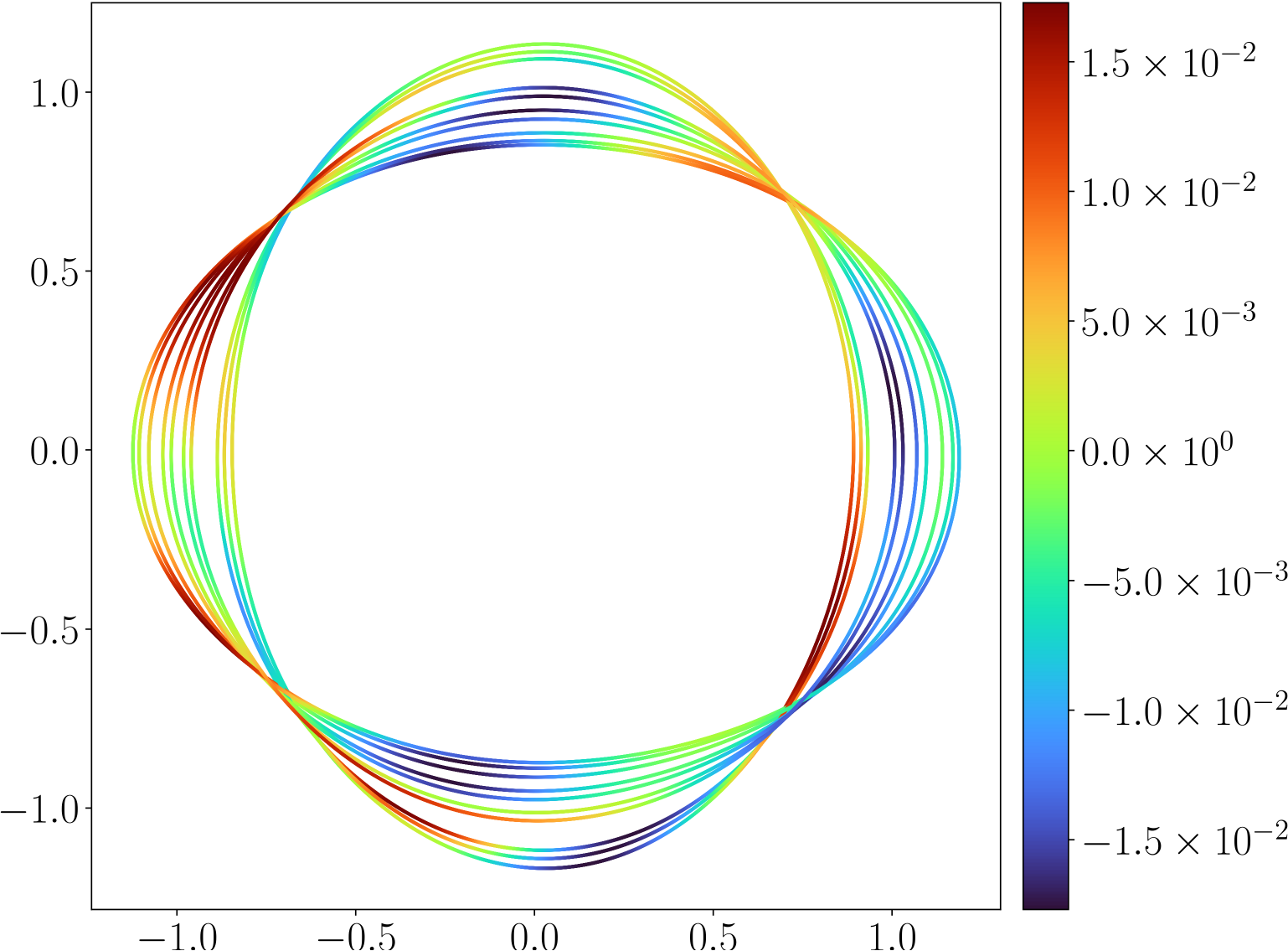}};
        \node[yshift=-0.02\textwidth] at (bottom_right.south) {(d) deviation from the average of the optimality condition};

    \end{tikzpicture}

    \caption{%
        Shape optimization for the Dirichlet energy,
        for Poisson's equation with source term given by~\eqref{eq:parametric_source},
        in the context of \cref{sec:shapo_f_exp_param}.
        For three values of $\mu$, we display the approximate solution
        (top left, top right, and bottom left panels)
        and the deviation from the average of the optimality condition for $10$ random values of $\mu$ (bottom right panel).
    }
    \label{fig:shapo_f_exp_param}
\end{figure}

We now go back to the Poisson problem with
source term given by \eqref{eq:parametric_source}.
We learn the optimal shape for $\mu \in \mathbb{M} = (0.5, 2)$,
and display the results on \cref{fig:shapo_f_exp_param}.
Namely, we observe that the optimality condition is almost satisfied
for all values of $\mu$.
Then, we also report statistics on the optimality condition
and, since no exact solution is available,
the value of the integral in \eqref{eq:fv_diff_forme_for_stats}.
We also observe a good quantitative behavior of our approach.

\begin{table}[!ht]
    \centering
    \caption{%
        Statistics of relevant metrics in the case of a
        parametric source term given by \eqref{eq:parametric_source},
        obtained by computing each metric for $10^3$ values of $\mu$.%
    }
    \label{tab:stats_DeepShape_exp_param}
    \begin{tabular}{ccccc}
        \toprule
        Metric
         & Mean
         & Max
         & Min
         & Standard deviation    \\
        \cmidrule(lr){1-5}
        optimality error
         & $1.37 \times 10^{-2}$
         & $2.23 \times 10^{-2}$
         & $8.31 \times 10^{-3}$
         & $9.19 \times 10^{-3}$ \\
        variational formulation
         & $5.89 \times 10^{-2}$
         & $1.58 \times 10^{-1}$
         & $3.21 \times 10^{-3}$
         & $2.89 \times 10^{-2}$ \\
        \bottomrule
    \end{tabular}
\end{table}

\subsubsection{Parametric problem: results with a non-radial source term}
\label{sec:shapo_f_bizaroid_param}

\begin{figure}[!ht]
    \centering
    \begin{tikzpicture}
        \node[xshift=-0.25\textwidth] (top_left) at (0,0)
        {\includegraphics[width=0.45\textwidth]{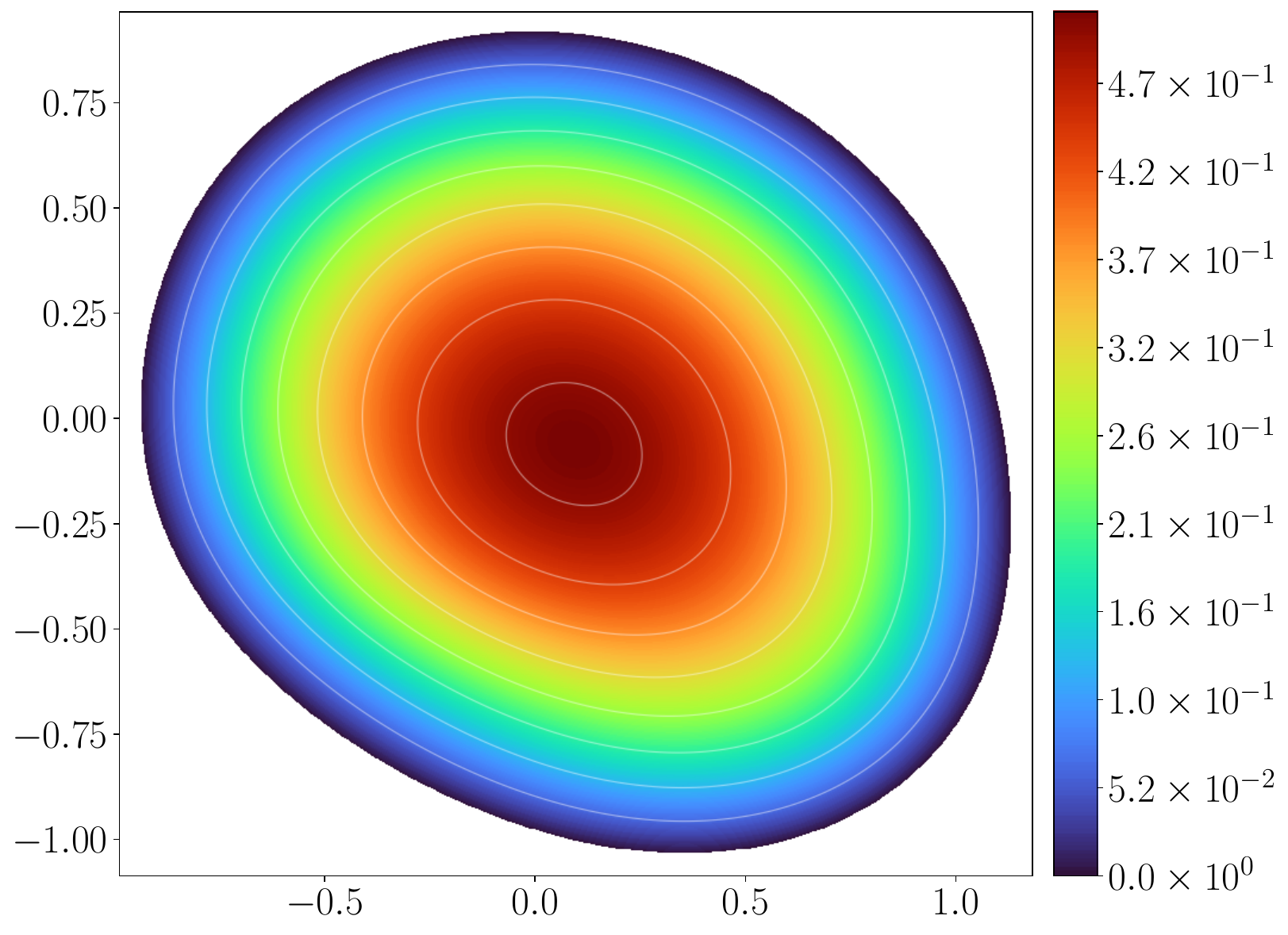}};
        \node[yshift=-0.02\textwidth] at (top_left.south) {(a) solution, $\lambda=0.69$};

        \node[xshift=0.25\textwidth] (top_right) at (top_left.east)
        {\includegraphics[width=0.45\textwidth]{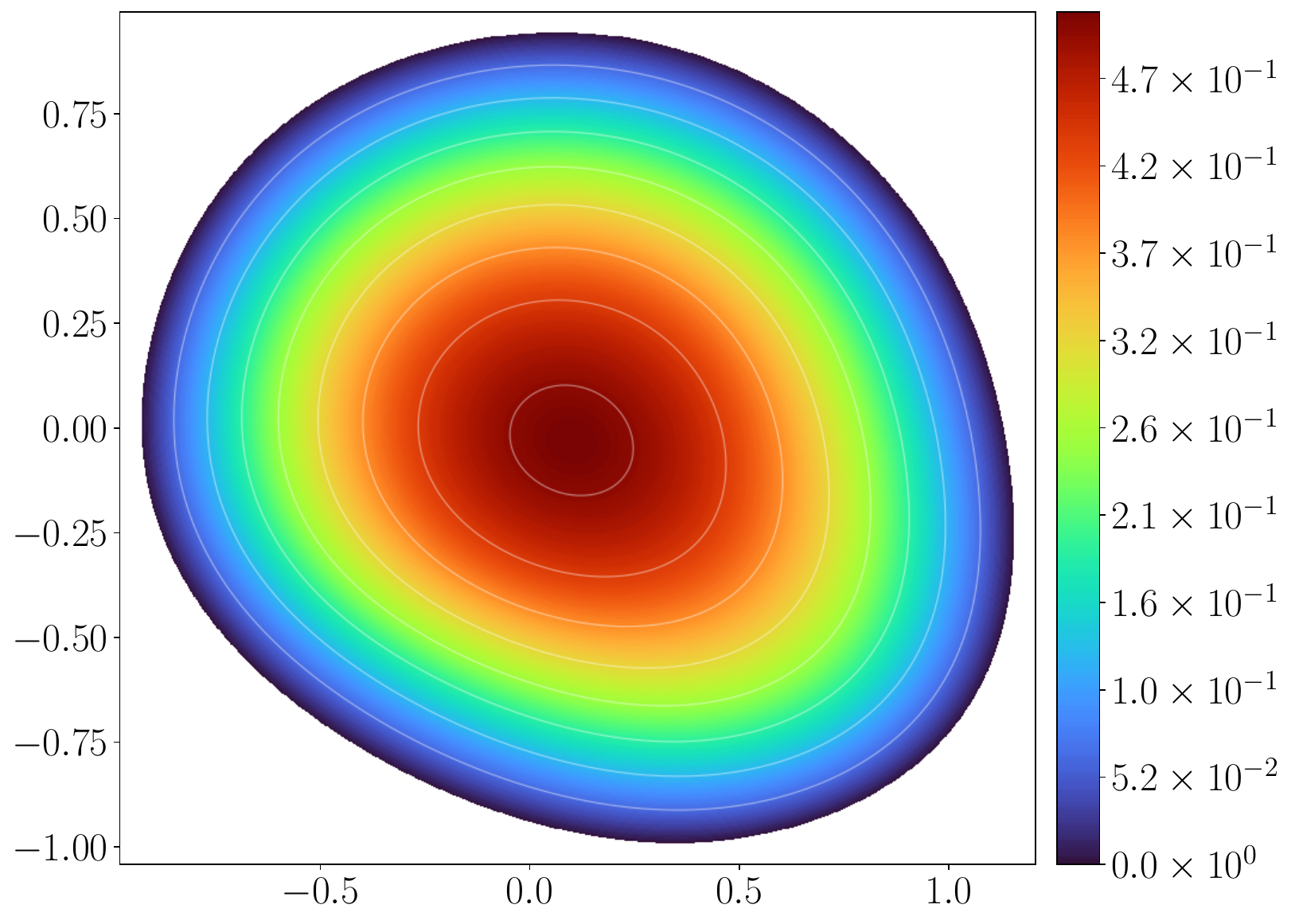}};
        \node[yshift=-0.02\textwidth] at (top_right.south) {(b) solution, $\lambda=1.08$};

        \node[yshift=-0.23\textwidth] (bottom_left) at (top_left.south)
        {\includegraphics[width=0.45\textwidth]{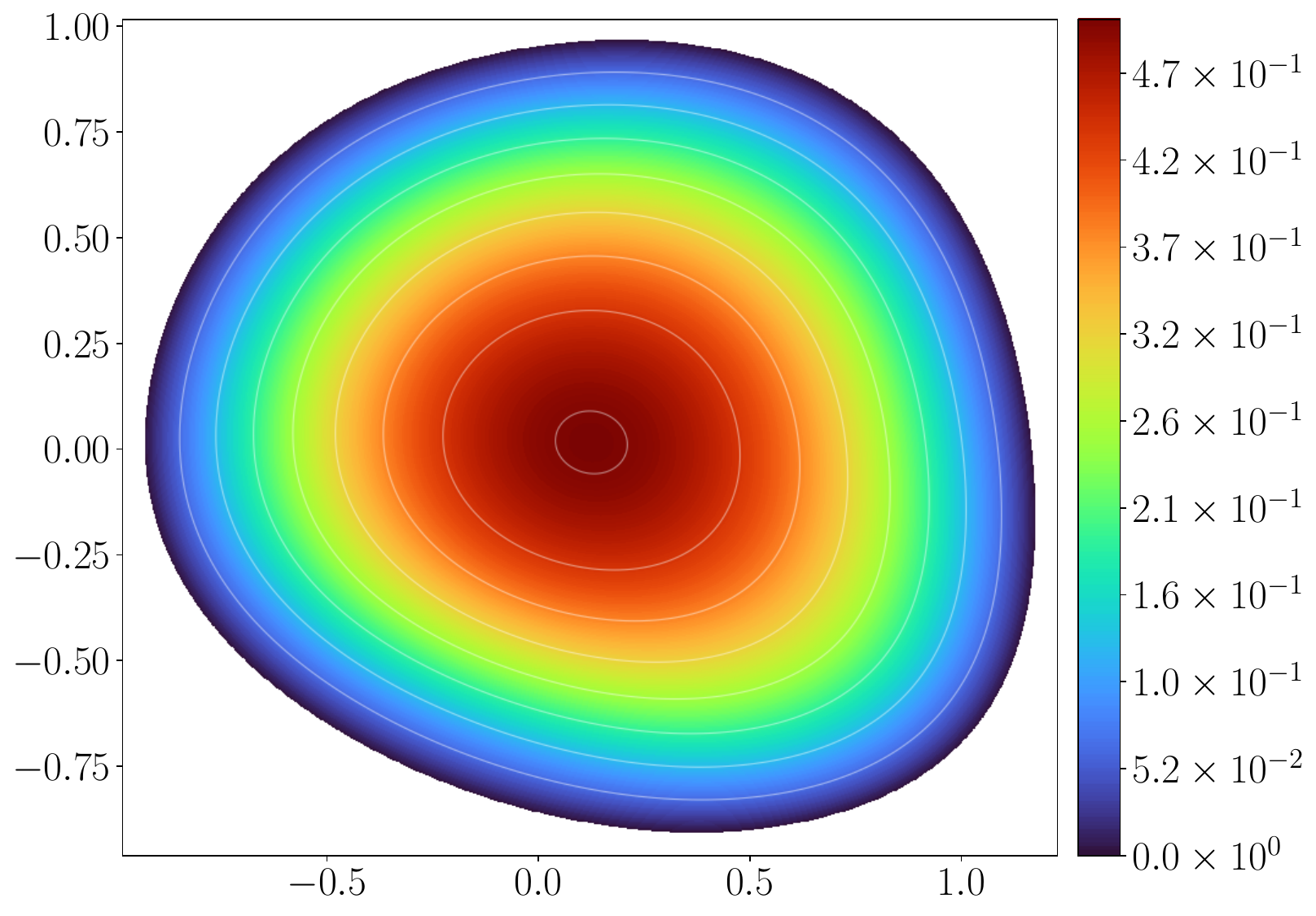}};
        \node[yshift=-0.02\textwidth] at (bottom_left.south) {(c) solution, $\lambda=1.82$};

        \node[yshift=-0.23\textwidth] (bottom_right) at (top_right.south)
        {\includegraphics[width=0.45\textwidth]{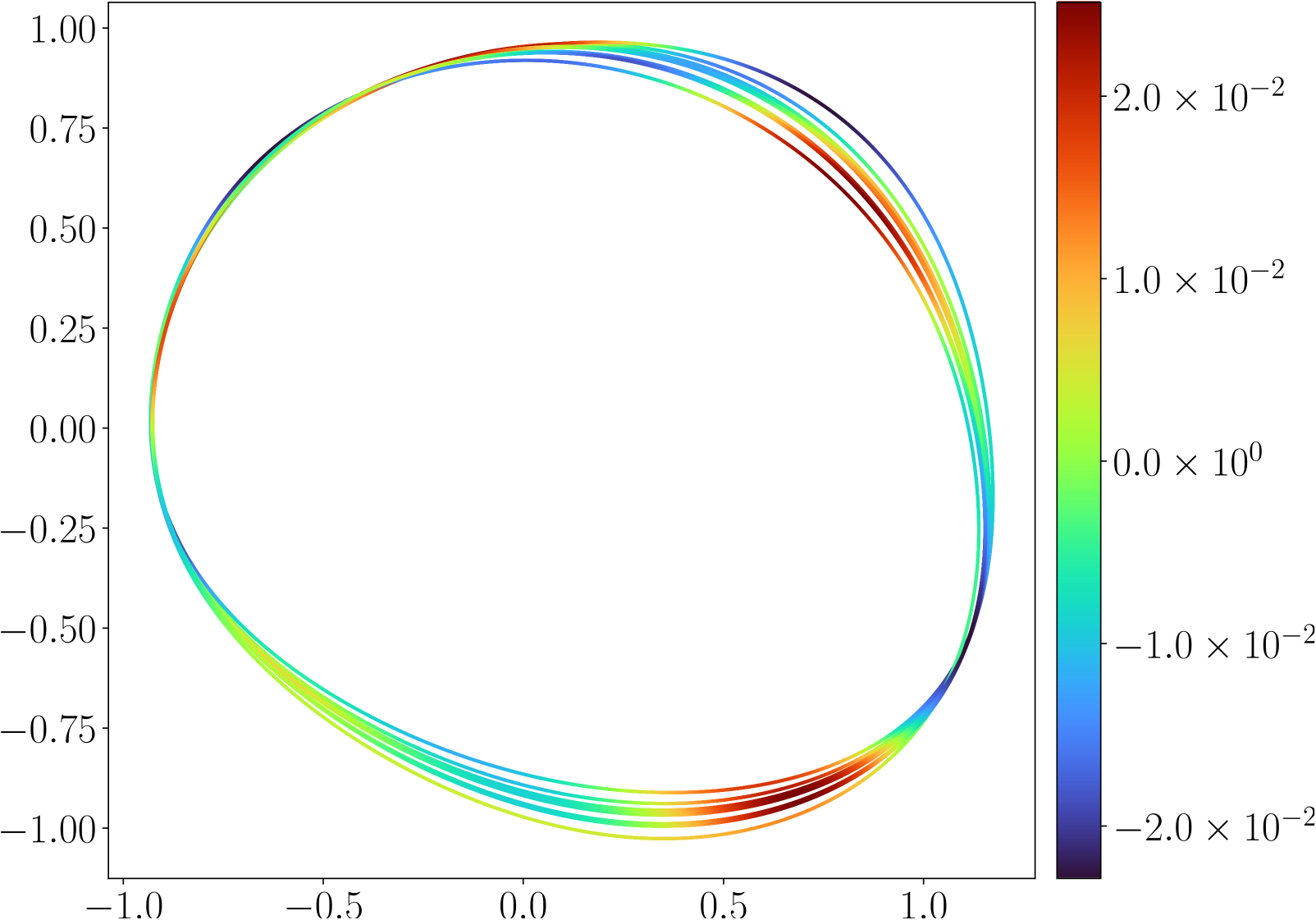}};
        \node[yshift=-0.02\textwidth] at (bottom_right.south) {(d) deviation from the average of the optimality condition};

    \end{tikzpicture}

    \caption{%
        Parametric shape optimization for the Dirichlet energy
        in the context of \cref{sec:shapo_f_bizaroid_param}.
        For three values of $\lambda$, we display the approximate solution
        (top left, top right, and bottom left panels)
        and the deviation from the average of the optimality condition for $10$ random values of $\lambda$ (bottom right panel).
    }
    \label{fig:shapo_f_bizaroid_param}
\end{figure}

For this last test case concerning the Poisson problem with Dirichlet boundary conditions,
we take a source term given by
$f(x, y; \lambda) = \exp (1 - \| \mathcal T_\lambda(x, y) \|^2)$,
with $\mathcal T_\lambda$ the symplectic map defined in \eqref{eq:bizaroid}
and $\lambda \in (0.5, 2)$.
The results are depicted on \cref{fig:shapo_f_bizaroid_param},
with a good observed behavior of our approach.
The same statistics as in \cref{sec:shapo_f_exp_param}
are reported in \cref{tab:stats_DeepShape_bizaroid_param},
which further confirm the relevance of our method.

\begin{table}[!ht]
    \centering
    \caption{%
        Statistics of relevant metrics
        for the parametric source term given in \cref{sec:shapo_f_bizaroid_param},
        obtained by computing each metric for $10^3$ values of $\lambda$.%
    }
    \label{tab:stats_DeepShape_bizaroid_param}
    \begin{tabular}{ccccc}
        \toprule
        Metric
         & Mean
         & Max
         & Min
         & Standard deviation    \\
        \cmidrule(lr){1-5}
        optimality error
         & $2.04 \times 10^{-2}$
         & $3.36 \times 10^{-2}$
         & $1.24 \times 10^{-2}$
         & $1.57 \times 10^{-2}$ \\
        variational formulation
         & $3.36 \times 10^{-2}$
         & $1.12 \times 10^{-1}$
         & $7.27 \times 10^{-5}$
         & $2.15 \times 10^{-2}$ \\
        \bottomrule
    \end{tabular}
\end{table}

\subsection{The exterior Bernoulli free-boundary problem}\label{sec:bernoulli}
To check the relevance of our methodology in a more complex case,
we consider the exterior Bernoulli free-boundary problem,
described for instance in \cite{henrot2000existence}.
Let $\Upomega$ be an open bounded connected set in $\mathbb R^n$,
such that $\partial \Upomega = \Gamma_e \cup \Gamma_i$,
and $K \subset \Upomega$ a compact set
such that $\partial K = \Gamma_i$.
These sets are represented to the right of equation \eqref{eq:bernoulli}.
The Bernoulli problem is an overdetermined PDE,
whose unique solution $u_\Upomega^K \in H^1_0(\Upomega)$ satisfies:
\noindent\begin{minipage}{0.3\textwidth}
    \begin{equation}
        \label{eq:bernoulli}
        \begin{cases}
            -\Updelta u^K_\Upomega = 0 & \text{in }\Upomega, \\
            u^K_\Upomega=1             & \text{on }\Gamma_i, \\
            u^K_\Upomega=0             & \text{on }\Gamma_e, \\
            |\nabla u^K_\Upomega | = c & \text{on }\Gamma_e,
        \end{cases}
    \end{equation}
\end{minipage}%
\begin{minipage}{0.3\textwidth}
    \
\end{minipage}%
\begin{minipage}{0.4\textwidth}

    \tikzset{every picture/.style={line width=0.75pt}} 

    \begin{tikzpicture}[x=0.75pt,y=0.75pt,yscale=-1,xscale=1,scale=0.6]
    
        \draw  [fill={rgb, 255:red, 255; green, 193; blue, 193 }  ,fill opacity=1 ] (307,11) .. controls (408,36) and (422,-33) .. (425,65) .. controls (425.25,73.07) and (425.62,80.37) .. (426.1,86.97) .. controls (426.23,88.8) and (426.37,90.57) .. (426.52,92.3) .. controls (432.41,159.94) and (451.1,149.16) .. (469,176) .. controls (489,206) and (434,213) .. (371,230) .. controls (308,247) and (270,228) .. (229,184) .. controls (188,140) and (172.25,131.63) .. (209,74) .. controls (245.75,16.38) and (206,-14) .. (307,11) -- cycle ;
        \draw  [fill={rgb, 255:red, 221; green, 221; blue, 221 }  ,fill opacity=1 ] (284,87) .. controls (304,77) and (394,67) .. (374,87) .. controls (354,107) and (319,124) .. (374,147) .. controls (429,170) and (304,177) .. (284,147) .. controls (264,117) and (264,97) .. (284,87) -- cycle ;
    
        \draw (246,21.4) node [anchor=north west][inner sep=0.75pt]    {$\textcolor[rgb]{1,0,0}\upOmega$};
        \draw (184,39.4) node [anchor=north west][inner sep=0.75pt]    {$\Gamma_{e}$};
        \draw (326,78.4) node [anchor=north west][inner sep=0.75pt]    {$K$};
        \draw (285,51.4) node [anchor=north west][inner sep=0.75pt]    {$\Gamma_{i}$};

    \end{tikzpicture}
    \end{minipage}

\noindent for a given constant $c>0$. This problem has been extensively studied. For instance, it can be demonstrated that if $K$ is star-shaped (resp. convex), there exists a unique regular solution $\Omega$ that is also star-shaped (resp. convex), see \cite[Chapter~6]{henrot-pierre}.

The solution $u_\Upomega^K$ of \eqref{eq:bernoulli} remains a minimum of the Dirichlet energy
$\mathcal{E}_K(\Upomega)=\frac 12\int_\Upomega |\nabla u_\Upomega^K|^2$.
Therefore, our method based on PINN and SympNets
remains relevant to solve this problem.
We require adapting this strategy to the problem at hand;
this is described in \cref{sec:numerical_strategy_bernoulli}.
Numerical results are then presented in \cref{sec:numerical_results_bernoulli}.

\subsubsection{\texorpdfstring{Numerical strategy for solving the Bernoulli problem with {\rb neural networks}}{Numerical strategy for solving the Bernoulli problem with neural networks}}
\label{sec:numerical_strategy_bernoulli}

To implement Bernoulli's overdetermined problem, we need to know $\varphi_i$ and $\varphi_e$, the respective level-set functions of $\partial K = \Gamma_i$ and $\Gamma_e$.
Here, we assume $\varphi_i$ to be known, since $K$ is a given set.
As an example and for simplicity,~$K$ is taken as an ellipse centered at the origin, with width $2a$ and height $2b$.

Like in \cref{sec:strategy_shapo}, we want to transform an initial shape $\Upomega$ by a symplectic map~$\mathcal T$ into an optimal shape $\Upomega^*$. Following \eqref{eq:bernoulli}, we denote by $\Gamma_i$ the inner boundary of $\Upomega$ and $\Gamma_e$ the outer boundary of $\Upomega$.
The main numerical challenge lies in mapping~$\Gamma_i$ onto $\partial K$ by the symplectic map $\mathcal T$. A solution could be using penalization in the loss function, but it leads to a numerically ill-posed problem. Therefore, we look for a more intrinsic method.

\begin{figure}[!ht]
    \centering

    \tikzset{every picture/.style={line width=0.75pt}} 

    \begin{tikzpicture}[x=0.75pt,y=0.75pt,yscale=-1,xscale=1]
    
    \draw   (58,78) .. controls (58,38.79) and (89.79,7) .. (129,7) .. controls (168.21,7) and (200,38.79) .. (200,78) .. controls (200,117.21) and (168.21,149) .. (129,149) .. controls (89.79,149) and (58,117.21) .. (58,78) -- cycle ;
    \draw    (203.5,52) .. controls (263.4,39.13) and (404.15,65.46) .. (449.65,91.22) ;
    \draw [shift={(451,92)}, rotate = 210.58] [color={rgb, 255:red, 0; green, 0; blue, 0 }  ][line width=0.75]    (10.93,-3.29) .. controls (6.95,-1.4) and (3.31,-0.3) .. (0,0) .. controls (3.31,0.3) and (6.95,1.4) .. (10.93,3.29)   ;
    \draw   (464,94) .. controls (484,84) and (580,100) .. (560,120) .. controls (540,140) and (540,150) .. (560,180) .. controls (580,210) and (481,237) .. (461,207) .. controls (441,177) and (444,104) .. (464,94) -- cycle ;
    \draw  [fill={rgb, 255:red, 0; green, 0; blue, 0 }  ,fill opacity=1 ] (480,138) .. controls (490,110) and (509,121) .. (514,139) .. controls (519,157) and (483,146) .. (515,170) .. controls (547,194) and (505,189) .. (491,177) .. controls (477,165) and (470,166) .. (480,138) -- cycle ;
    \draw   (56,230) .. controls (56,190.79) and (87.79,159) .. (127,159) .. controls (166.21,159) and (198,190.79) .. (198,230) .. controls (198,269.21) and (166.21,301) .. (127,301) .. controls (87.79,301) and (56,269.21) .. (56,230) -- cycle ;
    \draw  [fill={rgb, 255:red, 0; green, 0; blue, 0 }  ,fill opacity=1 ] (108,216) .. controls (79,193) and (125,203) .. (142,217) .. controls (159,231) and (111,224) .. (143,248) .. controls (175,272) and (142,248) .. (119,255) .. controls (96,262) and (137,239) .. (108,216) -- cycle ;
    \draw    (206,268) .. controls (257.48,273.94) and (406.97,234.8) .. (448.77,204.9) ;
    \draw [shift={(450,204)}, rotate = 143.13] [color={rgb, 255:red, 0; green, 0; blue, 0 }  ][line width=0.75]    (10.93,-3.29) .. controls (6.95,-1.4) and (3.31,-0.3) .. (0,0) .. controls (3.31,0.3) and (6.95,1.4) .. (10.93,3.29)   ;
    \draw    (209.22,228.42) .. controls (259.3,236.74) and (399.6,195.55) .. (439,166) ;
    \draw [shift={(207,228)}, rotate = 11.77] [color={rgb, 255:red, 0; green, 0; blue, 0 }  ][line width=0.75]    (10.93,-3.29) .. controls (6.95,-1.4) and (3.31,-0.3) .. (0,0) .. controls (3.31,0.3) and (6.95,1.4) .. (10.93,3.29)   ;
    
    \draw (61,7.4) node [anchor=north west][inner sep=0.75pt]    {$\mathbb S^1$};
    \draw (322,34.4) node [anchor=north west][inner sep=0.75pt]    {$\mathcal T$};
    \draw (570,106.4) node [anchor=north west][inner sep=0.75pt]    {$\mathcal T\mathbb S^1$};
    \draw (512,113.4) node [anchor=north west][inner sep=0.75pt]    {$K$};
    \draw (308,187.4) node [anchor=north west][inner sep=0.75pt]    {$\mathcal T^{-1}$};
    \draw (346,250.4) node [anchor=north west][inner sep=0.75pt]    {$\mathcal T$};
    \draw (65,154.4) node [anchor=north west][inner sep=0.75pt]    {$\mathbb S^1$};
    \draw (121,183.4) node [anchor=north west][inner sep=0.75pt]    {$\mathcal T^{-1} K$};

    \end{tikzpicture}
    \caption{Sampling method for the numerical resolution of the Bernoulli overdetermined PDE.}
    \label{fig:bernoulli_compute}
\end{figure}
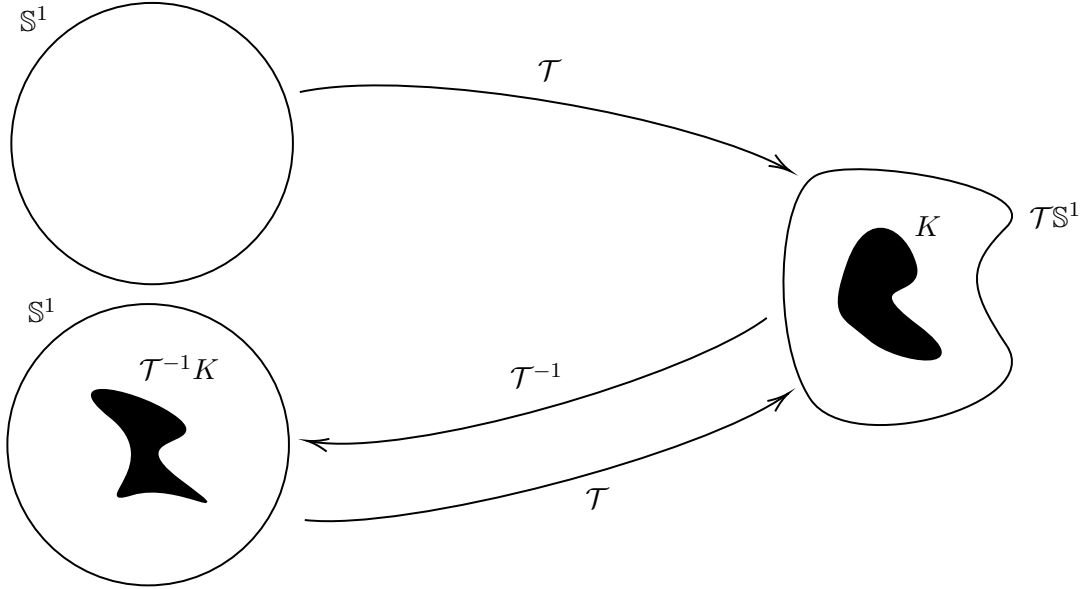

On \cref{fig:bernoulli_compute}, the numerical implementation is explained. We start from an initial shape $\Upomega$, which is the unit ball $\mathbb B^2$. After applying the symplectic map $\mathcal T$ on $\mathbb B^2$ to obtain a computational domain $\mathcal T\mathbb B^2$, a mask is applied in order to delete the collocation points mapped within $K$.

To implement the boundary conditions, we follow the strategy presented in \cref{sec:strategy_shapo}. Equipped with a PINN $u_{\theta}:\mathbb B^2 \to \mathbb R$ and a SympNet $T_{\omega}:\mathbb R^2 \to \mathbb R^2$, the approximate solution $v_{\theta, \omega}$ is represented as follows:
\begin{equation*}
    \begin{array}{r c c l}
        v_{\theta, \omega}: & \mathbb B^2 & \to     & \mathbb R                                                                           \\
                            & (x_1,x_2)   & \mapsto & \alpha_\omega(x_1,x_2)(u_\theta \circ T_{\omega})(x_1,x_2) + \beta_\omega(x_1,x_2), \\
    \end{array}
\end{equation*}
with $\alpha_\omega = 0$ on $\Gamma_i \cup \Gamma_e$,
$\beta_\omega = 0$ on $\Gamma_e$ and $\beta_\omega = 1$ on $\Gamma_i$.
To obtain functions satisfying these criteria, we take $\alpha_\omega,\beta_\omega : \mathbb B^2 \setminus T_\omega^{-1}K \to T_\omega\mathbb B^2 \setminus K$ such that
\begin{equation*}
    \alpha_\omega = \varphi_\omega \varphi_{\mathbb B^2}
    \text{\qquad and \qquad}
    \beta_\omega = \frac{\varphi_{\mathbb B^2}}{\varphi_\omega - \varphi_{\mathbb B^2}},
\end{equation*}
where, for $(x_1, x_2)\in \mathbb B^2\setminus T^{-1}_\omega K$,
\begin{equation*}
    \begin{aligned}
        \varphi_\omega(x_1,x_2)        & =
        \left(\frac{x_1^{T_\omega}}{a}\right)^2 + \left(\frac{x_2^{T_\omega}}{b}\right)^2 - 1
        \text{, \quad with }
        (x_1^{T_\omega}, x_2^{T_\omega}) = T_\omega(x_1, x_2); \\
        \varphi_{\mathbb B^2}(x_1,x_2) & =
        x^2+y^2 - \rho_\text{max}^2
        \text{,\quad with }\rho_\text{max}\text{ the radius of the ball }\mathbb B^2.
    \end{aligned}
\end{equation*}

\subsubsection{Numerical results}
\label{sec:numerical_results_bernoulli}

\begin{figure}[!ht]
    \centering
    \begin{tikzpicture}
        \node[xshift=-0.25\textwidth] (top_left) at (0,0)
        {\includegraphics[width=0.45\textwidth]{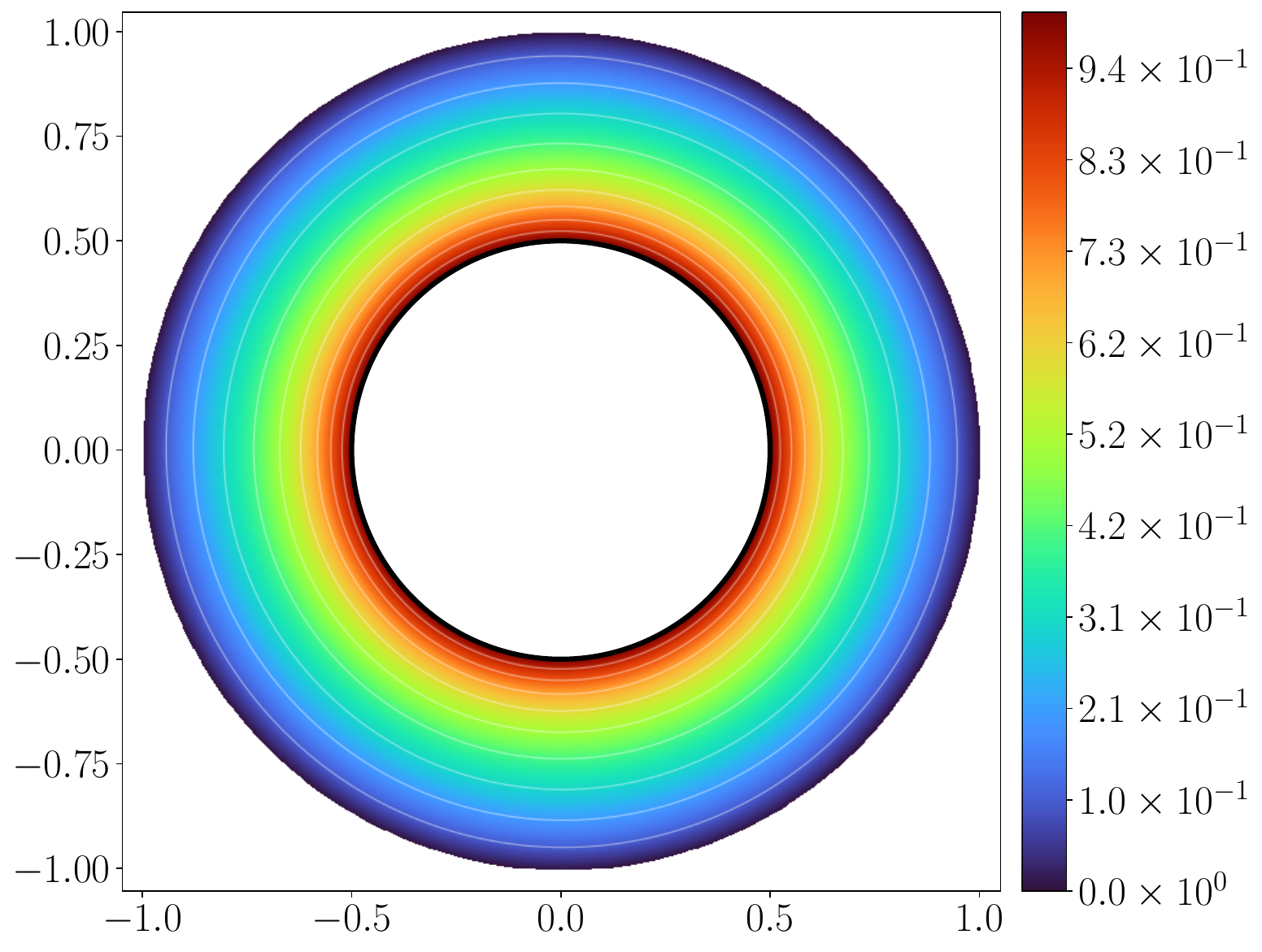}};
        \node[yshift=-0.02\textwidth] at (top_left.south) {\small (a) approximate solution ($a=0.5$)};

        \node[xshift=0.25\textwidth] (top_right) at (top_left.east)
        {\includegraphics[width=0.4625\textwidth]{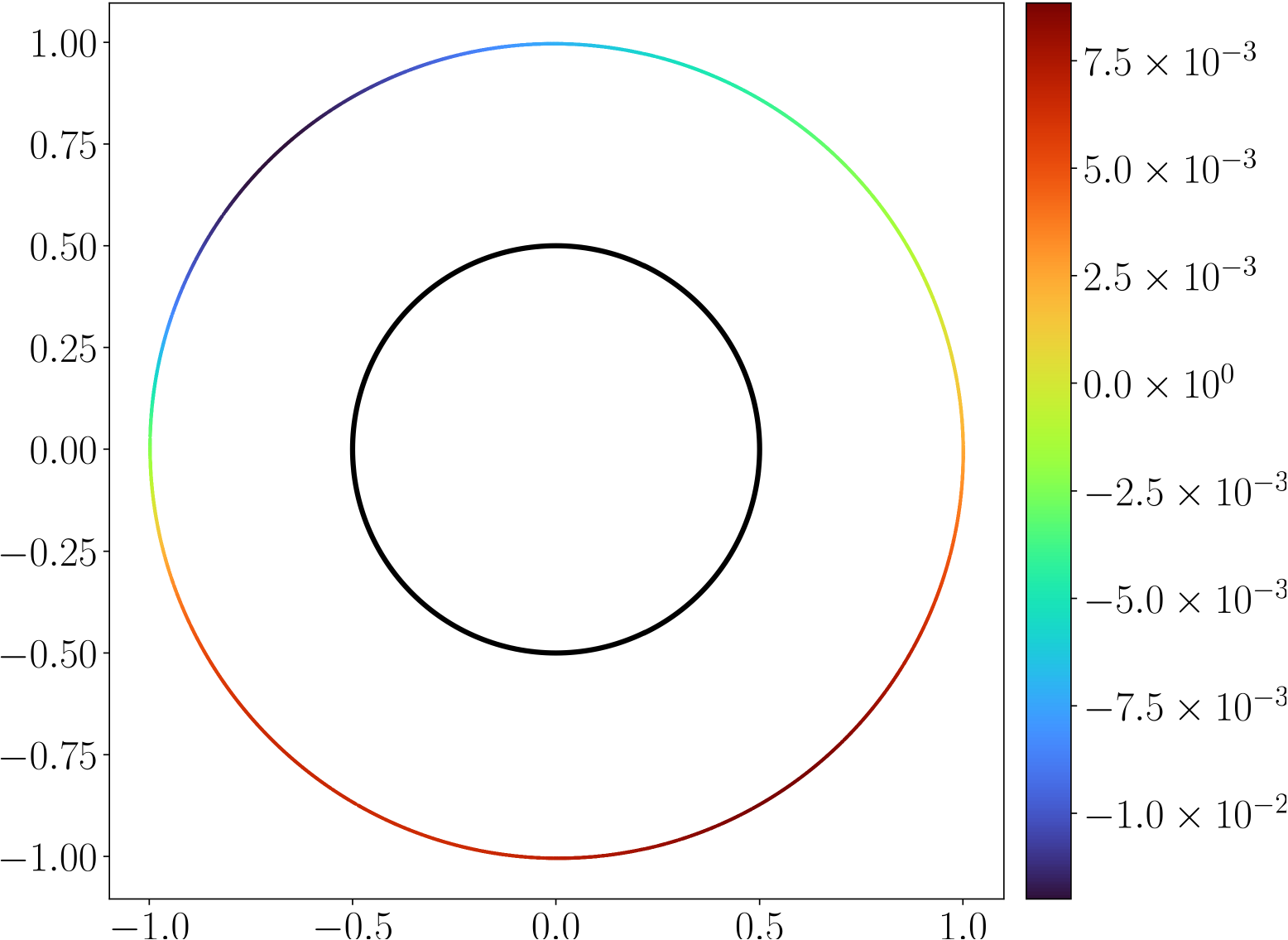}};
        \node[yshift=-0.02\textwidth] at (top_right.south) {\small (b) deviation from the average of the optimality condition ($a=0.5$)};

        \node[yshift=-0.21\textwidth] (bottom_left) at (top_left.south)
        {\includegraphics[width=0.45\textwidth]{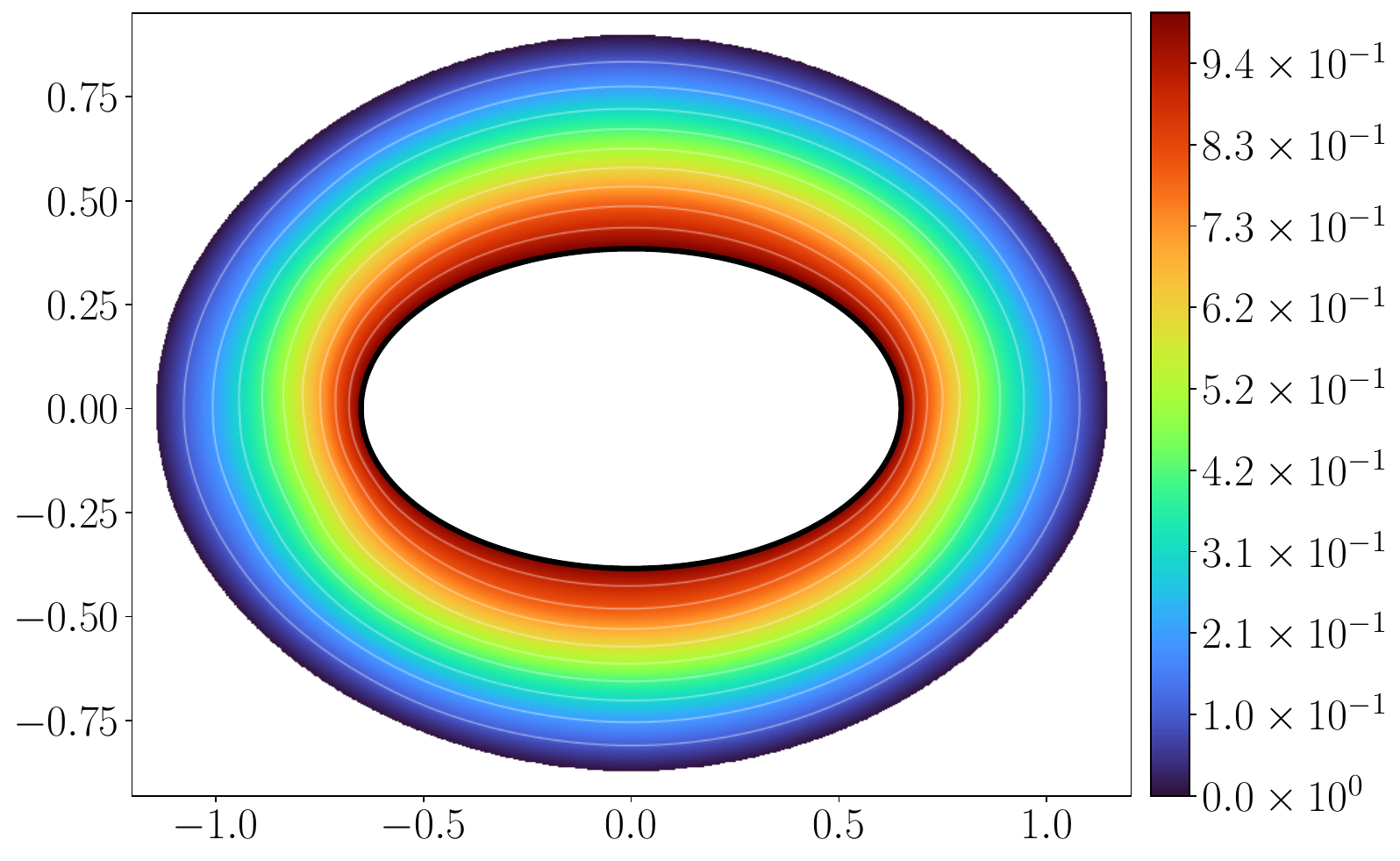}};
        \node[yshift=-0.02\textwidth] at (bottom_left.south) {\small (c) approximate solution ($a=0.65$)};

        \node[yshift=-0.21\textwidth] (bottom_right) at (top_right.south)
        {\includegraphics[width=0.4625\textwidth]{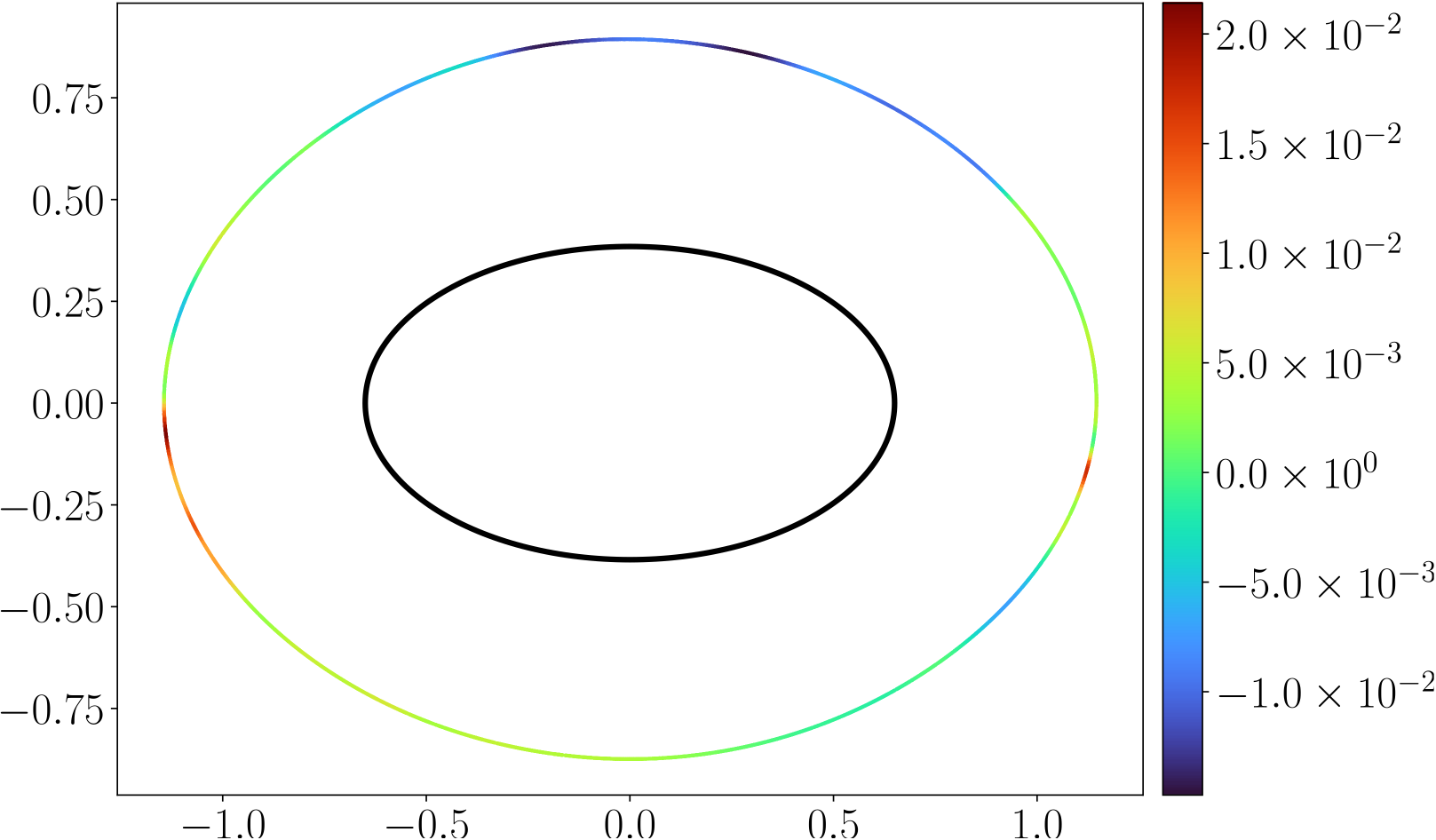}};
        \node[yshift=-0.02\textwidth] at (bottom_right.south) {\small (d) deviation from the average of the optimality condition ($a=0.65$)};
    \end{tikzpicture}

    \caption{%
        Approximate solutions (left panels) and optimality conditions (right panels)
        for the Bernoulli problem with an ellipsoid-shaped obstacle
        (top panels: $a=0.5$, i.e., the case of a disk;
        bottom panels: $a=0.65$).
        On all panels, the black line corresponds to $\partial K = \Gamma_i$.
    }
    \label{fig:bernoulli}
\end{figure}
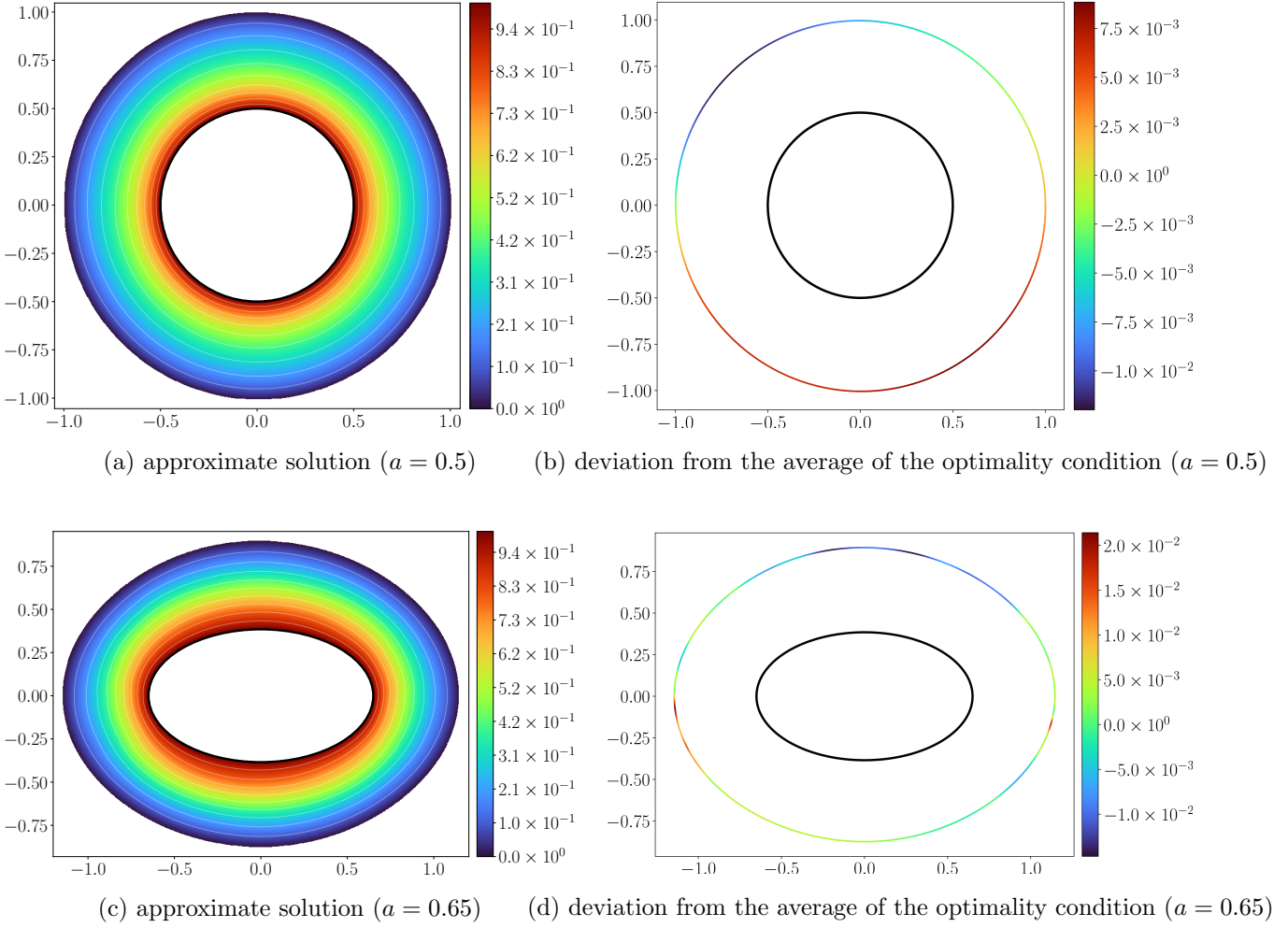

We start with a case where the exact solution is known:
with a centered disk-shaped obstacle (i.e., with $a=0.5$),
the optimal shape will also be a disk.
Indeed, assume that $\Upomega$ is a centered disk with radius $R>0$.
Then, one easily check that a particular solution of the
Bernoulli exterior problem \eqref{eq:bernoulli} is
$u(x,y)=(\log R-\log \sqrt{x^2+y^2})/\log (2R)$,
where $R$ is the unique positive number such that $R\, |\!\log (2R)|=1/c$.
We conclude by uniqueness of the solution
to \eqref{eq:bernoulli} (see e.g.~\cite[Theorem~1.1]{henrot2000existence}).
The results are displayed in the top panels of \cref{fig:bernoulli},
where we observe that the optimality condition is satisfied
(the optimality error is $\num{6.83e-3}$).
Moreover, the Hausdorff distance between the learned shape and $\mathbb B^2$
is $\num{6.23e-3}$, which further confirms the relevance of our approximation.

Then, we tackle the case of an elongated ellipse ($a=0.65$).
This case is harder to solve, as the optimal shape is not $\mathbb B^2$,
and the SympNet has to learn a complex transformation,
involving constraints on both $\Gamma_i$ and $\Gamma_e$.
In this case, to help train the model,
we have added a term in the loss function
to penalize the optimality condition on $\Gamma_e$.
This term is activated after \num{22000} epochs,
and is accompanied by a reduction in the learning rate,
from \num{1e-2} to \num{1e-4}.
The results are displayed in the bottom panels of \cref{fig:bernoulli}.
As we do not know the optimal shape in this case, we cannot report the Hausdorff distance;
however, the optimality error is \num{7.04e-3},
which is a good indication that the final shape is close to being optimal.

In both cases, we also report in \cref{tab:stats_bernoulli} the values of
the variational formulation~\eqref{eq:fv_diff_forme},
like in the previous sections,
to further confirm the relevance of our approach.
We observe that the values of the variational formulation
are slightly larger than for Poisson's equation,
but that is attributable to the difficulty of the problem
as well as the larger values taken by $u$ in this case.

\begin{table}[!ht]
    \centering
    \caption{%
        Statistics of the variational formulation
        for the overdetermined Bernoulli problem,
        described in \cref{sec:bernoulli},
        for different values of $a$.%
    }
    \label{tab:stats_bernoulli}
    \begin{tabular}{ccccc}
        \toprule
        Value of $a$   &
        Mean           &
        Max            &
        Min            &
        Standard deviation \\
        \cmidrule(lr){1-5}
        $a=0.5$        &
        \num{1.04e-01} &
        \num{2.45e-01} &
        \num{1.50e-04} &
        \num{5.18e-02}     \\
        $a=0.65$       &
        \num{8.95e-02} &
        \num{2.29e-01} &
        \num{5.61e-04} &
        \num{4.60e-02}     \\
        \bottomrule
    \end{tabular}
\end{table}

\section{Minimizing the Dirichlet Energy with prescribed volume for Robin boundary conditions}
\label{sec:robin}

To show how this methodology is generalizable to other boundary conditions, we introduce again $D$ a compact set of $\mathbb{R}^n$, as well as $\Upomega$ an open subset of $D$, $f\in H^{-1}(D)$, and a Robin coefficient $\kappa \in \mathbb R^+$.
We denote by $u^{f,\kappa}_\Upomega\in H^1(\Upomega)$ the unique solution of Poisson's equation in $\Upomega$ with Robin boundary conditions
\begin{equation}
    \label{eq:poisson_robin}
    \begin{cases}
        -\Updelta u^{f,\kappa}_\Upomega = f                               & \text{in }\Upomega,         \\
        \partial_n u^{f,\kappa}_\Upomega + \kappa u^{f,\kappa}_\Upomega=0 & \text{on }\partial\Upomega.
    \end{cases}
\end{equation}
Moreover, $u^{f,\kappa}_\Upomega$ is the minimizer to an energy $\mathcal J$ defined for $v \in H^1(\Upomega)$ by
\begin{equation}
    \label{fv_robin}
    \mathcal{J}(v) = \frac{1}{2} \int_\Upomega |\nabla v|^2 + \frac{\kappa}{2}\int_{\partial \Upomega}v^2 - \int_\Upomega fv.
\end{equation}
This leads to the definition of a new Dirichlet energy, defined for $\Upomega \subset \mathbb R^n$ by
\begin{equation*}
    \mathcal{E}(\Upomega) = \frac{1}{2} \int_\Upomega |\nabla u^{f,\kappa}_\Upomega|^2 + \frac{\kappa}{2}\int_{\partial \Upomega}{u^{f,\kappa}_\Upomega}^2 - \int_\Upomega fu^{f,\kappa}_\Upomega
\end{equation*}
For a prescribed volume $V_0 > 0$, the new shape optimization problem for the Dirichlet energy of the Robin-Poisson \eqref{eq:poisson_robin} reads
\begin{equation}
    \label{eq:optim_ener_robin}
    \inf_{} \lbrace\mathcal E(\Upomega),
    \, \Upomega \text{ bounded open set of }D, \text{ such that }
    \lvert \Upomega \rvert=V_0\rbrace.
\end{equation}

\begin{remark}
    Shape optimization problems involving elliptic equations with Robin boundary conditions are generally more difficult to study than those with Dirichlet conditions. A series of recent works introduced a suitable relaxation of this type of problem in the space of special functions with bounded variation, as developed by De Giorgi and Ambrosio. This approach yields the existence of a generalized minimizer \cite{bucur2015faber,bucur2016robin}. In certain problems, the authors conducted further analysis to prove the existence of a ``regular'' optimal domain (such as an open set).
    Let us highlight that the case where $f(\cdot )=1$ has been fully studied: Saint-Venant's inequality indicates that the optimal domain is a ball \cite{bucur2015saint}. Generally, determining the regularity of the minimizers for this problem is not straightforward. In the scope of this article, we do not dwell on these questions, especially since we are seeking a local minimizer for this shape optimization problem.
\end{remark}

\begin{theorem}
    \label{thm:optimality_condition_robin}
    Assume that there exists a $C^2$ optimal shape $\Upomega^*$.
    If $f\in H^1_\text{loc}(\mathbb R^n)$, then the solution of \eqref{eq:poisson_robin} $u^*$ belongs to $H^3(\Upomega^*)$. Moreover, the first order optimality condition on $\partial \Upomega^*$ reads
    \begin{equation}
        \label{eq:optimality_condition_robin}
        \exists c\in \mathbb{R} \quad \text{such that} \quad \frac12 |\nabla u^*|^2 + \frac{u^*}{2}(\kappa H^* - 2\kappa^2) - f = c,
    \end{equation}
    where $H^*$ is the mean curvature of $\partial \Upomega^*$.
\end{theorem}
For the sake of completeness, the proof of this result, although standard, is provided in \cref{sec:proof_optimality_condition_robin}.

\begin{remark}
    The mean curvature $H$ of $\mathcal{T}\mathbb{B}^2$ is given by
    $H=\mathrm{div}n_{\mathcal{T}}$,
    with $n_{\mathcal{T}}$ the outwards unit normal vector of $\mathbb B^2$.
    To compute this curvature, define
    $\tau(x=(x_1, x_2))=(x_2, -x_1)$ the unit tangent vector of $\mathbb B^2$
    and $\tau_\mathcal{T}$ the unit tangent vector of $\mathcal{T}\mathbb B^2$.
    Then, with $J_\mathcal{T}$ the Jacobian matrix of $\mathcal{T}$,
    we know that for all $x\in \partial\mathbb{B}^2$, the identity
    \begin{equation*}
        \tau_\mathcal{T}(\mathcal{T} x)
        =
        \frac{J_\mathcal{T}(x)\tau(x)}{\lVert J_\mathcal{T}(x)\tau(x) \rVert}
    \end{equation*}
    is satisfied.
    This leads to very easy computation of the mean curvature
    $H$ of $\mathcal T \mathbb B^2$.
\end{remark}

Equipped with these results,
we now explain how to solve the shape optimization problem
in \cref{sec:solving_robin}.
Then, we validate our methodology on several
numerical experiments in \cref{sec:numerics_robin}.

\subsection{Solving the shape optimization problem with Robin boundary conditions}
\label{sec:solving_robin}

To solve the shape optimization problem \eqref{eq:optim_ener_robin} associated to the Robin-Poisson problem \eqref{eq:poisson_robin}, the same change of variables as in \cref{lem:formulation_with_metric_tensor} is required.
We introduce $w = u\circ \mathcal T$ with $\mathcal T$ a symplectic map and $u$ the solution of the Poisson-Robin problem in $\mathcal T \mathbb B^2$.
The following result states the PDE satisfied by $w$ in $\mathbb B^2$.

\begin{lemma}
    \label{lem:formulation_with_metric_tensor_robin}
    Let $\mathcal{T}$ be a symplectic map on $\mathbb R^2$ and $u_\mathcal{T}\in H^1(\mathcal T\mathbb B^2)$, such that $w=u_\mathcal{T}\circ\mathcal{T}\in H^1(\mathbb B^2)$. If $u_\mathcal{T}\in H^1_0(\mathcal T \mathbb B^2)$ is the solution of the Robin-Poisson problem~\eqref{eq:poisson_robin}, then $w\in H^1(\mathbb B^2)$ is the weak solution of
    \begin{equation} \label{eq:diff_form_robin}
        \inf_{v\in H^1(\mathbb{B}^2)}
        \bigg\lbrace\frac12 \int_{\mathbb B^2}
        \left(A \nabla v\cdot \nabla v\right)
        + \frac\kappa2 \int_{\partial\mathbb B^2} |J_{\mathcal T}^{-\intercal}n|v^2
        - \int_{\mathbb B^2}\tilde{f}v
        \bigg\rbrace,
    \end{equation}
    where $A:\mathbb B^2 \to \mathcal{M}_2(\mathbb{R})$ is defined by
    $A = J_\mathcal{T}^{-1} J_\mathcal{T}^{-\intercal}$ (with $J_\mathcal{T}$ the Jacobian matrix of $\mathcal T$, {\rb and $J_\mathcal{T}^{-\intercal} = \left(J_\mathcal{T}^{-1}\right)^\intercal = \left(J_\mathcal{T}^\intercal\right)^{-1}$}), $\tilde{f} = f \circ \mathcal{T}$, and $n$ is the outwards unit normal vector to $\mathbb S^1$.
\end{lemma}
\begin{proof}
    The proof of \cref{lem:formulation_with_metric_tensor_robin} follows the same reasoning as the proof of \cref{lem:formulation_with_metric_tensor}.
    Therefore, we omit it here.
    Nevertheless, we point out that the variational formula of \eqref{eq:poisson_robin} has a boundary term that involves the tangential Jacobian determinant $|J_{\mathcal T}^{-\intercal}n|$ during the change of variables.
\end{proof}

\begin{remark}
    The boundary conditions of \eqref{eq:poisson_robin} can be derived from the variational formulation associated with the energetic formulation \eqref{eq:diff_form_robin}. Notably, they do not affect the space of test functions. Consequently, for numerical methods, minimizing \eqref{fv_robin} is sufficient to satisfy the Robin boundary conditions. Thus, the methodology for the numerical resolution of the Dirichlet energy minimization presented in the previous section can be directly adapted for \eqref{eq:poisson_robin} without the need for any additional function to enforce the boundary condition.
\end{remark}

The loss function to be minimized with respect to the PINN and the SympNet then becomes
\begin{equation}\label{loss_robin}
    \begin{aligned}
        \mathcal{J}_{\mathrm{Robin}}\big(\theta, \omega; \lbrace x^i_\Upomega, x^i_\Gamma, \mu^i \rbrace_{i=1}^{N}\big)
         & = \frac{V_0}{N} \sum_{i=1}^{N} \left\lbrace \frac12 \mathrm{A}_{\omega} \nabla v_{\theta,\omega} \cdot \nabla v_{\theta,\omega}\right\rbrace (x_\Upomega^i;\mu^i) \\
         & + \kappa \frac{2\sqrt{\pi V_0}}{N} \sum_{i=1}^N \frac12 \left\lbrace\mathscr{J}^\Gamma_{\omega} v_{\theta,\omega}^2\right\rbrace(x_\Gamma^i;\mu^i)                \\
         & - \frac{V_0}{N} \sum_{i=1}^N \left\lbrace\tilde{f}_\omega v_{\theta,\omega}\right\rbrace(x_\Upomega^i;\mu^i),
    \end{aligned}
\end{equation}
with
\begin{itemize}
    \item $\theta$ the trainable weights of the PINN, $\omega$ the trainable weights of the SympNet;
    \item $T_{\omega}:\mathbb{R}^{2d} \to \mathbb{R}^{2d}$ the SympNet;
    \item $u_{\theta}:{T}_{\omega}\mathbb B^2\to \mathbb{R}$ the PINN;
    \item $v_{\theta,\omega}=u_\theta \circ T_\omega$ the solution of the Poisson problem set in $T_\omega\mathbb B^2$;
    \item $\mathrm{A}_\omega$ the diffusion matrix defined by $\mathrm{A}_\omega = \mathrm{J}_{T_\omega}^{-1}\mathrm{J}_{T_\omega}^{-\intercal}$;
    \item $\mathrm{J}_\omega^\Gamma = |\mathrm{J}_{\mathrm{T}_\omega}^{-\intercal}n|$ the tangential Jacobian determinant defined on $\mathbb S^1$, with $n$ the normal outward vector of $\mathbb S^1$, and $|.|$ the euclidean norm in $\mathbb R^2$;
    \item $\tilde{f}_\omega:x \in \mathbb B^2 \mapsto (f\circ T_\omega)(x) \in \mathbb R$;
    \item $N$ the number of collocation points;
    \item $\lbrace x^i_\Upomega, x^i_\Gamma, \mu^i \rbrace_{i=1}^{N}\in (\mathbb B^2\times \mathbb S^1 \times \mathbb M)^N$ a set of random collocation points and values of the parameters.
\end{itemize}

\subsection{Numerical experiments}
\label{sec:numerics_robin}

Now, we proceed with numerical experiments.
We first treat non-parametric shape optimization problems,
where the source term $f$ is either constant
(in \cref{sec:non_parametric_robin_one})
or given by a radial function
(in \cref{sec:non_parametric_robin_exp}).
Then, in \cref{sec:parametric_robin_one},
we consider a parametric shape optimization problem,
where the Robin coefficient is parametric.

\begin{remark}\label{rem:optim_cond}
    For the remainder of the article, and with a slight abuse of notation, the wording ``optimality error'' will refer to the standard deviation of \eqref{eq:optimality_condition_robin}. More precisely, we will verify the local optimality of the obtained shape by examining if
    \begin{equation*}
        \text{optimality error}:=\sqrt{\frac{1}{|\partial\Upomega|}\int_{\partial\Upomega}
            \left(\rho(x) - \overline{\rho}\right)^{\!2} \,
            \mathrm{d}\upsigma} = 0,
    \end{equation*}
    where $\rho = \frac12 |\nabla u^*|^2 + \frac{u^*}{2}(\kappa H^* - 2\kappa^2) - f$ on $\partial\Upomega$, and  $\overline{\rho}$ denotes the average value
    of $\rho$ on $\partial\Upomega$.
\end{remark}

\subsubsection{
    \texorpdfstring
    {Non-parametric problem: results with $f=1$ and $\kappa=1$}
    {Non-parametric problem: results with constant f and fixed kappa}
}
\label{sec:non_parametric_robin_one}

\begin{figure}[!ht]
    \centering
    \begin{tikzpicture}
        \node[xshift=-0.25\textwidth] (top_left) at (0,0)
        {\includegraphics[width=0.353\textwidth]{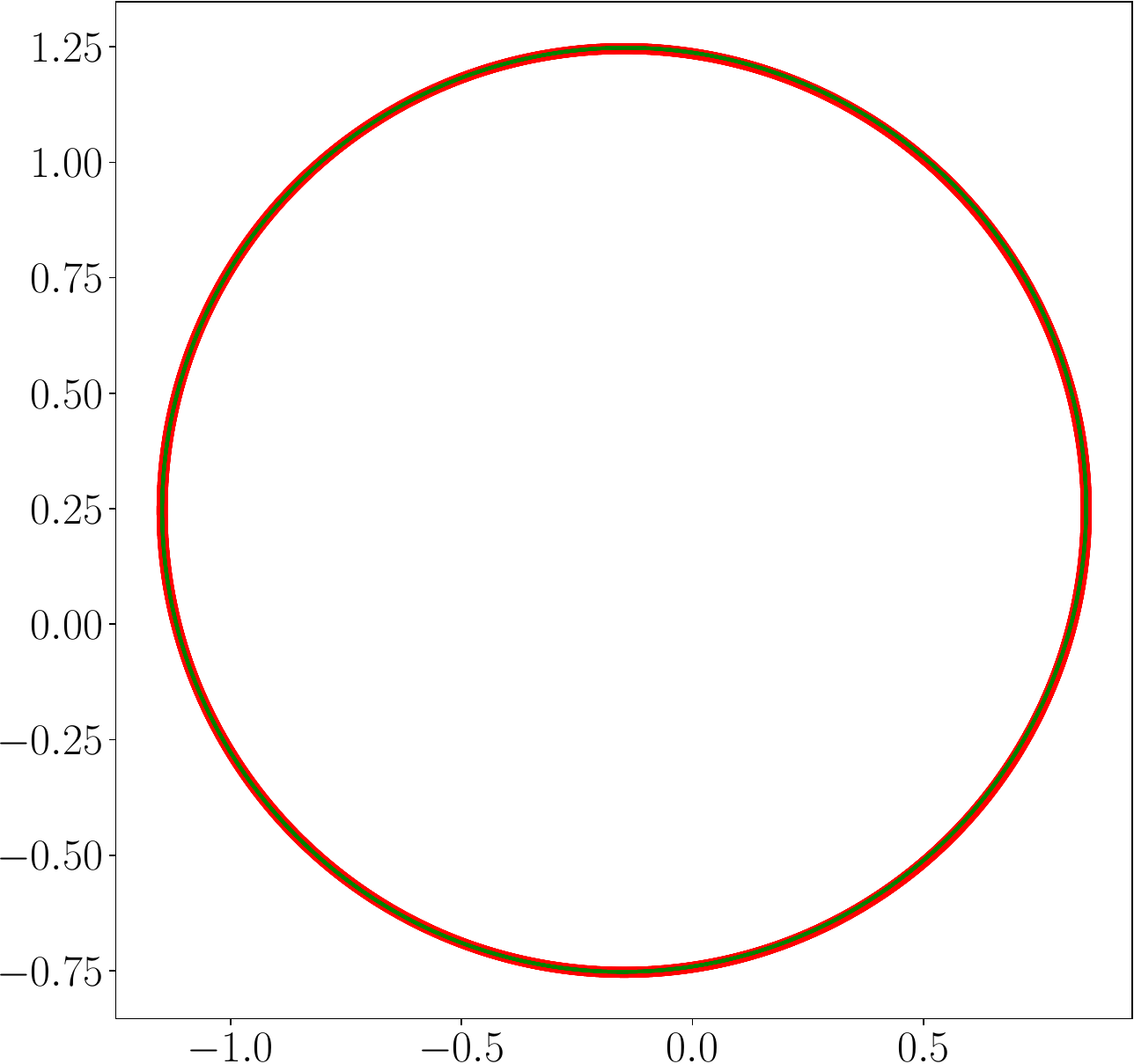}};

        \node[xshift=0.34\textwidth] (top_right) at (top_left.east)
        {\includegraphics[width=0.45\textwidth]{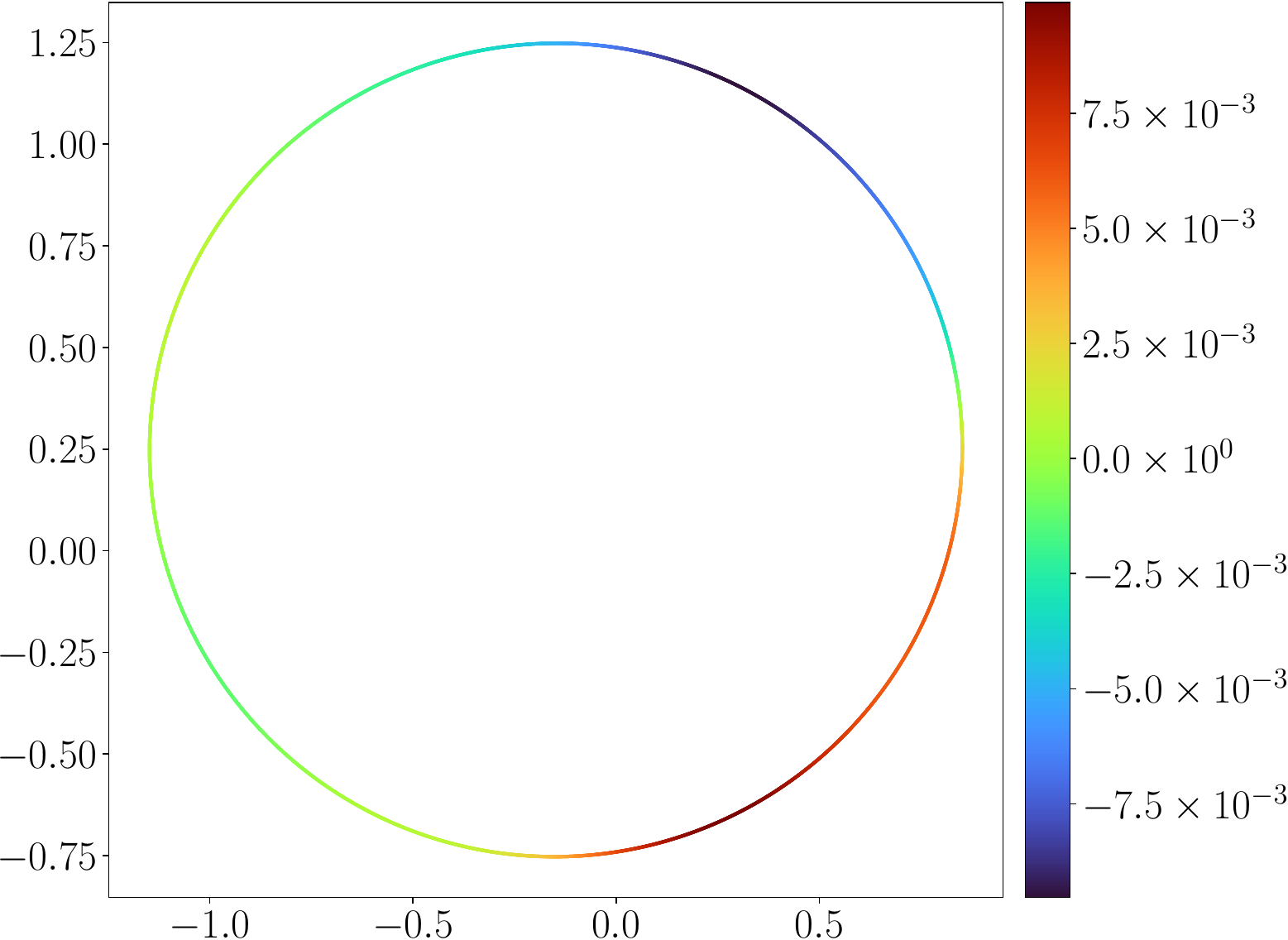}};

        \node[yshift=-0.2\textwidth,xshift=-0.005\textwidth] (bottom_right) at (top_right.south)
        {\includegraphics[width=0.44\textwidth]{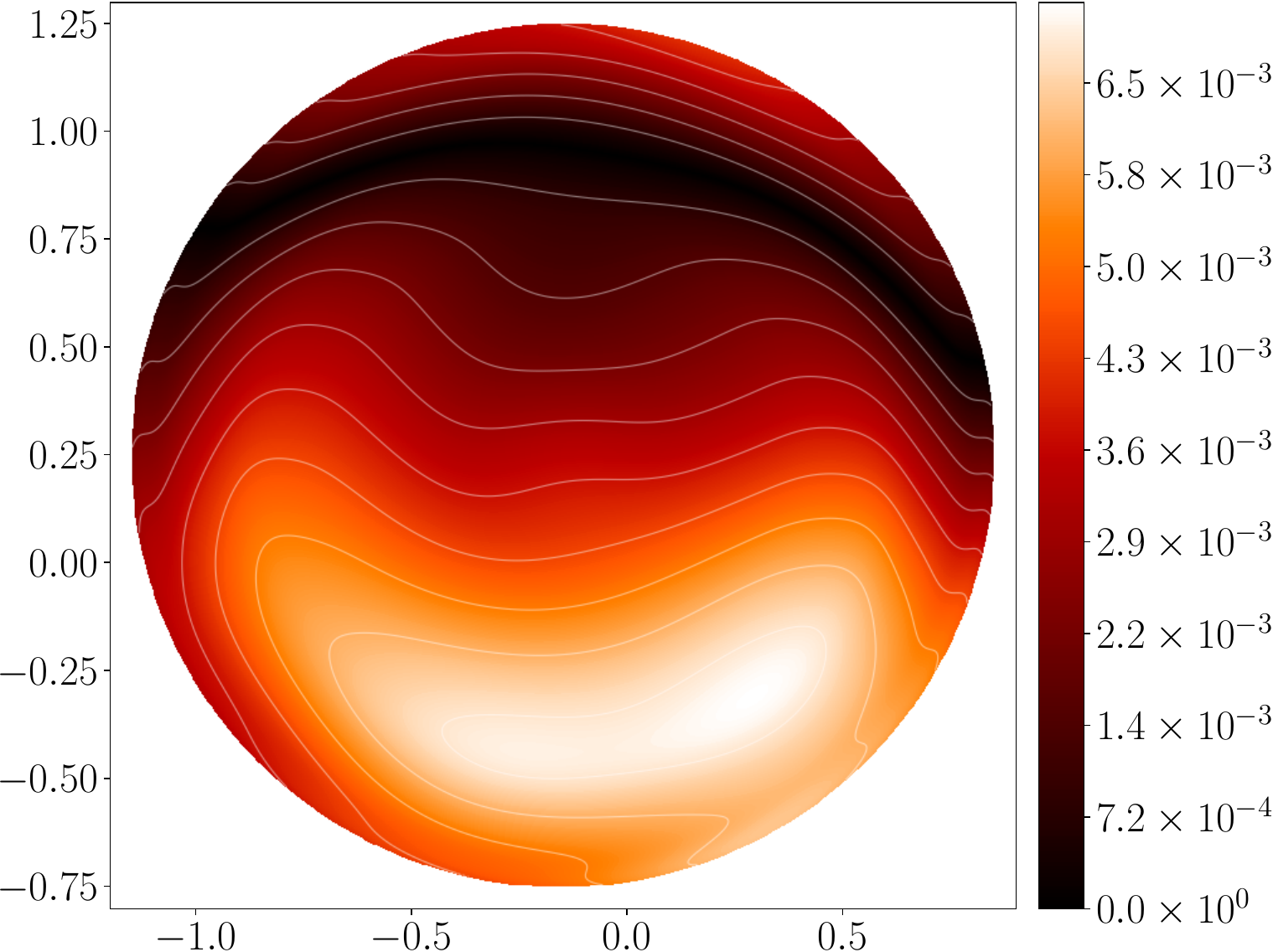}};

        \node[xshift=-0.25\textwidth] (bottom_left) at (bottom_right.west)
        {\includegraphics[width=0.44\textwidth]{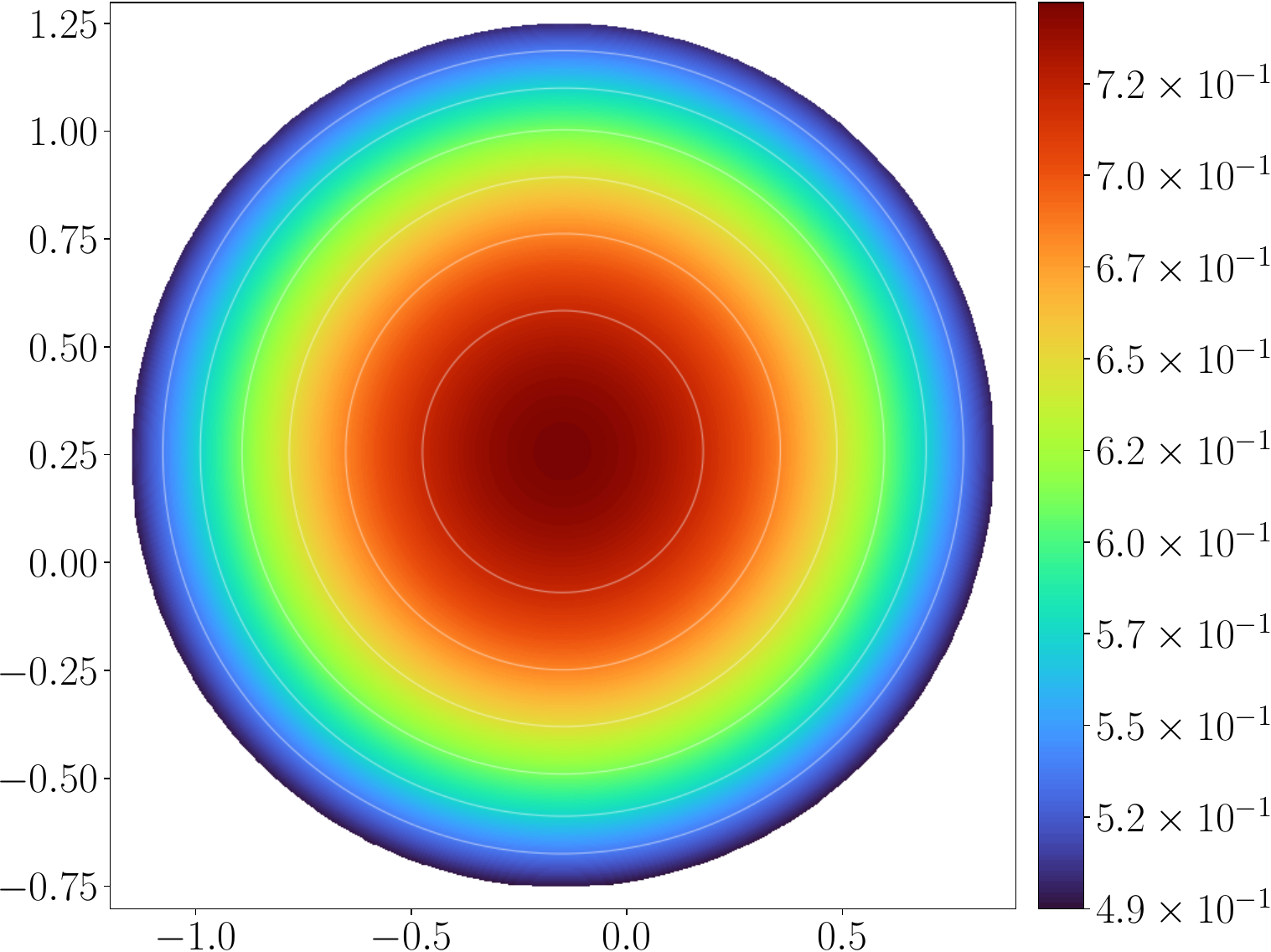}};
    \end{tikzpicture}

    \caption{%
        Shape optimization for the Dirichlet energy,
        for the Robin-Poisson equation with source term $f=1$ and $\kappa=1$.
        Top left panel: learned shape (green line)
        and reference shape (red line).
        Top right panel: deviation from the average of the optimality condition in the learned shape.
        Bottom left panel: approximate PDE solution.
        Bottom right panel: pointwise error
        between the true and approximate PDE solutions.
    }
    \label{fig:shapo_f_1_robin}
\end{figure}

The results with $f=1$ are displayed on \cref{fig:shapo_f_1_robin}.
In this case, the optimal shape is known to be the unit sphere and the solution of the Robin-Poisson problem \eqref{eq:poisson_robin} reads \smash{$u(x,y) = \frac{3-x^2-y^2}{4}$}, see \cite{bucur2015saint}.
Remark that, since the source term is constant, the optimal shape is not necessary centered at the origin, like in \cref{sec:shapo_f_1}.
To get an idea of the quality of our results presented in \cref{fig:shapo_f_1_robin}, we looked at the distance between the learned and the optimal shape (top left panel), as well as the deviation from the average of the optimality condition (top right panel), and the pointwise error between the exact solution and the approximated solution (bottom right panel). Of course, the approximate solution to the Robin-Poisson problem \eqref{eq:poisson_robin} is also displayed (bottom left panel).

As with Dirichlet boundary conditions, we report some metrics
in \cref{tab:stats_DeepShape_constant_robin}.
Remark that the Hausdorff distance, the standard deviation of the optimality condition and the $L^2$ error reported in \cref{tab:stats_DeepShape_constant_robin} are of the same order as the one obtained in \cref{tab:stats_DeepShape_constant}. This is a good indicator of the robustness of our approach with respect to the boundary conditions

\begin{table}[!ht]
    \centering
    \caption{%
        Relevant metrics related to the shape optimization
        of the Robin-Poisson equation with $f=1$ and $\kappa=1$.%
    }
    \label{tab:stats_DeepShape_constant_robin}
    \begin{tabular}{ccc}
        \toprule
        Hausdorff distance
         & Optimality error
         & $L^2$ error           \\
        \cmidrule(lr){1-3}
        $9.86 \times 10^{-3}$
         & $4.99 \times 10^{-3}$
         & $4.13 \times 10^{-3}$
        \\
        \bottomrule
    \end{tabular}
\end{table}

\subsubsection{
    \texorpdfstring
    {Non-parametric problem: results with $f$ given by \eqref{eq:parametric_source} with $\mu=0.8$ and with $\kappa=1$}
    {Non-parametric problem: results with non-constant f and fixed kappa}
}
\label{sec:non_parametric_robin_exp}

\begin{figure}[!ht]
    \centering
    \begin{tikzpicture}
        \node[xshift=-0.25\textwidth] (top_left) at (0,0)
        {\includegraphics[height=0.23\textheight]{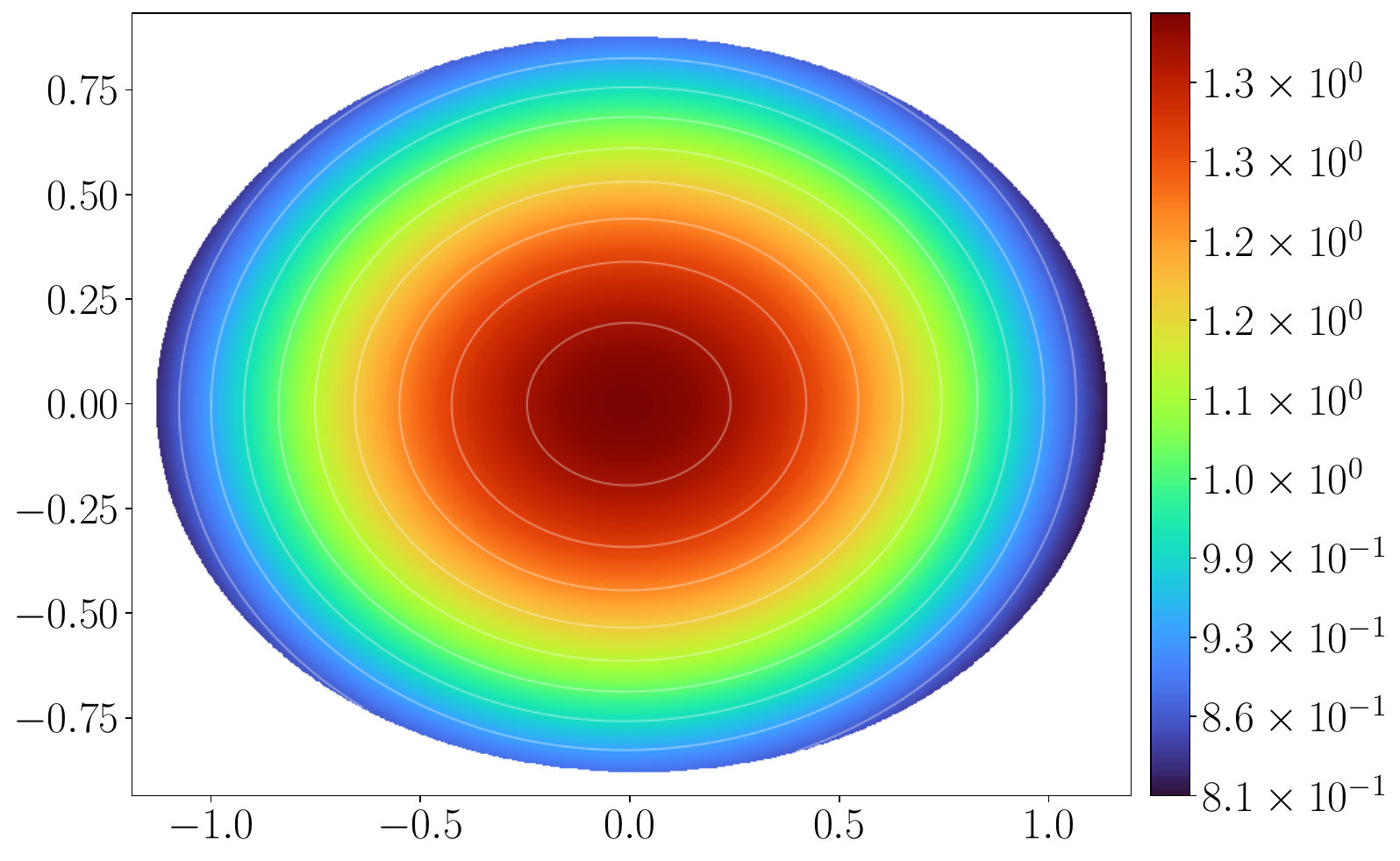}};

        \node[xshift=0.25\textwidth] (top_right) at (top_left.east)
        {\includegraphics[height=0.23\textheight]{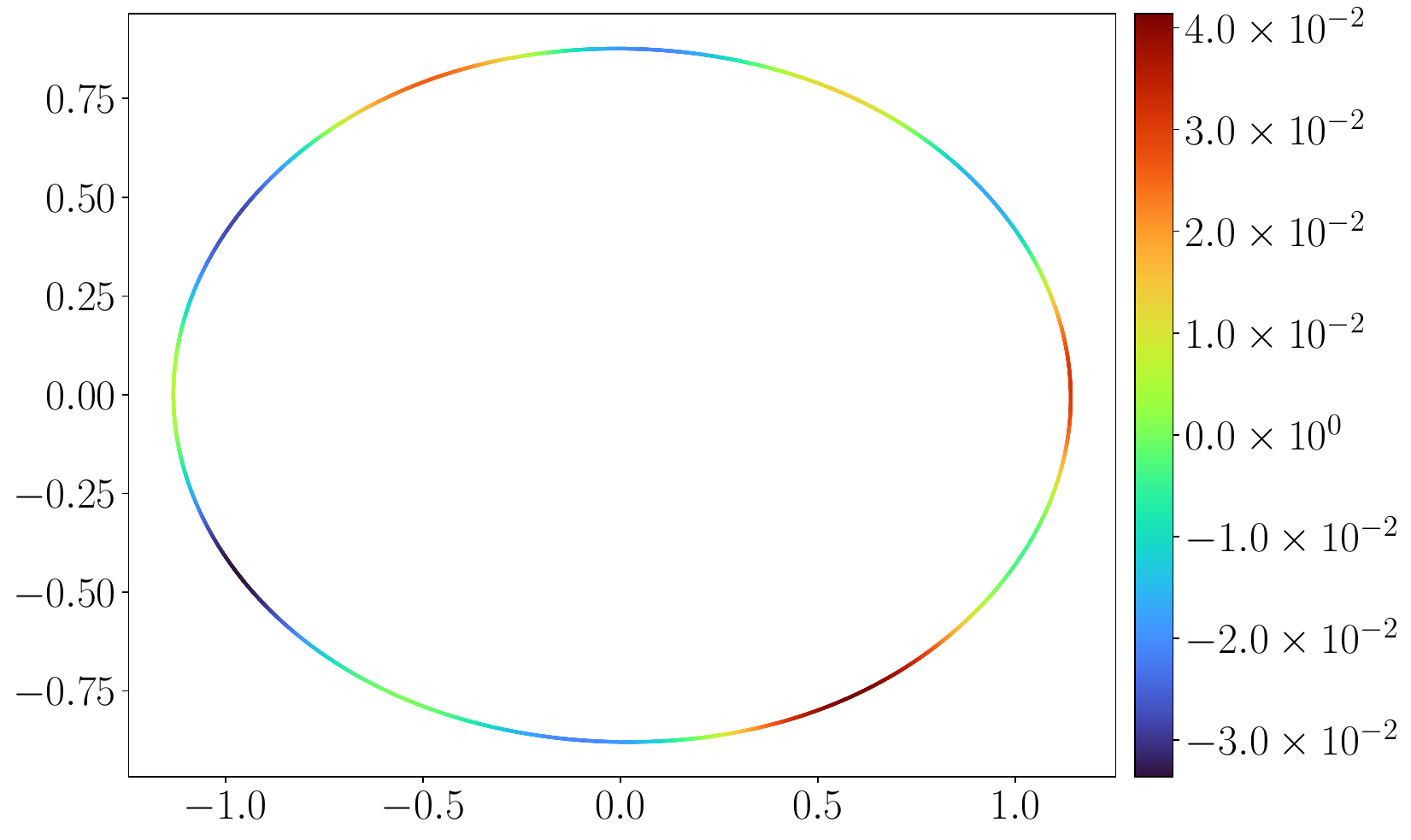}};

    \end{tikzpicture}

    \caption{%
        Shape optimization for the Dirichlet energy,
        for Poisson's equation with source term $f$
        given by~\eqref{eq:parametric_source} with $\mu=0.8$ and $\kappa=1$ and with Robin boundary conditions.
        Left panel: approximate PDE solution in the learned shape.
        Right panel: deviation from the average of the optimality condition.
    }
    \label{fig:shapo_f_exp_Robin}
\end{figure}

To increase the difficulty of the numerical resolution, we turn to an exponential source term $f$ given by~\eqref{eq:parametric_source} with $\mu=0.8$,
and display the results on \cref{fig:shapo_f_exp_Robin}.
Because of the lack of knowledge we have on the potential optimal solution, the only relevant metric we can look at is the optimality condition \eqref{eq:optimality_condition_robin}. For this numerical experiment the optimality error is $1.87 \times 10^{-2}$, which is a satisfying result.

\subsubsection{
    \texorpdfstring
    {Parametric problem: results with $f=1$ and $\kappa \in (0.5, 1.5)$}
    {Parametric problem: results with fixed f and varying kappa}
}
\label{sec:parametric_robin_one}

\begin{figure}[!ht]
    \centering
    \begin{tikzpicture}
        \node[xshift=-0.25\textwidth] (top_left) at (0,0)
        {\includegraphics[width=0.45\textwidth]{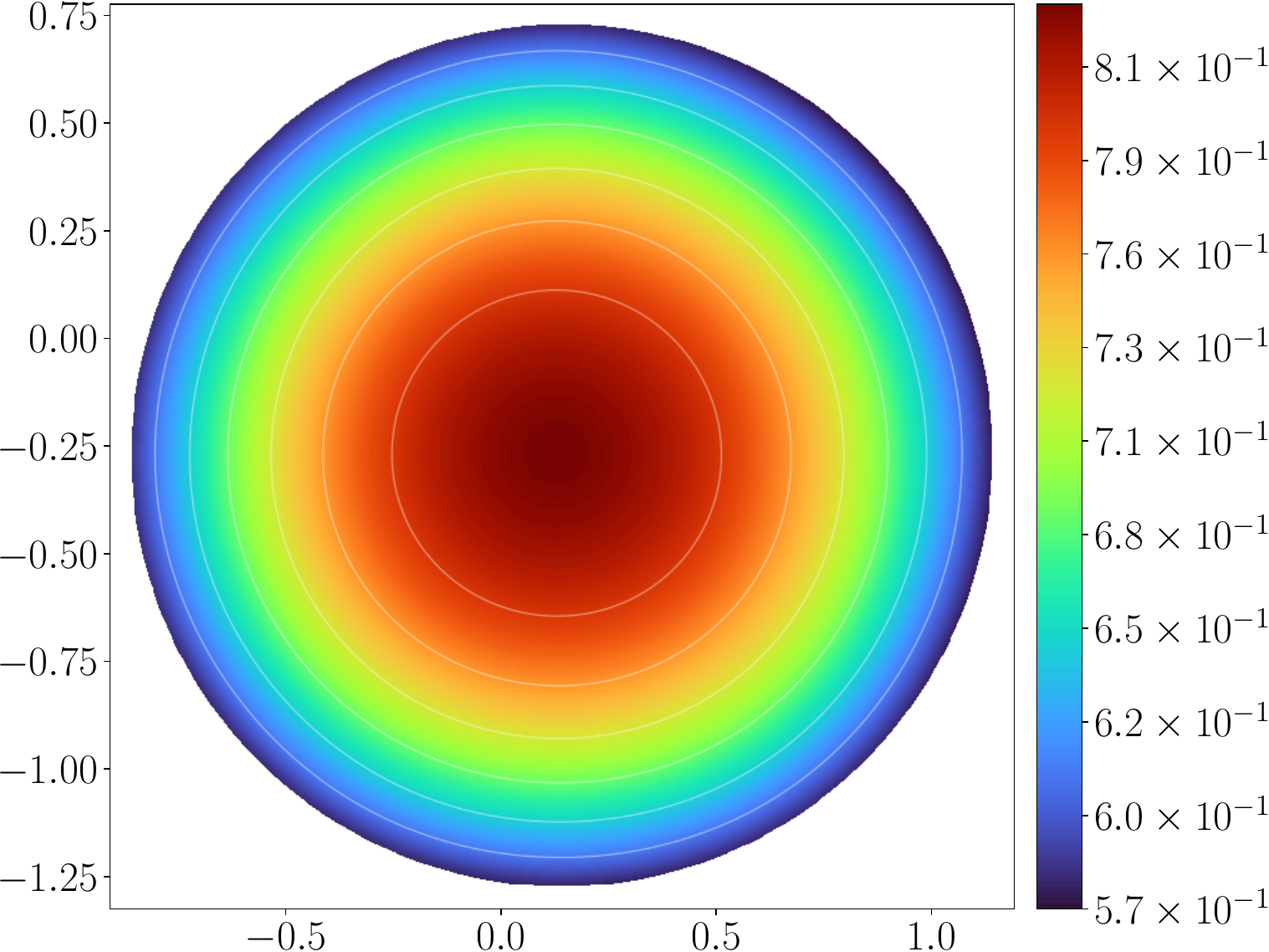}};
        \node[yshift=-0.02\textwidth] at (top_left.south) {(a) solution, $\kappa=0.87$};

        \node[xshift=0.25\textwidth] (top_right) at (top_left.east)
        {\includegraphics[width=0.45\textwidth]{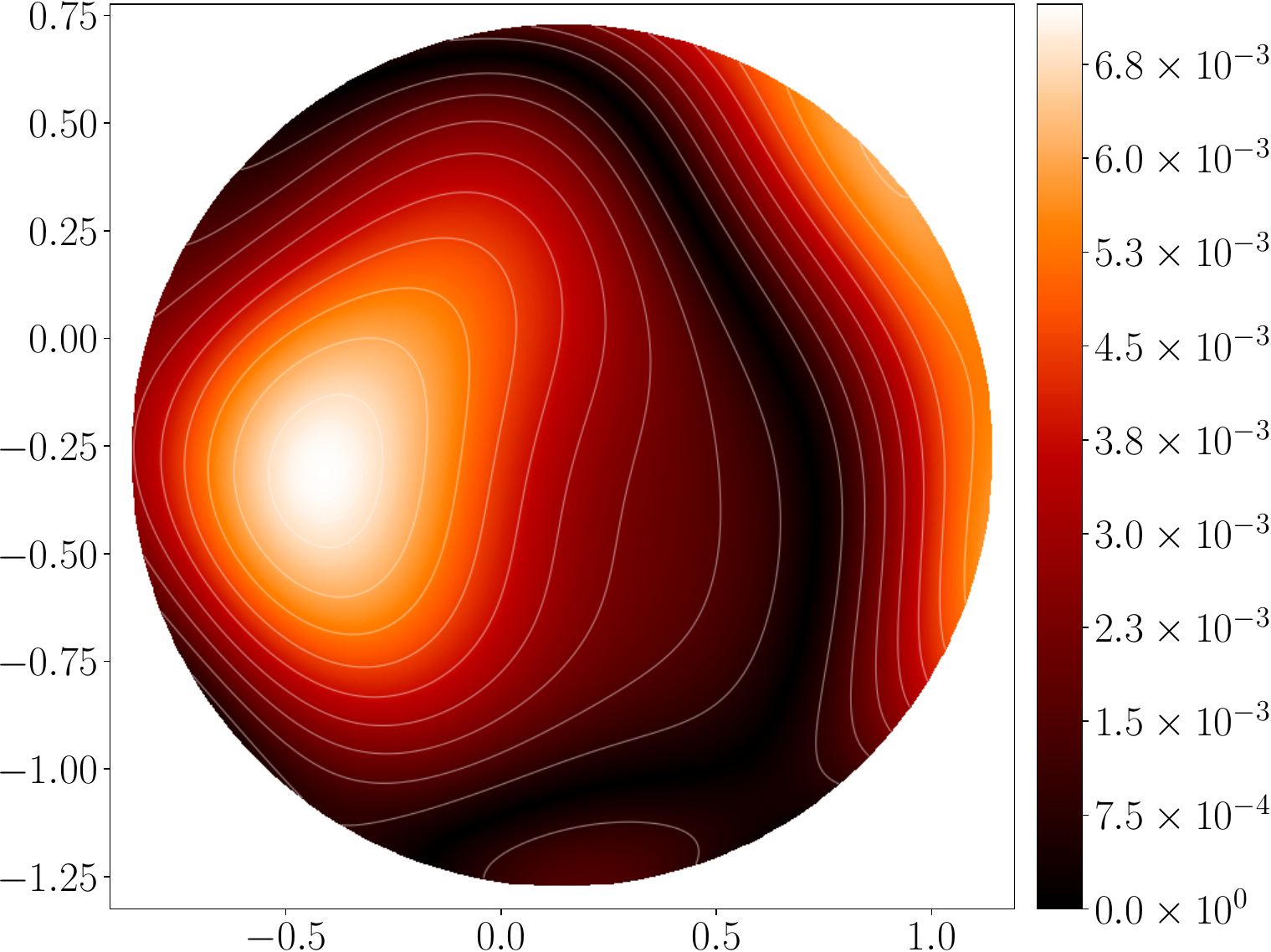}};
        \node[yshift=-0.02\textwidth] at (top_right.south) {(b) error, $\kappa=0.87$};

        \node[yshift=-0.23\textwidth] (bottom_left) at (top_left.south)
        {\includegraphics[width=0.45\textwidth]{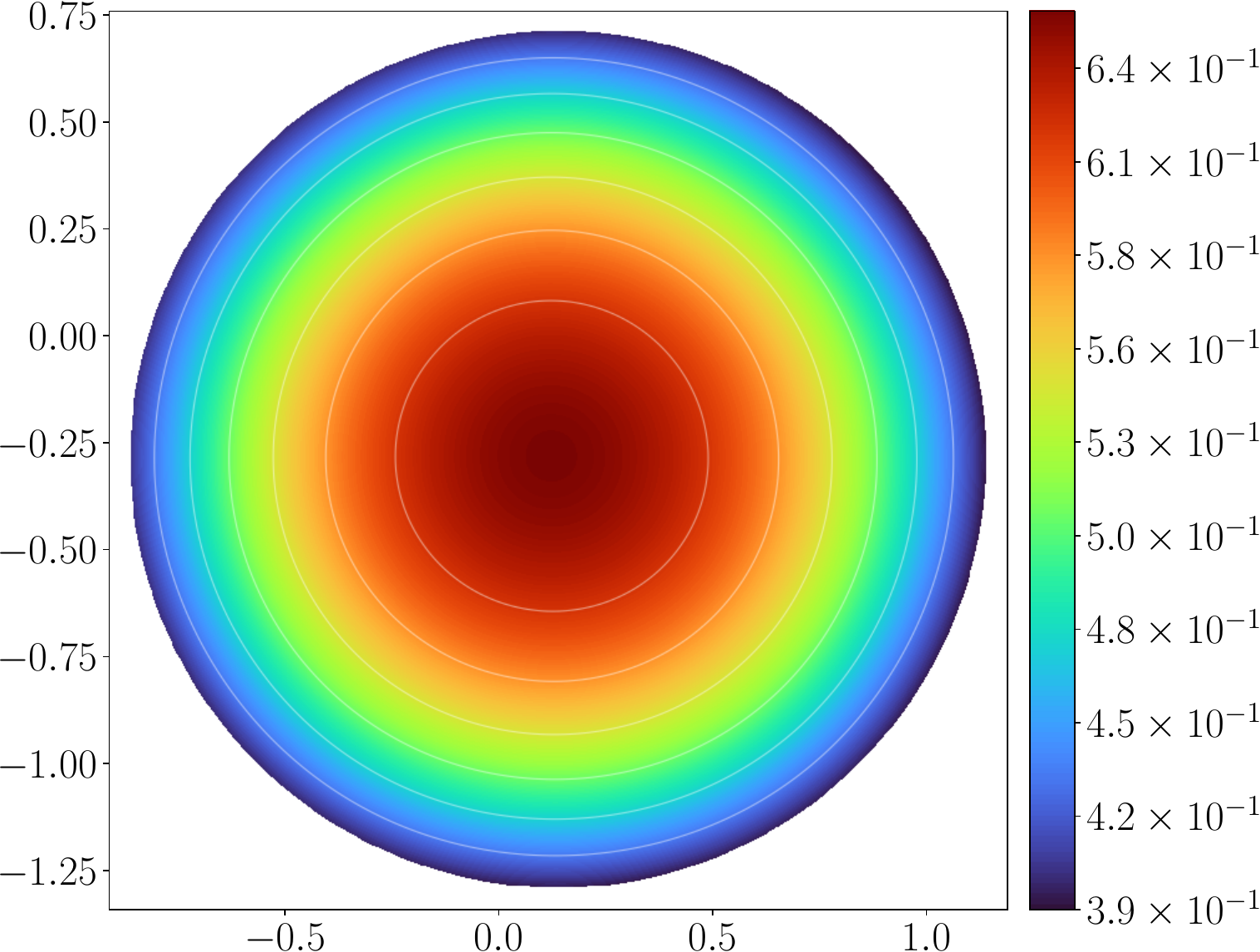}};
        \node[yshift=-0.02\textwidth] at (bottom_left.south) {(c) solution, $\kappa=1.24$};

        \node[yshift=-0.23\textwidth] (bottom_right) at (top_right.south)
        {\includegraphics[width=0.45\textwidth]{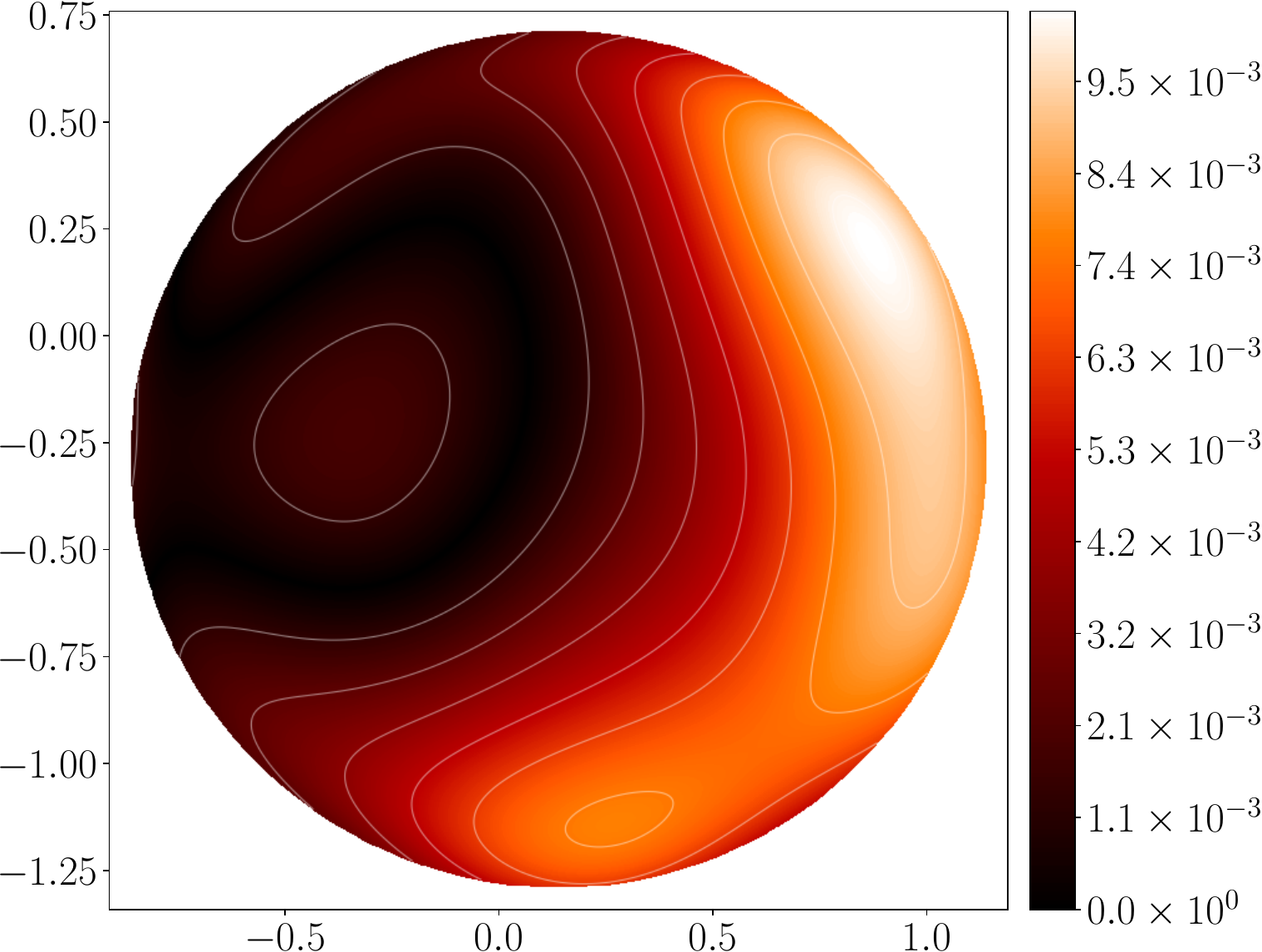}};
        \node[yshift=-0.02\textwidth] at (bottom_right.south) {(d) error, $\kappa=1.24$};

    \end{tikzpicture}

    \caption{%
        Shape optimization for the Dirichlet energy,
        for the Robin-Poisson equation with source term $f=1$ and $\kappa\in(0.5,1.5)$.
        For two values of $\kappa$, we display the approximate solution (left column)
        and the pointwise error between the approximate and exact solutions (right column).
    }
    \label{fig:shapo_f_1_param_robin}
\end{figure}

\begin{figure}[!ht]
    \centering

    \includegraphics[width=0.45\textwidth]{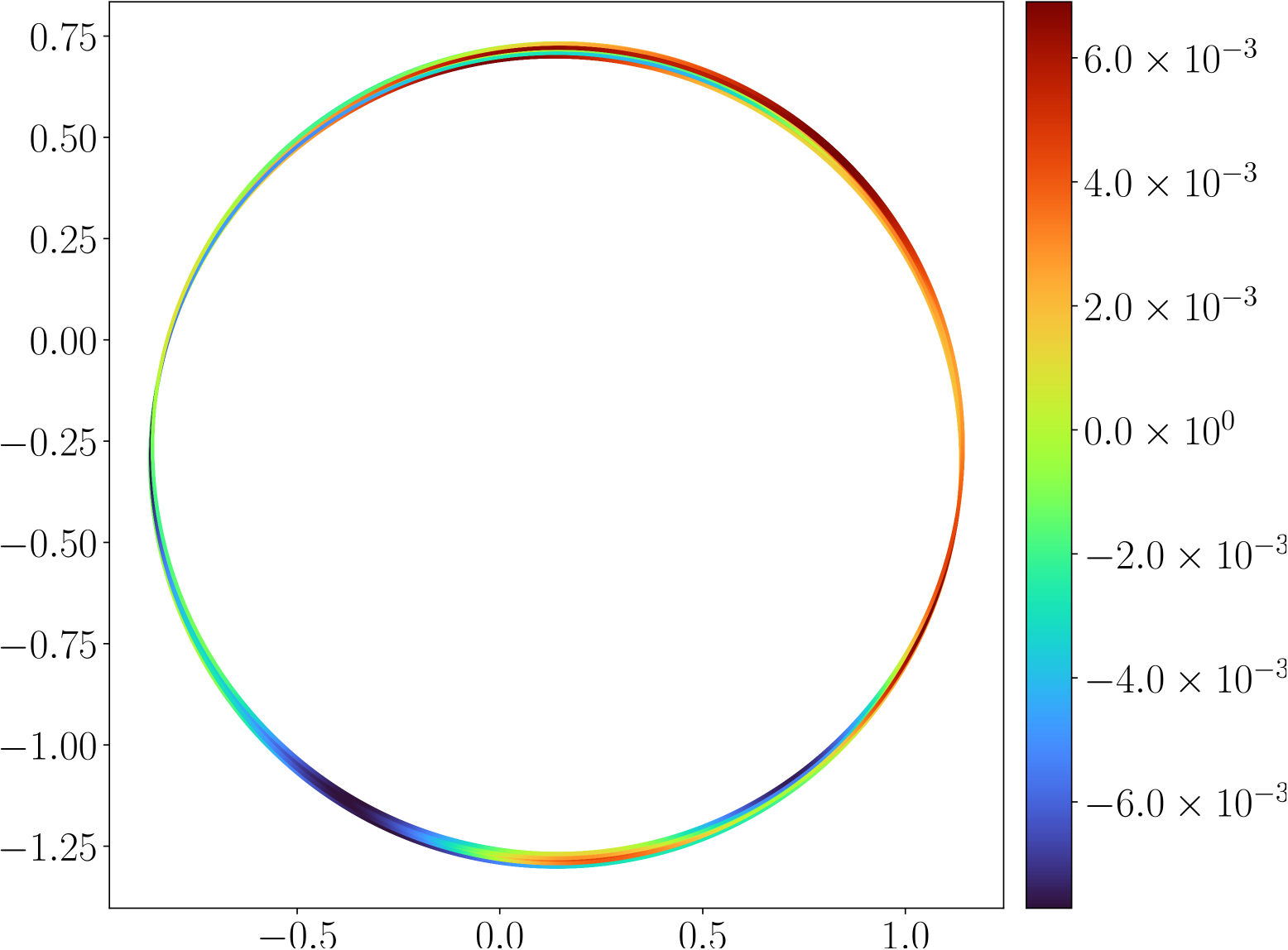}

    \caption{%
        Shape optimization for the Dirichlet energy,
        for the Robin-Poisson equation with source term $f=1$ and $\kappa\in[0.5,1.5]$:
        optimal shapes and deviation from the average of the optimality condition for $10$ random values of~$\kappa$.
    }
    \label{fig:shapo_f_1_param_optimality_conditions_robin}
\end{figure}

Lastly, we tackle a parametric shape optimization problem for $f=1$ and $\kappa \in \mathbb{M} = (0.5, 1.5)$.
As described in~\cite{bucur2015saint}, the optimal shape is a disk, and thus the solution of the Robin-Poisson problem is \smash{$u(x,y) = \frac{2+\kappa}{4\kappa}-\frac{x^2+y^2}{4}$}.
The results and errors are displayed on \cref{fig:shapo_f_1_param_robin},
where there is a good agreement between approximate and exact solution.
In addition, we display the optimal shapes and optimality error for
$10$ random values of $\kappa$ on \cref{fig:shapo_f_1_param_optimality_conditions_robin}.
Finally, we report some statistics on our three main metrics
in \cref{tab:stats_DeepShape_constant_param_robin},
which confirm the relevance of our approach for a parametric problem.

\begin{table}[!ht]
    \centering
    \caption{Statistics of relevant metrics in the case of a constant source term $f=1$ and a parametric constant Robin coefficient $\kappa$, obtained by computing each metric for $10^3$ values of $\kappa$.}
    \label{tab:stats_DeepShape_constant_param_robin}
    \begin{tabular}{ccccc}
        \toprule
        Metric
         & Mean
         & Max
         & Min
         & Standard deviation    \\
        \cmidrule(lr){1-5}
        Hausdorff distance
         & $5.96 \times 10^{-3}$
         & $1.08 \times 10^{-2}$
         & $2.11 \times 10^{-3}$
         & $2.87 \times 10^{-3}$ \\
        optimality error
         & $4.22 \times 10^{-3}$
         & $8.81 \times 10^{-3}$
         & $2.71 \times 10^{-3}$
         & $3.50 \times 10^{-3}$ \\
        $L^2$ error
         & $4.20 \times 10^{-3}$
         & $1.05 \times 10^{-2}$
         & $3.27 \times 10^{-3}$
         & $8.08 \times 10^{-4}$ \\
        \bottomrule
    \end{tabular}
\end{table}

\section{Conclusion}
In this work, we developed a flexible Proof of Concept for geometric shape optimization using neural networks. We minimized the Dirichlet energy under volume constraints by approximating the solution to Poisson's equation and shape transformations with neural networks. Our combined optimization algorithm successfully handled Dirichlet and Robin boundary conditions, parametric right-hand sides or Robin coefficients, and Bernoulli-type free boundary problems. Moreover, it is very easily parallelizable.

A significant advantage of our methodology is that it does not require the calculation of the shape derivative. Additionally, our approach has the potential to address non-self-adjoint cases and could be adapted to handle non-energetic criteria. This flexibility makes our method a promising tool for various complex shape optimization problems.

This work should be regarded as an initial step,
with many questions still open for future exploration.
One significant challenge is extending our approach to
PDEs in higher space dimensions.
While our current work addresses a two-dimensional case,
our goal is to develop methods capable of handling problems
in dimensions larger than three.
    {\ra%
        Indeed, testing our method with
        other neural network architectures,
        such as Hénon networks~\cite{BurTanMau2020}
        or volume preserving networks\cite{ZhuZhuZhaTanLiu2022},
        or exploring relevant applications like eigenvalues of the Laplacian
        or the Stokes problem,
        or even extending to 3D with other neural network architectures,
        are all interesting directions. %
    }%
This challenge is intricately linked to incorporating volume constraints
(using e.g. volume-preserving neural networks \cite{ZhuZhuZhaTanLiu2022}),
or more broadly, managing general manufacturing constraints in any dimension.
A simplistic approach might involve penalizing constraints,
but this often conflicts with learning algorithms,
as it can lead to poorly conditioned systems.

A particularly ambitious direction for future research is to tackle physical problems that traditional shape optimization methods find challenging. For instance, optimizing criteria related to solutions of turbulent fluid dynamics equations is a complex area that warrants further investigation. This will be a key focus of our future efforts.

\section*{Data Availability Statement}
No datasets were generated or analysed during the current study.
\section*{Conflict of interest}
None of the authors have a conflict of interest to disclose.
\section*{Replication of results}
The source code for solving the shape optimization problem is open and available at the following link:
\begin{center}
    \url{https://github.com/belieresfrendo/GeSONN}
\end{center}
A comprehensive documentation explaining how to use it is also provided.

\section*{Acknowledgments}
We would like to sincerely thank Caroline Vernier and Killian Lutz for helpful conversations and their comments on the problem.
We would like to thank Konrad Janik as well, for enlightening discussions on parametric SympNets.

This work was supported by PEPR PDE-AI and by PEPR Numpex Exa-MA.
The last author were partially supported by the ANR Project “STOIQUES”.

%
%

\appendix

\section{
  \texorpdfstring
  {Necessary first order optimality condition for Problem~\eqref{eq:optim_ener_robin}}{Necessary first order optimality condition for the Robin problem}
 }
\label{sec:proof_optimality_condition_robin}

Let us compute the shape derivative of the functional $\mathscr{E}$ defined by
\begin{equation*}
    \label{shape_functional_naive}
    \mathscr{E}(\Upomega) = \frac{1}{2}\int_\Upomega|\nabla u|^2 + \frac{\kappa}{2}\int_{\partial \Upomega}u^2 - \int_\Upomega fu,
\end{equation*}
with $u$ solution of the Robin-Poisson problem \eqref{eq:poisson_robin}.
In what follows, the notation $\nabla_\Gamma$ stands for the tangential gradient.

Multiplying the main equation of \eqref{eq:poisson_robin} by $u$ and integrating then by parts leads to
\begin{equation*}
    \int_\Upomega|\nabla u|^2 + \kappa\int_{\partial \Upomega}u^2 = \int_\Upomega fu,
\end{equation*}
so that
\begin{equation*}
    \label{shape_functional}
    \mathscr{E}(\Upomega) = -\frac{1}{2}\int_\Upomega fu.
\end{equation*}
Let $V\in W^{2,\infty}(\mathbb{R}^2,\mathbb{R}^2)$ denote a compactly supported vector field. In what follows, we extend the normal vector $n$ to a neighborhood of $\partial \Upomega$, in such a way that $n$ remains unitary in a neighborhood of $\partial \Upomega$. This way, all the upcoming calculations will make sense and, in particular, will allow us to define the quantity $\nabla_\Gamma(V\cdot n)$.
Therefore, $\lVert n\rVert^2=1$ in this neighborhood which implies in particular that $(Dn)n=0$.

Then, according to \cite[\rb Theorems~5.2.2 and~5.3.1, Chapter~5]{henrot-pierre},
\begin{equation}
    \label{shape_derivative_1}
    \left\langle\mathrm{d}\mathscr{E}(\Upomega),V\right\rangle =-\frac12\left( \int_\Upomega fu' + \int_{\partial\Upomega}fu (V\cdot n)\right),
\end{equation}
where the Eulerian derivative $u'$ of $u$ solves
\begin{equation}
    \label{eq:poisson_robin_u_prime}
    \begin{cases}
        -\Updelta u' = 0                                                                                                                                                                         & \mathrm{in} \, \Upomega;         \\
        \frac{\partial u'}{\partial n} + \kappa u' = + \nabla u \cdot \nabla_\Gamma(V\cdot n) - \left(\frac{\partial^2u}{\partial n^2} + \kappa \frac{\partial u}{\partial n}\right) (V \cdot n) & \mathrm{on} \, \partial\Upomega.
    \end{cases}
\end{equation}
We now multiply the equation \eqref{eq:poisson_robin} solved by $u$ by $u'$,
and we integrate by parts, to obtain, arguing Green's formula:
\begin{equation}
    \label{fv_u_u_prime}
    \int_\Upomega\nabla u \cdot \nabla u' + \kappa \int_{\partial \Upomega}uu' = \int_\Upomega fu'.
\end{equation}
Symmetrically, by reversing the roles of $u$ and $u'$,
and arguing \eqref{eq:poisson_robin_u_prime}, we get
\begin{equation}
    \label{fv_u_prime_u}
    \int_\Upomega\nabla u' \cdot \nabla u + \kappa \int_{\partial \Upomega}u'u - \int_{\partial \Upomega}u\nabla u \cdot \nabla_\Gamma(V\cdot n) + \int_{\partial \Upomega}u\left(\frac{\partial^2u}{\partial n^2} + \kappa \frac{\partial u}{\partial n}\right) (V \cdot n) = 0.
\end{equation}
By combining \eqref{fv_u_u_prime} and \eqref{fv_u_prime_u}, we obtain
\begin{equation}
    \label{fu_prime}
    \int_\Upomega fu' =  + \int_{\partial \Upomega}u\nabla u \cdot \nabla_\Gamma(V\cdot n) - \int_{\partial \Upomega}u\left(\frac{\partial^2u}{\partial n^2} + \kappa \frac{\partial u}{\partial n}\right) (V \cdot n).
\end{equation}
By simplifying $\int_{\partial \Upomega}u\nabla u \cdot \nabla_\Gamma(V\cdot n)$ with an integration by parts on $\partial \Upomega$ and using an orthogonal decomposition of the gradient, we obtain
\begin{equation*}
    \int_{\partial \Upomega}u\nabla u \cdot \nabla_\Gamma(V\cdot n) = \int_{\partial \Upomega}\frac12 \nabla_\Gamma \left(u^2\right) \cdot \nabla_\Gamma(V\cdot n).
\end{equation*}
Since $\partial \Upomega$ is assumed to be $C^2$, since $V\cdot n \in H^2(\Upomega)$ and $u\in H^3(\Upomega)$, we have \cite[Section 5.4, eq. 5.59]{henrot-pierre}
\begin{equation*}
    \int_{\partial \Upomega}u\nabla u \cdot \nabla_\Gamma(V\cdot n) = \int_{\partial \Upomega} (V\cdot n)\left(-\frac12\Updelta_\Gamma \left(u^2\right)\right),
\end{equation*}
where $\Updelta_\Gamma$ denotes the Laplace-Beltrami operator. Using
\begin{equation*}
    \frac12 \Updelta_\Gamma (u^2) = u\Updelta_\Gamma u + |\nabla_\Gamma u|^2,
\end{equation*}
and the orthogonal decomposition of the Laplacian
\begin{equation*}
    \Updelta_\Gamma u = \Updelta u - H\frac{\partial u}{\partial n} - \frac{\partial^2u}{\partial n^2} = - f - H\frac{\partial u}{\partial n} - \frac{\partial^2u}{\partial n^2},
\end{equation*}
where $H=\mathrm{div}_\Gamma (n)$ denotes the mean curvature of $\partial \Upomega$, and $\mathrm{div}_\Gamma$ the tangential divergence, we come back to \eqref{fu_prime} and infer that
\begin{equation}
    \label{fv_u_prime_2}
    \int_\Upomega fu' =  \int_{\partial \Upomega}(V\cdot n)
    \left(
    u\left(f + H\frac{\partial u}{\partial n} + \frac{\partial^2u}{\partial n^2} \right)
    -|\nabla_\Gamma u|^2
    -u\left(\frac{\partial^2u}{\partial n^2} + \kappa \frac{\partial u}{\partial n}\right)
    \right).
\end{equation}
Note that
\begin{equation*}
    \frac{\partial u}{\partial n}=-\kappa u \qquad \text{a.e. on }\partial \Upomega,
\end{equation*}
so that
\begin{equation*}
    |\nabla_\Gamma u|^2 = |\nabla u|^2 - \left|\frac{\partial u}{\partial n}\right|^2 = |\nabla u|^2 - \kappa^2 u^2.
\end{equation*}
Now, \eqref{fv_u_prime_2} rewrites
\begin{equation*}
    \label{fv_u_prime_3}
    \int_\Upomega fu' =  \int_{\partial \Upomega}(V\cdot n)
    \left(fu - u^2(\kappa H - 2\kappa^2) - |\nabla u|^2
    \right).
\end{equation*}
Coming back to \eqref{shape_derivative_1}, the above computations yield
\begin{equation*}
    \label{shape_derivative_2}
    \left\langle\mathrm{d}\mathscr{E}(\Upomega),V\right\rangle = \int_{\partial\Upomega}(V\cdot n)
    \left(\frac12 |\nabla u|^2 + \frac{u^2}2 (\kappa H - 2\kappa^2) - fu
    \right).
\end{equation*}
Thus, because of the volume constraint, there exists a Lagrange multiplier $c\in \mathbb{R}$ such that
\begin{equation*}
    \frac12 |\nabla u|^2 + \frac{u^2}2(\kappa H - 2\kappa^2) - fu = c,\quad \text{on }\partial \Upomega,
\end{equation*}
we refer for instance to \cite[Section~6.1.3]{henrot-pierre}.

\section{Network hyperparameters}
\label{appendix:hyperparameters}
This appendix provides the hyperparameters used for the networks  in the numerical experiments.
In \cref{tab:hyperparameters},
$\eta_\text{P}$ and $\eta_\text{S}$ are the learning rates
of the PINN and SympNet, respectively,
$\sigma_\text{P}$ and $\sigma_\text{S}$ are their respective activation functions,
$N$ is the number of collocation points,
and $q$ and $\delta$ are respectively the number of layers and width of the SympNet.

\begin{table}[!ht]
    \footnotesize
    \centering
    \label{tab:hyperparameters}
    \caption{%
        Hyperparameters of the networks used in \cref{sec:numerics,sec:robin}.%
    }
    \begin{tabular}{rccccccccc}
        \toprule
        Figure
         & layers                   & $\eta_\text{P}$    & $\sigma_\text{P}$
         & $q$                      & $\delta$           & $\eta_\text{S}$   & $\sigma_\text{S}$
         & $N$                      & epochs                                                                                                          \\
        \cmidrule(lr){1-10}
        \cref{fig:PINN_approximation}
         & [20, 40, 40, 20]         & \num{1e-3}         & $\tanh$
         & ---                      & ---                & ---               & ---
         & \num{20000}              & \num{1500}                                                                                                      \\
        \cref{fig:sympnet_approximation}
         & ---                      & ---                & ---               & 8                 & 10 & \num{1e-2}
         & $\mathrm{sigmoid}$       & \num{10000}        & \num{1000}                                                                                 \\
        \cref{fig:shapo_f_1}
         & [10, 20, 20, 10]         & \num{5e-3}         & $\tanh$           & 4                 & 4
         & \num{5e-3}               & $\mathrm{sigmoid}$ & \num{5000}
         & \num{750}                                                                                                                                  \\
        \cref{fig:shapo_f_exp}
         & [10, 20, 40, 40, 20, 10] & \num{5e-3}         & $\tanh$           & 8                 & 4  & \num{5e-3} & $\mathrm{sigmoid}$ & \num{5000}
         & \num{9000}                                                                                                                                 \\
        \cref{fig:shapo_f_1_param}
         & [10, 20, 40, 40, 20, 10] & \num{5e-3}         & $\tanh$           & 8                 & 4  & \num{5e-3} & $\tanh$            & \num{10000}
         & \num{1000}                                                                                                                                 \\
        \cref{fig:shapo_f_exp_param}
         & [10, 20, 40, 40, 20, 10] & \num{5e-3}         & $\tanh$           & 8                 & 6  & \num{5e-3} & $\tanh$            & \num{10000}
         & \num{2000}                                                                                                                                 \\
        \cref{fig:shapo_f_bizaroid_param}
         & [10, 20, 40, 40, 20, 10] & \num{5e-3}         & $\tanh$           & 8                 & 6  & \num{5e-3} & $\tanh$            & \num{10000}
         & \num{2000}                                                                                                                                 \\
        \cref{fig:bernoulli}
         & [10, 20, 40, 40, 20, 10] & \num{1e-2}         & $\tanh$           & 8                 & 4  & \num{1e-2} & $\mathrm{sigmoid}$ & \num{2500}
         & \num{15000}                                                                                                                                \\
        \cref{fig:shapo_f_1_robin}
         & [10, 20, 40, 40, 20, 10] & \num{1e-2}         & $\tanh$           & 8                 & 5  & \num{1e-2} & $\tanh$            & \num{10000}
         & \num{10000}                                                                                                                                \\
        \cref{fig:shapo_f_exp_Robin}
         & [10, 20, 40, 40, 20, 10] & \num{1e-2}         & $\tanh$           & 8                 & 5  & \num{1e-2} & $\tanh$            & \num{10000}
         & \num{10000}                                                                                                                                \\
        \cref{fig:shapo_f_1_param_robin}
         & [10, 20, 40, 40, 20, 10] & \num{1e-2}         & $\tanh$           & 8                 & 5  & \num{1e-2} & $\tanh$            & \num{10000}
         & \num{10000}                                                                                                                                \\
        \bottomrule
    \end{tabular}
\end{table}


\begin{thebibliography}{10}

    \bibitem{allaire2007conception}
    G.~Allaire.
    \newblock {\em Conception optimale de structures}, volume~58 of {\em
      Mathématiques et Applications}.
    \newblock Springer Berlin Heidelberg, 2006.
    
    \bibitem{allaire2021shape}
    G.~Allaire, C.~Dapogny, and F.~Jouve.
    \newblock Shape and topology optimization.
    \newblock In {\em Handbook of numerical analysis}, volume~22, pages 1--132.
      Elsevier, 2021.
    
    \bibitem{allaire2002level}
    G.~Allaire, F.~Jouve, and A.-M. Toader.
    \newblock A level-set method for shape optimization.
    \newblock {\em C. R. Math. Acad. Sci. Paris}, 334(12):1125--1130, 2002.
    
    \bibitem{AnsYan2024}
    J.~Ansel, E.~Yang, et~al.
    \newblock {PyTorch 2: Faster Machine Learning Through Dynamic Python Bytecode
      Transformation and Graph Compilation}.
    \newblock In {\em Proceedings of the 29th ASPLOS, Volume 2}, volume~5 of {\em
      ASPLOS '24}, pages 929--947. ACM, 2024.
    
    \bibitem{Arnold}
    V.~I. Arnold.
    \newblock {\em {Mathematical Methods of Classical Mechanics}}.
    \newblock Graduate Texts in Mathematics. Springer New York, 1989.
    
    \bibitem{bucur2015faber}
    D.~Bucur and A.~Giacomini.
    \newblock {Faber--Krahn inequalities for the Robin-Laplacian: A free
      discontinuity approach}.
    \newblock {\em Arch. Ration. Mech. An.}, 218(2):757--824, 2015.
    
    \bibitem{bucur2015saint}
    D.~Bucur and A.~Giacomini.
    \newblock {The Saint-Venant inequality for the Laplace operator with Robin
      boundary conditions}.
    \newblock {\em Milan J. Math.}, 83(2):327--343, 2015.
    
    \bibitem{bucur2016robin}
    D.~Bucur, A.~Giacomini, and P.~Trebeschi.
    \newblock {The Robin--Laplacian problem on varying domains}.
    \newblock {\em Calc. Var. Partial Dif.}, 55:1--29, 2016.
    
    \bibitem{BurTanMau2020}
    J.~W. Burby, Q.~Tang, and R.~Maulik.
    \newblock Fast neural {P}oincaré maps for toroidal magnetic fields.
    \newblock {\em Plasma Phys. Contr. F.}, 63(2):024001, 2020.
    
    \bibitem{amaury_belieres_frendo_2024_13133269}
    A.~Bélières-Frendo and V.~Michel-Dansac.
    \newblock {belieresfrendo/GeSONN: Geometric Shape Optimization with Neural
      Networks}, 2024.
    \newblock GitHub repository.
    
    \bibitem{caflisch1998monte}
    R.~E. Caflisch.
    \newblock {Monte Carlo and quasi-Monte Carlo methods}.
    \newblock {\em Acta Numer.}, 7:1--49, 1998.
    
    \bibitem{chambolle:hal-04140177}
    A.~Chambolle, I.~Mazari-Fouquer, and Y.~Privat.
    \newblock {Stability of optimal shapes and convergence of thresholding
      algorithms in linear and spectral optimal control problems}.
    \newblock preprint, 2023.
    
    \bibitem{ChaSur2020}
    A.~Chandrasekhar and K.~Suresh.
    \newblock {TOuNN: Topology Optimization using Neural Networks}.
    \newblock {\em Struct. Multidiscip. O.}, 63(3):1135--1149, 2020.
    
    \bibitem{CuoGiaIzzNitPicTro2022}
    S.~Cuomo, F.~Giampaolo, S.~Izzo, C.~Nitsch, F.~Piccialli, and C.~Trombetti.
    \newblock A physics-informed learning approach to {B}ernoulli-type free
      boundary problems.
    \newblock {\em Comput. Math. Appl.}, 128:34--43, 2022.
    
    \bibitem{delfour2011shapes}
    M.~C. Delfour and J.-P. Zol{\'e}sio.
    \newblock {\em Shapes and geometries: metrics, analysis, differential calculus,
      and optimization}.
    \newblock SIAM, 2011.
    
    \bibitem{e2017deep}
    W.~E and B.~Yu.
    \newblock {The Deep Ritz Method: A Deep Learning-Based Numerical Algorithm for
      Solving Variational Problems}.
    \newblock {\em Commun. Math. Stat.}, 6(1):1--12, 2018.
    
    \bibitem{franck2023approximately}
    E.~Franck, V.~Michel-Dansac, and L.~Navoret.
    \newblock {Approximately well-balanced Discontinuous Galerkin methods using
      bases enriched with Physics-Informed Neural Networks}.
    \newblock {\em J. Comput. Phys.}, 512:113144, 2024.
    
    \bibitem{Goodfellow-et-al-2016}
    I.~Goodfellow, Y.~Bengio, and A.~Courville.
    \newblock {\em {Deep Learning}}.
    \newblock MIT Press, 2016.
    
    \bibitem{hairer2006structure}
    E.~Hairer, C.~Lubich, and G.~Wanner.
    \newblock {\em {Geometric Numerical Integration}}.
    \newblock Springer Series in Computational Mathematics. Springer-Verlag, 2006.
    
    \bibitem{MR3043640}
    F.~Hecht.
    \newblock {New Development in FreeFem++}.
    \newblock {\em J. Numer. Math.}, 20(3–4), 2012.
    
    \bibitem{antunesshape}
    A.~Henrot.
    \newblock {\em Shape optimization and spectral theory}.
    \newblock De Gruyter Open, 2017.
    
    \bibitem{henrot-pierre}
    A.~Henrot and M.~Pierre.
    \newblock {\em {Shape Variation and Optimization}}.
    \newblock European Mathematical Society Publishing House, 2018.
    
    \bibitem{henrot2000existence}
    A.~Henrot and H.~Shahgholian.
    \newblock Existence of classical solutions to a free boundary problem for the
      $p$-{L}aplace operator: ({I}) the exterior convex case.
    \newblock {\em J. reine angew. Math.}, 2000(521), 2000.
    
    \bibitem{JEONG2023115484}
    H.~Jeong, J.~Bai, C.~P. Batuwatta-Gamage, C.~Rathnayaka, Y.~Zhou, and Y.~Gu.
    \newblock {A Physics-Informed Neural Network-based Topology Optimization
      (PINNTO) framework for structural optimization}.
    \newblock {\em Eng. Struct.}, 278:115484, 2023.
    
    \bibitem{9716789}
    P.~Jin, Z.~Zhang, I.~G. Kevrekidis, and G.~E. Karniadakis.
    \newblock {Learning Poisson Systems and Trajectories of Autonomous Systems via
      Poisson Neural Networks}.
    \newblock {\em IEEE Trans. Neural Netw. Learn. Syst.}, 34(11):8271--8283, 2023.
    
    \bibitem{JIN2020166}
    P.~Jin, Z.~Zhang, A.~Zhu, Y.~Tang, and G.~E. Karniadakis.
    \newblock Symp{N}ets: {I}ntrinsic structure-preserving symplectic networks for
      identifying {H}amiltonian systems.
    \newblock {\em Neural Networks}, 132:166--179, 2020.
    
    \bibitem{kac1966can}
    M.~Kac.
    \newblock Can one hear the shape of a drum?
    \newblock {\em The american mathematical monthly}, 73(4P2):1--23, 1966.
    
    \bibitem{kingma2014adam}
    D.~Kingma and J.~Ba.
    \newblock {Adam: A Method for Stochastic Optimization}.
    \newblock In {\em International Conference on Learning Representations (ICLR)},
      San Diego, CA, USA, 2015.
    
    \bibitem{LagLikFot1998}
    I.~E. Lagaris, A.~Likas, and D.~I. Fotiadis.
    \newblock Artificial neural networks for solving ordinary and partial
      differential equations.
    \newblock {\em IEEE Trans. Neural Netw.}, 9(5):987--1000, 1998.
    
    \bibitem{mohammadi2009applied}
    B.~Mohammadi and O.~Pironneau.
    \newblock {\em Applied shape optimization for fluids}.
    \newblock OUP Oxford, 2009.
    
    \bibitem{nakahara2018geometry}
    M.~Nakahara.
    \newblock {\em Geometry, {T}opology and {P}hysics}.
    \newblock CRC Press, 2018.
    
    \bibitem{odot2023real}
    A.~Odot, G.~Mestdagh, Y.~Privat, and S.~Cotin.
    \newblock Real-time elastic partial shape matching using a neural network-based
      adjoint method.
    \newblock In {\em International Conference on Optimization and Learning}, pages
      137--147. Springer, 2023.
    
    \bibitem{raissi2019physics}
    M.~Raissi, P.~Perdikaris, and G.~E. Karniadakis.
    \newblock Physics-informed neural networks: {A} deep learning framework for
      solving forward and inverse problems involving nonlinear partial differential
      equations.
    \newblock {\em J. Comput. Phys.}, 378:686--707, 2019.
    
    \bibitem{serrin1971symmetry}
    J.~Serrin.
    \newblock A symmetry problem in potential theory.
    \newblock {\em Arch. Ration. Mech. An.}, 43(4):304--318, 1971.
    
    \bibitem{taha2015efficient}
    A.~A. Taha and A.~Hanbury.
    \newblock {An Efficient Algorithm for Calculating the Exact Hausdorff
      Distance}.
    \newblock {\em IEEE Trans. Pattern Anal. Mach. Intell.}, 37(11):2153--2163,
      2015.
    
    \bibitem{tonnel2023avaframe}
    M.~Tonnel, A.~Wirbel, F.~Oesterle, and J.-T. Fischer.
    \newblock {AvaFrame com1DFA (v1.3): a thickness-integrated computational
      avalanche module -- theory, numerics, and testing}.
    \newblock {\em Geosci. Model Dev.}, 16(23):7013--7035, 2023.
    
    \bibitem{WolAagBaeSig2022}
    R.~V. Woldseth, N.~Aage, J.~A. Bærentzen, and O.~Sigmund.
    \newblock On the use of artificial neural networks in topology optimisation.
    \newblock {\em Struct. Multidiscip. O.}, 65(10), 2022.
    
    \bibitem{ZehLiCorTho2021}
    J.~Zehnder, Y.~Li, S.~Coros, and B.~Thomaszewski.
    \newblock {NT}opo: mesh-free topology optimization using implicit neural
      representations.
    \newblock In {\em Proceedings of the 35th International Conference on Neural
      Information Processing Systems}, NIPS '21, Red Hook, NY, USA, 2021. Curran
      Associates Inc.
    
    \bibitem{ZhaYaoLiZhoChe2023}
    Z.~Zhang, W.~Yao, Y.~Li, W.~Zhou, and X.~Chen.
    \newblock Topology optimization via implicit neural representations.
    \newblock {\em Comput. Method. Appl. M.}, 411:116052, 2023.
    
    \bibitem{ZhuZhuZhaTanLiu2022}
    A.~Zhu, B.~Zhu, J.~Zhang, Y.~Tang, and J.~Liu.
    \newblock {VPNets: Volume-preserving neural networks for learning source-free
      dynamics}.
    \newblock {\em J. Comput. Appl. Math.}, 416:114523, 2022.
    
    \end{thebibliography}
\end{document}